\documentclass{m2an}

\usepackage{amsmath,amsfonts,amssymb,mathtools}
\usepackage{graphicx,subcaption}
\usepackage{braket,dsfont,url,enumerate}
\usepackage{tikz}
\usepackage{pdfpages} 
\usepackage{cmap}
\usepackage{pgfplots}
\usepackage{hyperref}
\usepackage{color,dsfont,braket}
\usepackage{bbm}
\usepackage{fancyhdr}
\usepackage{times}
\usepackage{cmap}
\usepackage{pgf}
\usepackage{comment}
\usepackage{algorithm}
\usepackage{algorithmic}
\makeatletter
\@namedef{subjclassname@1991}{2010 Mathematics Subject Classification}
\makeatother
\definecolor{darkred}{rgb}{.7,0,0}

\definecolor{green}{rgb}{0,0.7,0}

\newtheorem{assumption}{Assumption}

\newcommand{\supp}{\operatorname{supp}}

\newcommand{\M}{\mathcal{M}}

\newcommand{\na}{\nabla}
\newcommand{\pa}{\partial}
\newcommand{\eps}{\varepsilon}

\newcommand{\Om}{\Omega}

\newcommand{\IOm}{I\times \Om}

\newcommand{\vertiii}[1]{{\left\vert\kern-0.25ex\left\vert\kern-0.25ex\left\vert #1
    \right\vert\kern-0.25ex\right\vert\kern-0.25ex\right\vert}}

\newcommand{\norm}[1]{\lVert#1\rVert}
\newcommand{\mnorm}[1]{\norm{#1}_{\M(\Omega)}}
\newcommand{\ltwonorm}[1]{\lVert#1\rVert_{L^2(\Omega)}}
\newcommand{\lpnorm}[1]{\lVert#1\rVert_{L^p(\Omega)}}
\newcommand{\linfnorm}[1]{\lVert#1\rVert_{L^\infty(\Omega)}}
\newcommand{\linfnormn}[1]{\lVert#1\rVert_{L^\infty(\Omega_0)}}

\newcommand{\abs}[1]{\lvert#1\rvert}
\newcommand{\Abs}[1]{\left\lvert#1\right\rvert}

\newcommand{\half}{\frac{1}{2}}

\newcommand{\ou}{\bar u}

\newcommand{\oq}{\bar q}
\newcommand{\oz}{\bar z}
\newcommand{\ob}{\bar \beta}
\newcommand{\oxi}{\bar x_i}
\newcommand{\oxki}{\bar x_{k,i}}
\newcommand{\oxkhi}{\bar x_{kh,i}}
\newcommand{\oxkhij}{\bar x_{kh,ij}}
\newcommand{\oukh}{\bar{u}_{kh}}
\newcommand{\oqkh}{\bar{q}_{kh}}

\newcommand{\R}{\mathbb{R}}
\newcommand{\N}{\mathbb{N}}
\newcommand{\lh}{\abs{\ln{h}}}
\newcommand{\lk}{\ln{\frac{T}{k}}}

\newcommand{\Xk}{X_k^r}
\newcommand{\Xkh}{X^{r,1}_{k,h}}
\DeclareMathOperator*{\argmax}{arg\,max}
\newcommand{\tXk}{\widetilde X_k^r}
\newcommand{\hXk}{\widehat X_k^r}

\definecolor{green}{rgb}{0,0.7,0}

\begin{document}
%
\author{Dmitriy Leykekhman}\address{Department of Mathematics,
               University of Connecticut,
              Storrs,
              CT~06269, USA (dmitriy.leykekhman@uconn.edu).}
\author{Boris Vexler}\address{Chair of Optimal Control, Technical University of Munich,
Department of Mathematics , Boltzmannstra{\ss}e 3, 85748 Garching b. Munich, Germany
(vexler@ma.tum.de).}
\author{Daniel Walter}\address{Johann Radon Institute for Computational and Applied Mathematics, {\"O}AW, Altenbergerstra{\ss}e 69, 4040 Linz, Austria (daniel.walter@oeaw.ac.at). The third author gratefully acknowledges support from
	the International Research Training Group IGDK, funded by the
	German Science Foundation (DFG) and the Austrian Science
	Fund (FWF).}
\begin{abstract} In this paper we consider a problem of initial data identification from the final time observation for homogeneous parabolic problems. It is well-known that such problems are exponentially ill-posed due to the  strong smoothing property of parabolic equations.  We are interested in a situation when the initial data we intend to recover is known to be sparse, i.e. its support has Lebesgue measure zero. We formulate the problem as an optimal control problem and incorporate the information on the sparsity of the unknown initial data into the structure of the objective functional. In particular, we are looking for the control variable in the space of regular Borel measures and use the corresponding norm  as a regularization term in the objective functional. This leads to a convex but non-smooth optimization problem. For the discretization we use continuous piecewise linear finite elements in space and discontinuous Galerkin finite elements of arbitrary degree in time. For the general case we establish error estimates for the state variable. Under a certain structural assumption, we show that the control variable consists of a finite linear combination of Dirac measures. For this case we obtain  error estimates for the locations of Dirac measures as well as for the corresponding coefficients.  The key to the numerical analysis are the sharp smoothing type pointwise finite element error estimates for homogeneous parabolic problems, which are of independent interest. Moreover, we discuss an efficient algorithmic approach to the problem and show several numerical experiments illustrating our theoretical results.  \end{abstract}
\title{Numerical Analysis of Sparse Initial Data Identification\\ for Parabolic Problems}
\maketitle

\section{Introduction}
In this paper we consider a problem of identification of an unknown initial data $q$ for a homogenous parabolic equation
\begin{equation}\label{eq:state}
\begin{aligned}
\pa_t u-\Delta u &= 0 &&\text{in} \quad (0,T)\times \Om,\;  \\
u                &= 0 &&\text{on}\quad  (0,T)\times\pa\Omega, \\
u(0)             &= q &&\text{in}\quad \Omega,
\end{aligned}
\end{equation}
from a given (measured) data $u_d \approx u(T)$ of the terminal state $u(T)$ for some $T>0$. In general, this problem is known to be exponentially ill-posed, see, e.g.,~\cite{IsakovV_2017}. We are interested in the situation, where the initial data we are looking for, is known to be sparse, i.e. to have a support of Lebesgue measure zero. The strong smoothing property of the above equation makes it difficult to identify such  sparse initial data. The remedy is the incorporation of the information that the unknown $q$ should be sparse in the optimal control formulation. Following the idea for measure valued formulation of sparse control problems, see, e.g.,~\cite{CasasE_VexlerB_ZuazuaE_2015, CasasZuazua:2013, ClasonKunisch:2011,KunischPieperVexler:2014, PieperVexler:2013}, we will look for the initial state $q$ in the space of regular Borel measures $\M(\Omega)$ on the domain $\Omega$, which is known to be isomorphic to the dual space of continuous functions which are zero on~$\partial \Omega$, $C_0(\Omega)^*$.

The corresponding optimal control formulation reads as follows
\begin{equation}\label{eq:problem}
\text{Minimize }\; J(q,u) = \half \norm{u(T) - u_d}^2_{L^2(\Omega)} + \alpha \mnorm{q}, \; {q \in {\mathcal M}(\Omega)}, \; \text{subject to}\;\eqref{eq:state}.
\end{equation}
Here and in what follows, $\Om$ is a convex polygonal/polyhedral domain in $\mathbb{R}^N$, $N=2,3$, $I=(0,T]$ is the time interval, $u_d \in L^2(\Omega)$ is the given (desired /measured) final state, and $\alpha>0$ is the regularization parameter. A very similar problem is considered in~\cite{CasasE_VexlerB_ZuazuaE_2015}. There, the initial state $q$ is also searched for in the space $\M(\Omega)$. For given $\varepsilon>0$ and $u_d \in L^2(\Omega)$ the optimal control problem in~\cite{CasasE_VexlerB_ZuazuaE_2015} is formulated as follows:
\begin{equation}\label{eq:problem_from_CVZ}
\text{Minimize }\; \mnorm{q}\;\text{subject to}\; \norm{u(T) - u_d}_{L^2(\Omega)}\le \varepsilon  \; \text{and}\;\eqref{eq:state}.
\end{equation}
One can directly show, that problems~\eqref{eq:problem} and~\eqref{eq:problem_from_CVZ} are equivalent by appropriate choices of $\alpha$ and $\varepsilon$.

The optimal control problem~\eqref{eq:problem} possesses a unique solution $(\oq,\ou)$, see next section for details. 
For a numerical solution of the optimal control problem under consideration we will use discontinuous Galerkin methods dG($r$) of order $r$ for temporal and linear (conforming) finite elements for spatial discretizations of the state equation~\eqref{eq:state} leading to the discrete optimal solution $(\oq_{kh},\ou_{kh})$. The same type of discretization (with $r=0$) is used in~\cite{CasasE_VexlerB_ZuazuaE_2015}, where weak-star convergence $\oq_{kh} \overset{\ast}{\rightharpoonup} \oq$ in $\M(\Omega)$ for the control and strong convergence $\ou_{kh}(T) \to \ou(T)$ in $L^\infty(\Omega)$ is shown for the discretization parameters $k$ and $h$ tending to zero. However, no convergence rates with respect to $k$ or $h$ are derived in~\cite{CasasE_VexlerB_ZuazuaE_2015}. The main goal of this paper is to close this gap and obtain precise error estimates. In addition, in the case when the optimal control is in the form of linear combination of Diracs, we obtain convergence rates for the source locations and the corresponding coefficients. We illustrate the theoretical results with numerical experiments. 

 For the general case (i.e. without any further assumptions) we will prove the following error estimate
\[
\ltwonorm{(\ou-\ou_{kh})(T)} \le c  (k^{r+\frac{1}{2}} + \ell_{kh}h),
\]
where $k$ denotes the maximal time step,  $h$ is the spatial mesh size, and $\ell_{kh}$ is a logarithmic term, see Theorem~\ref{theorem:general_estimate} for details.

From the optimality system (see next section) we will deduce, that the support of the optimal control (optimal initial state) $\oq$ is contained in the set of maxima and minima of the adjoint state $\oz(0)$, see Corollary~\ref{cor:structure_opt_control}. Under additional assumptions (Assumption~\ref{ass:finite}) on this set, which implies that the optimal control $\oq$ consists of finitely many Dirac measures, i.e.
\[
\oq = \sum_{i=1}^K \ob_i \delta_{\oxi},
\]
we will show, that the discrete optimal control $\oq_{kh}$ has a similar structure, i.e.
\[
\oq_{kh} = \sum_{i=1}^K \sum_{j=1}^{n_i} \ob_{kh,ij} \delta_{\oxkhij},
\]
where each Dirac measure $\delta_{\oxi}$ on the continuous level is approximated by $n_i\ge 1$ Dirac measures $\delta_{\oxkhij}$ on the discrete level, see Lemma~\ref{lemma:oq_finite_structure_h} for details. In this setting we will provide (see Theorem~\ref{th:improved_k_error} and Theorem~\ref{th:improved_h_error}) an improved error estimate for the optimal states, i.e.
\[
\ltwonorm{(\ou-\ou_{kh})(T)} \le c  (k^{2r+1} + \ell_{kh} h).
\]
Moreover, we will prove an estimate for the error in position of the support points, 
\[
\abs{\oxi-\oxkhij} \le c(k^{2r+1}+\ell_{kh}^\frac{1}{2} h)
\]
for all $1 \le i \le K$ and $1 \le j \le n_i$ and a corresponding estimates for the coefficients.
As a corollary we obtain an error estimate for the discrete optimal solutions in the norm on the topological dual of the Sobolev space~$W^{1,\infty}(\Omega)$. This also implies the same rate of convergence for~$\oqkh$ with respect to the Kantorovich-Rubinshtein norm,~\cite[Section~8.3]{bogachev}, given by
\begin{align}
\|q\|_{\text{KR}}=\sup \Set{\langle q,  \varphi \rangle |\varphi\in \mathcal{C}(\Omega),~\sup_{\substack{x_1,x_2 \in \Omega,\\x_1 \neq x_2}} \frac{|\varphi(x_1)-\varphi(x_2)|}{|x_1-x_2|} \leq 1,~|\varphi(x)| \leq 1,~x \in \Omega}
\end{align}
for~$q \in \M(\Omega)$.
In fact, we readily verify that this norm is equiavalent to the $(W^{1,\infty})^*$ norm.
Roughly speaking, the metric induced by the Kantorovich-Rubinshtein norm can be interpreted as an extension of the well-known Wasserstein-1 distance,~\cite{kantorovich}, which is defined for probability measures, to signed measures with different mean values.

In order to obtain such convergence rates we need to revise fully discrete pointwise smoothing error estimates for a homogeneous parabolic problem 
\begin{equation}\label{eq:eq_aux}
\begin{aligned}
\pa_t v-\Delta v &= 0 &&\text{in} \quad (0,T)\times \Om,\;  \\
v                &= 0 &&\text{on}\quad  (0,T)\times\pa\Omega, \\
v(0)             &= v_0 &&\text{in}\quad \Omega,
\end{aligned}
\end{equation}
with a general initial condition $v_0\in L^2(\Om)$. This means that for the fully discrete approximation $v_{kh}$ we need optimal pointwise spatial  error estimates for $(v-v_{kh})(T)$ in terms of the $L^2(\Omega)$ norm of the initial data. This problem is classical and was considered in a number of papers,  we only mention the most relevant ones to our presentation. Global $L^\infty(\Om)$ error estimates for smooth domains and uniform time steps were established in \cite{HansboA_2002a}, on the other hand superconvergent results at time nodes in $L^2(\Om)$ norm, again on smooth domains were established in \cite{ErikssonK_JohnsonC_ThomeeV_1985}. One of the main contributions of our paper is the derivation of superconvergent in time and pointwise in space interior error estimates on convex polygonal/polyhedral domains. More precisely, we establish the following result
\begin{equation}\label{eq: smoothing global}
|(v-v_{kh})(T,x_0)|\le C(T)\left({k^{2r+1}}+\ell_{kh}{h^2}\right)\ltwonorm{v_0},
\end{equation}
where $x_0 \in \Omega$ is an interior point. The precise form of the constants and the logarithmic terms are given in the statements of the Theorem \ref{thm:smoothing_est_time} and Theorem~\ref{thm:smoothing_est_space}. This result is required for our error analysis for the problem~\eqref{eq:problem} and is also of independent interest.

Throughout the paper we use $\abs{\cdot}$ for the absolute value and also for the Euclidian norm of a vector in ${\mathbb R}^n$. We employ the usual notation for the Lebesgue and Sobolev
spaces. We denote by $(\cdot,\cdot)$ the inner product in $L^2(\Omega)$, by $\langle\cdot,\cdot \rangle$ the duality product between $\M(\Omega)$ and $C_0(\Omega)$, and by $(\cdot,\cdot)_{J\times \Omega}$ the inner product in $L^2(J\times \Omega)$ with a subinterval $J \subset I$. With $W(0,T)$ we denote the usual space
\[
W(0,T) = L^2(I;H^1_0(\Omega)) \cap H^1(I;H^{-1}(\Omega)).
\]

The paper is organized as follows. In the next section we introduce the optimal control problem, derive first order optimality conditions and discuss structural properties of the optimal solutions. In section~\ref{sec:discretization}, we present a fully discrete scheme for the homogeneous parabolic equation~\eqref{eq:eq_aux} and state key smoothing error estimates, the proofs of which are postponed until sections~\ref{sec:proof_thm_smoothing_time} and~\ref{sec:proof_thm_smoothing_space}. In section~\ref{sec:disc_optimal_control}, we look separately at the time semidiscretization and the full discretization of the optimal control problem and derive some preliminary results. In section~\ref{sec:general_ee} we first obtain suboptimal error estimates for the general case which under additional assumptions we improve in section~\ref{sec:improved_ee}. Finally, the last two sections are devoted to the description of the algorithm and numerical illustrations of our theoretical results.

\section{Optimal control problem}\label{sec:optimal_control}

To introduce the precise formulation of the optimal control problem under the consideration we first discuss the solution of the state equation~\eqref{eq:state}. For a given $q \in \M(\Omega)$ we define a (very weak) solution $u=u(q) \in L^1(I \times \Omega)$ of~\eqref{eq:state} if the following identity holds
\[
(\psi, u)_{I \times \Omega} = \langle q,\varphi(0) \rangle
\]
for all $\psi \in L^\infty(I \times \Omega)$, where $\varphi \in W(0,T)$ is the weak solution of
\[
\begin{aligned}
-\pa_t \varphi-\Delta \varphi &= \psi &&\text{in} \quad (0,T)\times \Om,\;  \\
\varphi                &= 0 &&\text{on}\quad  (0,T)\times\pa\Omega, \\
\varphi(T)             &= 0 &&\text{in}\quad \Omega.
\end{aligned}
\]
It is well known, that $\varphi \in C(\bar I\times \bar \Omega)$ for $\psi \in  L^\infty(I \times \Omega)$, see, e.g.,~\cite[Theorem 6.8]{Griepentrog:2007} on general Lipschitz domains, or~\cite[Theorem 5.1]{Casas:1997}. Therefore, $\varphi(0) \in C_0(\Omega)$ and the solution $u$ is well defined. There holds the following proposition, see~\cite[Lemma 2.2]{CasasE_VexlerB_ZuazuaE_2015}.

\begin{prpstn}\label{prop:state}
For each $q \in \M(\Omega)$ there exists a unique solution $u$ of~\eqref{eq:state} in the above sense. Moreover, there holds $u\in L^r(I;W_0^{1,p}(\Omega))$ for all $r,p \in [1,2)$ with 
\[
\frac{2}{r} + \frac{N}{p} > N+1
\]
and $u(T) \in L^2(\Omega)$ with the corresponding estimates
\[
\norm{u}_{L^r(I;W_0^{1,p}(\Omega))} \le c \mnorm{q}
\]
and
\[
\ltwonorm{u(T)} \le c \mnorm{q}.
\]
\end{prpstn}
\begin{rmrk}\label{remark:H2_reg}
The final state $u(T)$ has more regularity. There holds $(-\Delta)^k u(T) \in L^2(\Omega)$ for any natural number $k$. For example by taking $k=1$, we obtain $u(T) \in H^2(\Omega)\cap H^1_0(\Omega)$ using the convexity of the domain.
\end{rmrk}

The unique solvability of the state equation allows us to introduce the control-to-state mapping $S \colon \M(\Omega) \to L^2(\Omega)$ with $S(q) = u(q)(T)$. By the discussion above this operator is linear continuous and due to $S(q)\in H^2(\Omega)$ it maps every weakly star converging sequence $\{q_n\}\subset \M(\Omega)$ to a strongly converging sequence in $L^2(\Omega)$. Based on this operator we define the reduced cost functional $j \colon \M(\Omega) \to \R$ by
\[
j(q) = \half \ltwonorm{S(q)-u_d}^2 + \alpha \mnorm{q}.
\]
The optimal control problem~\eqref{eq:problem} can then be formulated as
\begin{equation}\label{eq:red_problem}
\text{Minimize }\; j(q), q \in \M(\Omega).
\end{equation}

\begin{thrm}\label{theorem:existence}
The problem~\eqref{eq:red_problem} possesses a unique solution $\oq \in \M(\Omega)$. There holds the estimates
\[
\ltwonorm{\ou(T)} \le 2 \ltwonorm{u_d} \quad\text{and}\quad \alpha\mnorm{\oq} \le \frac{1}{2}\ltwonorm{u_d}^2,
\]
where $\ou = u(\oq)$ is the corresponding optimal state.
\end{thrm}
\begin{proof}
	The existence follows by standard arguments, cf., e.g,~\cite[Proposition 2.2.]{ClasonKunisch:2011}. The uniqueness follows as in~\cite[Theorem 2.4]{CasasE_VexlerB_ZuazuaE_2015} using density of the range of the semigroup
generated by the heat equation~\cite{FabrePUelZuazua:1995}, which is equivalent to the backward uniqueness property of the heat equation. The estimates follow from $j(\oq)\le j(0)$.
\end{proof}

The unique solution $\oq$ and the corresponding optimal state $\ou$ can be characterized by the following optimality conditions.

\begin{thrm}\label{theorem:opt_cond}
The control $\oq \in \M(\Omega)$ is the solution of~\eqref{eq:red_problem} if and only if the triple $(\oq,\ou,\oz)$ satisfies the following conditions:
\begin{itemize}
	\item state equation, $\ou = u(\oq)$ in the sense of Proposition~\ref{prop:state}.
	\item adjoint equation for $\oz \in W(0,T)$ being the weak solution of
	\[
	\begin{aligned}
	-\pa_t \oz-\Delta \oz &= 0, &&\text{in} \quad (0,T)\times \Om,\;  \\
	\oz                &= 0, &&\text{on}\quad  (0,T)\times\pa\Omega, \\
	\oz(T)             &= \ou(T)-u_d, &&\text{in}\quad \Omega.
	\end{aligned}
	\]
	\item variational inequality
	\[
	-\langle q-\oq, \oz(0)\rangle \le \alpha \left(\mnorm{q}-\mnorm{\oq}\right) \quad \text{for all }\; q \in \M(\Omega).
	\]
\end{itemize}
\end{thrm}
\begin{proof}
The proof is similar to~\cite[Theorem 2.1]{CasasClasonKunisch:2012}. Note, that $\oz(0)\in C_0(\Omega)$, which makes the duality product in the variational inequality well defined.
\end{proof}

The next lemma states additional regularity for $\oz(0)$.

\begin{lmm}\label{lemma:higher_reg_z}
Let $\oq \in \M(\Omega)$ be the solution of~\eqref{eq:red_problem}, $\ou$ be the corresponding state and $\oz$ the corresponding adjoint state.  Let $\Omega_0$ be an interior subdomain of $\Omega$, i.e. $\bar \Omega_0 \subset \Omega$. Then there holds $\oz(0) \in H^4(\Omega_0) \hookrightarrow C^2(\Omega_0)$ with
\[
\norm{\oz(0)}_{H^4(\Omega_0)} \le c \ltwonorm{u_d},
\] 
where the constant $c$ depends on $\Omega$, $T$ and $\Omega_0$.
\end{lmm}
\begin{proof}
As in Remark~\ref{remark:H2_reg},  one shows directly $-\Delta \oz(0) \in H^2(\Omega)\cap H^1_0(\Omega)$ with
\[
\norm{\Delta \oz(0)}_{H^2(\Omega)} \le c \ltwonorm{\Delta^2 \oz(0)} \le c \ltwonorm{\ou(T)-u_d},
\]
cf. also~\eqref{eq: continuous smoothing} below. Then the elliptic interior regularity result from~\cite[Chapter 6.3,Theorem 2]{LCEvans_2010} implies
\[
\norm{\oz(0)}_{H^4(\Omega_0)} \le c \norm{\Delta \oz(0)}_{H^2(\Omega)} \le c \ltwonorm{\ou(T)-u_d} \le c \ltwonorm{u_d},
\]
where in the last estimate we used Theorem~\ref{theorem:existence}.
\end{proof}

From the above optimality condition we obtain the following structural properties of the optimal solution $\oq$ and the corresponding optimal adjoint state $\oz$.

\begin{crllr}\label{cor:structure_opt_control}
Let $\oq$ be the solution of ~\eqref{eq:red_problem}, $\ou$ be the corresponding state and $\oz$ the corresponding adjoint state. Then there hold
\begin{itemize}
\item[(a)] a bound for the adjoint state $\oz(0)$
\[
\abs{\oz(0,x)} \le \alpha \quad \text{for all }\: x \in \bar \Omega,
\]
\item[(b)] a support condition for the positive and the negative parts in the Jordan decomposition of $\oq = \oq^+-\oq^-$
\[
\supp \oq^+ \subset \Set{x\in \Omega| \oz(0,x)=-\alpha} \:\text{ and }\: \supp \oq^- \subset \Set{x\in \Omega| \oz(0,x)=\alpha}.
\]
Moreover there is a subdomain $\Omega_0$ with $\bar \Omega_0 \subset \Omega$ such that
\[
\supp \oq \subset \Omega_0.
\]
\end{itemize}
\end{crllr}
\begin{proof}
The proof is similar to~\cite{ClasonKunisch:2011} or~\cite{CasasClasonKunisch:2012}.
\end{proof}
\begin{rmrk}
The adjoint state $\oz(0)$ is analytic on $\Omega_0$, see~\cite{Kovats:2002}. This implies by the above corollary that Lebesgue measure of $\supp \oq$ is zero. 
\end{rmrk}
\section{Discretization and smoothing type error estimates}\label{sec:discretization}

In this section we describe the (fully discrete) finite element discretization of the (axillary) homogeneous equation~\eqref{eq:eq_aux} and present smoothing type error estimates. To discretize the problem we use  continuous linear Lagrange finite elements in space and discontinuous Galerkin methods of order $r$ in time. To be more precise, we partition  $I =(0,T]$ into subintervals $I_m = (t_{m-1}, t_m]$ of length $k_m = t_m-t_{m-1}$, where $0 = t_0 < t_1 <\cdots < t_{M-1} < t_M =T$. The maximal and minimal time steps are denoted by $k =\max_{m} k_m$ and $k_{\min}=\min_{m} k_m$, respectively.
We impose the following conditions on the time mesh (as in~\cite{LeykekhmanD_VexlerB_2016a} or~\cite{DMeidner_RRannacher_BVexler_2011a}):
\begin{enumerate}[(i)]
  \item There are constants $c,\beta>0$ independent on $k$ such that
    \[
      k_{\min}\ge ck^\beta.
    \]
  \item There is a constant $\kappa>0$ independent on $k$ such that for all $m=1,2,\dots,M-1$
    \[
    \kappa^{-1}\le\frac{k_m}{k_{m+1}}\le \kappa.
    \]
  \item It holds $k\le\frac{T}{2r+2}$.
\end{enumerate}
The semidiscrete space $\Xk$ of piecewise polynomial functions in time is defined by
\[
\Xk=\Set{\varphi_k\in L^2(I;H^1_0(\Om)) | \varphi_k|_{I_m}\in \mathbb{P}_{r}(I_m;H^1_0(\Om)), \ m=1,2,\dots,M},
\]
where $\mathbb{P}_{r}(I_m;V)$ is the space of polynomial functions of degree $r$ in time om $I_m$ with values in a Banach space $V$.
We will employ the following notation for functions with possible discontinuities at the nodes $t_m$:
\begin{equation}\label{def: time jumps}
w^+_m=\lim_{\eps\to 0^+}w(t_m+\eps), \quad w^-_m=\lim_{\eps\to 0^+}w(t_m-\eps), \quad [w]_m=w^+_m-w^-_m.
\end{equation}
Next we define the following bilinear form
\begin{equation}\label{eq: bilinear form B}
 B(w,\varphi)=\sum_{m=1}^M \langle w_t,\varphi \rangle_{I_m \times \Omega} + (\na w,\na \varphi)_{\IOm}+\sum_{m=2}^M([w]_{m-1},\varphi_{m-1}^+)+(w_{0}^+,\varphi_{0}^+),
\end{equation}
where 
$\langle \cdot,\cdot \rangle_{I_m \times \Omega}$ is the duality product between $ L^2(I_m;H^{-1}(\Omega))$ and $ L^2(I_m;H^{1}_0(\Omega))$. We note, that the first sum vanishes for $w \in X^0_k$. The dG($r$) semidiscrete (in time) approximation $v_k\in X_k^q$ of \eqref{eq:eq_aux} is defined as
\begin{equation}\label{eq: semidiscrete heat}
B(v_k,\varphi_k)=(v_0,\varphi_{k,0}^+) \quad \text{for all }\; \varphi_k\in \Xk.
\end{equation}
Rearranging the terms in \eqref{eq: bilinear form B}, we obtain an equivalent (dual) expression for $B$:
\begin{equation}\label{eq:B_dual}
 B(w,\varphi)= - \sum_{m=1}^M \langle w,\varphi_t \rangle_{I_m \times \Omega} + (\na w,\na \varphi)_{\IOm}-\sum_{m=1}^{M-1} (w_m^-,[\varphi]_m) + (w_M^-,\varphi_M^-).
\end{equation}
In the sequel we require the projection operator $\pi_k$ for $w \in C(I,L^2(\Omega))$ with $\pi_k w|_{I_m} \in \mathbb{P}_r(I_m;L^2(\Omega))$ for $m=1,2,\dots,M$ on each subinterval $I_m$ by
\begin{subequations}\label{eq: projection pi_k}
\begin{equation}\label{eq: projection pi_k first}
(\pi_k w-w,\varphi)_{I_m \times \Omega}=0,\quad \text{for all }\; \varphi\in \mathbb{P}_{r-1}(I_m,L^2(\Omega)),\quad r>0,
\end{equation}
\begin{equation}\label{eq: projection pi_k second}
\pi_k w(t_m^-)=w(t_m^-).
\end{equation}
\end{subequations}
In the case $r = 0$, $\pi_k w$ is defined only by the second condition.

Next we define the fully discrete approximation scheme. For $h \in (0, h_0]$; $h_0 > 0$, let $\mathcal{T}$  denote  a quasi-uniform triangulation of $\Om$  with mesh size $h$, i.e., $\mathcal{T} = \{\tau\}$ is a partition of $\Om$ into cells (triangles or tetrahedrons) $\tau$ of diameter $h_\tau$ such that for $h=\max_{\tau} h_\tau$,
$$
\operatorname{diam}(\tau)\le h \le C |\tau|^{\frac{1}{N}}, \quad \text{for all }\; \tau\in \mathcal{T},
$$
hold. Let $V_h$ be the set of all functions in $H^1_0(\Om)$ that are affine linear on each cell $\tau$, i.e. $V_h$ is the usual space of linear  conforming finite elements. We define the following three operators to be used in the sequel: discrete Laplacian $\Delta_h \colon V_h \to V_h$ defined by
\[
(-\Delta_h v_h,w_h) =(\nabla v_h,\nabla w_h) \quad \text{for all }\; v_h,w_h \in V_h,
\]
the $L^2$ projection $P_h \colon L^2(\Omega) \to V_h$ defined by
\[
(P_h v, w_h) = ( v, w_h) \quad \text{for all }\; w_h \in V_h,
\]
and the Ritz projection $R_h \colon H^1_0(\Omega) \to V_h$ defined by
\[
(\nabla R_h v, \nabla w_h) = (\nabla v, \nabla w_h) \quad \text{for all }\; w_h \in V_h.
\]
To obtain the fully discrete approximation of~\eqref{eq:eq_aux} we consider the space-time finite element space
\begin{equation} \label{def: space_time}
\Xkh=\Set{v_{kh} \in \Xk | v_{kh}|_{I_m}\in \mathbb{P}_{q}(I_m;V_h), \ m=1,2,\dots,M}.
\end{equation}
We define a fully discrete cG($1$)dG($r$) approximation $v_{kh} \in \Xkh$ of~\eqref{eq:eq_aux} by
\begin{equation}\label{eq:fully discrete heat}
B(v_{kh},\varphi_{kh})=(v_0,\varphi_{kh}^+) \quad \text{for all }\; \varphi_{kh}\in \Xkh.
\end{equation}
Notice that we have the following orthogonality relations 
\begin{subequations}
\begin{equation}\label{eq: orthogonality semidiscrete}
B(v-v_k,\varphi_k)=0  \quad \text{for all }\; \varphi_k\in \Xk,
\end{equation}
\begin{equation}\label{eq: orthogonality fully discrete}
B(v-v_{kh},\varphi_{kh})=0  \quad \text{for all }\;  \varphi_{kh}\in \Xkh.
\end{equation}
\end{subequations}
In the proofs we will use the following truncation argument. For $w_k,\varphi_k \in \Xk$, we let $\tilde{w}_k=\chi_{(t_{\tilde{m}},T]}w_k$ and $\tilde{\varphi}_k=\chi_{(t_{\tilde{m}},T]}\varphi_{k}$, 
 where $\chi_{(t_{\tilde{m}},T]}$ is the characteristic function on the interval $(t_{\tilde{m}},T]$, for some $1\le\tilde{m}\le M$, i.e. $\tilde{w}_k=0$  on $I_1\cup \cdots \cup I_{\tilde{m}}$ for some $\tilde{m}$  and $\tilde{w}_k={w}_k$ on the remaining   time intervals.   Then from \eqref{eq: bilinear form B}, we have the identity
\begin{equation}\label{eq: Bilinear tilde}
B(\tilde{w}_{k},\varphi_{k})=B(w_{k},\tilde{\varphi}_{k})+(w_{k,\tilde{m}}^-,\varphi_{k,\tilde{m}}^+).
\end{equation}
Same identity holds of course for fully discrete functions  $w_{kh},\varphi_{kh}\in\Xkh$. The following smoothing properties of the continuous, semidiscrete and fully discrete solutions are essential in our arguments. 

\subsection{Parabolic smoothing}
It is well known that the solution $v$ to the homogeneous  problem~\eqref{eq:eq_aux} has the following smoothing property. 
\begin{equation}\label{eq: continuous smoothing}
\|\partial_t^l v(t)\|_{L^p(\Om)} + \|(-\Delta)^l v(t)\|_{L^p(\Om)}\le \frac{C}{t^l}\| v_0\|_{L^p(\Om)}\quad t>0,\quad 1\le p\le \infty,\quad l=0,1,\dots.
\end{equation}
To get smoothing estimates in some other norms, we will frequently use the Gagliardo-Nirenberg inequality
\begin{equation}\label{eq:general_Gagliardo-Nirenberg}
\norm{g}_{L^\infty(B)}\le C \norm{g}^{\frac{N}{4}}_{H^2(B)} \norm{g}^{1-\frac{N}{4}}_{L^2(B)},
\end{equation}
which holds for any subdomain $B \subset \Omega$ fulfilling the cone condition (in particular for $B=\Omega$) and for all $g \in H^2(B)$, see~\cite[Theorem 3]{AdamsFournier:1977}. For $B=\Omega$ it follows with the $H^2$-regularity
\begin{equation}\label{eq:Gagliardo-Nirenberg_Omega}
\norm{g}_{L^\infty(\Omega)}\le C \norm{\Delta g}^{\frac{N}{4}}_{L^2(\Omega)} \norm{g}^{1-\frac{N}{4}}_{L^2(\Omega)}.
\end{equation}

The following smoothing estimates can be  obtained from~\eqref{eq: continuous smoothing}.
\begin{lmm}\label{lemma:cont_smoothing}
Let $v_0 \in L^2(\Omega)$ and $v\in W(0,T)$ be the solution of~\eqref{eq:eq_aux}. Then $v(T) \in H^2(\Omega)\cap H^1_0(\Omega)$ and the following estimate holds
\[
\norm{v(T)}_{H^2(\Omega)} \le C T^{-1} \ltwonorm{v_0}.
\]
Moreover, for each interior subdomain $\Omega_0$ with $\bar \Omega_0 \subset \Omega$, the final state $v(T)$ is (real) analytic on such $\Omega_0$ and there hold
\[
\norm{\na v(T)}_{L^{\infty}(\Omega_0)} \le C T^{-\frac{1}{2}-\frac{N}{4}} \ltwonorm{v_0} \quad \text{and} \quad \norm{v(T)}_{C^2(\Omega_0)} \le C T^{-1-\frac{N}{4}} \ltwonorm{v_0}.
\]
\end{lmm}
\begin{proof}
The first inequality follows right the way from~\eqref{eq: continuous smoothing} with $l=1$ by $H^2$ regularity. The analyticity can be found, e.g., in~\cite{Kovats:2002}. To prove the second inequality we first observe that
\[
\ltwonorm{\nabla v(T)} \le C T^{-\frac{1}{2}} \ltwonorm{v_0} \quad \text{and} \quad \ltwonorm{\nabla \Delta v(T)} \le C T^{-\frac{3}{2}} \ltwonorm{v_0}.
\]
Then we use Gagliardo-Nirenberg inequality~\eqref{eq:general_Gagliardo-Nirenberg} for $g=\nabla v(T)$ and $B=\Omega_0$ resulting in
\[
\begin{aligned}
\norm{\na v(T)}_{L^{\infty}(\Omega_0)} & \le C \norm{\nabla v(T)}_{H^2(\Omega_0)}^{\frac{N}{4}} \ltwonorm{\nabla v(T)}^{1-\frac{N}{4}}\\
&\le C \ltwonorm{\nabla \Delta v(T)}^{\frac{N}{4}} \ltwonorm{\nabla v(T)}^{1-\frac{N}{4}}\\
&\le C T^{-\frac{3}{2}\cdot \frac{N}{4}-\frac{1}{2}(1-\frac{N}{4})} \ltwonorm{v_0}\\
&= C T^{-\frac{1}{2}-\frac{N}{4}}\ltwonorm{v_0},
\end{aligned}
\]
where we have used the interior regularity result~\cite[Chapter 6.3,Theorem 2]{LCEvans_2010}. To show the last inequality we use Gagliardo-Nirenberg  inequality~\eqref{eq:general_Gagliardo-Nirenberg} for $g=\nabla^2 v(T)$ and $B=\Omega_0$ resulting in
\[
\begin{aligned}
\norm{v(T)}_{C^2(\Omega_0)} &\le C\|\nabla^2 v(T)\|^{\frac{N}{4}}_{H^2(\Om_0)}\|\nabla^2 v(T)\|^{1-\frac{N}{4}}_{L^2(\Om_0)}\\
&\le C \norm{v(T)}_{H^4(\Omega_0)}^{\frac{N}{4}} \norm{v(T)}_{H^2(\Omega_0)}^{1-\frac{N}{4}}\\
& \le C \norm{\Delta v(T)}_{H^2(\Omega)}^{\frac{N}{4}} \ltwonorm{\Delta v(T)}^{1-\frac{N}{4}}\\
& \le C \ltwonorm{\Delta^2 v(T)}^{\frac{N}{4}} \ltwonorm{\Delta v(T)}^{1-\frac{N}{4}}\\
& \le C T^{-1-\frac{N}{4}} \ltwonorm{v_0},
\end{aligned}
\]
where we again have used the interior regularity result~\cite[Chapter 6.3,Theorem 2]{LCEvans_2010} and convexity of $\Omega$.
\end{proof}

For the discontinuous Galerkin methods similar smoothing type estimates also hold, see Theorems 3,4,5,10 in \cite{LeykekhmanD_VexlerB_2017a} for general $L^p$ norms, cf. also~\cite[Theorem~5.1]{ErikssonK_JohnsonC_LarssonS_1998a} for the case of the $L^2$ norm.
\begin{lmm}[Smoothing estimate]\label{lemma: homogeneous smoothing dG_r fully discrete}
 Let $v_k$ and $v_{kh}$ be the semidiscrete and fully discrete solutions  of  \eqref{eq: semidiscrete heat} and \eqref{eq:fully discrete heat}, respectively.  Then, there exists a constant $C$ independent of $k$ and $h$ such that
$$
\begin{aligned}
\sup_{t\in I_m}\|\pa_t v_{k}(t)\|_{L^p(\Om)}+\sup_{t\in I_m}\|\Delta v_{k}(t)\|_{L^p(\Om)}+k_m^{-1}\|[v_{k}]_{m-1}\|_{L^p(\Om)}&\le \frac{C}{t_m}\|v_0\|_{L^p(\Om)},\\
\sup_{t\in I_m}\|\pa_t v_{kh}(t)\|_{L^p(\Om)}+\sup_{t\in I_m}\|\Delta_h v_{kh}(t)\|_{L^p(\Om)}+k_m^{-1}\|[v_{kh}]_{m-1}\|_{L^p(\Om)}&\le \frac{C}{t_m}\|v_0\|_{L^p(\Om)},
\end{aligned}
$$
for $m=1,2,\dots,M$ and any $1\le p\le  \infty$. For $m=1$ the jump term is understood as $[v_{k}]_0 = v_{k,0}^+-v_0$ and $[v_{kh}]_0 = v_{kh,0}^+-P_h v_0$.
\end{lmm}

In addition the stability with respect to the $L^p(\Omega)$ norm is valid for the semidiscrete and fully discrete approximations of the heat equation. For the proof we refer to~\cite[Lemma~5]{LeykekhmanD_VexlerB_2017a}, see also~\cite{Palencia:1994}.
\begin{lmm}\label{lemma: stability dG_r fully discrete}
 Let $v_k$ and $v_{kh}$ be the semidiscrete and fully discrete solutions  of~\eqref{eq: semidiscrete heat} and~\eqref{eq:fully discrete heat}, respectively.  Then, there exists a constant $C$ independent of $k$ and $h$ such that
\[
\norm{v_k}_{L^\infty(I;L^p(\Omega))} \le C \norm{v_0}_{L^p(\Omega)} \quad\text{and}\quad \norm{v_{kh}}_{L^\infty(I;L^p(\Omega))} \le C \norm{v_0}_{L^p(\Omega)}
\]
holds for any $1\le p\le  \infty$.
\end{lmm}

From Lemma~\ref{lemma: homogeneous smoothing dG_r fully discrete}  we immediately obtain the following corollary. Note, that the corresponding estimate is not true on the continuous level, which explains the presence of the logarithmic term.
\begin{crllr}\label{cor: maximal parabolic initial in L1}
Under the assumptions of Lemma \ref{lemma: homogeneous smoothing dG_r fully discrete}, for  any  $1\le p\le \infty$ we have
$$
\sum_{m=1}^M \left( \|\pa_t v_{k}\|_{L^1(I_m;L^p(\Om))}+\|\Delta v_{k}\|_{L^1(I_m;L^p(\Om))}+ k_m \lpnorm{\Delta v_{k,m}^+}  + \|[v_{k}]_{m-1}\|_{L^p(\Om)}\right)\le C\lk\|v_0\|_{L^p(\Om)}.
$$
and
\begin{multline*}
\sum_{m=1}^M \Bigl(\|\pa_t v_{kh}\|_{L^1(I_m;L^p(\Om))}+\|\Delta_h v_{kh}\|_{L^1(I_m;L^p(\Om))}
+k_m \lpnorm{\Delta_h v_{kh,m}^+}+\|[v_{kh}]_{m-1}\|_{L^p(\Om)}\Bigr)\le C\lk\|v_0\|_{L^p(\Om)}.
\end{multline*}
\end{crllr}
\begin{proof}
We only provide the proof for the semidiscrete case, the fully discrete case is identical. 
Using the above smoothing result from Lemma \ref{lemma: homogeneous smoothing dG_r fully discrete}, we have
$$
\begin{aligned}
&\sum_{m=1}^M\left(\int_{I_m}\|\pa_t v_{k}(t)\|_{L^p(\Om)}dt+\int_{I_m}\|\Delta v_{k}(t)\|_{L^p(\Om)}dt+k_m \lpnorm{\Delta v_{k,m}^+} +\|[v_{k}]_{m-1}\|_{L^p(\Om)}\right)\\
&\le\sum_{m=1}^M k_m\left(\sup_{t\in I_m}\|\pa_t v_{k}(t)\|_{L^p(\Om)}+\sup_{t\in I_m}\|\Delta v_{k}(t)\|_{L^p(\Om)}+k_m^{-1}\|[v_{k}]_{m-1}\|_{L^p(\Om)}\right)\\
& \le C  \sum_{m=1}^M \frac{k_m}{t_m}\|v_{0}\|_{L^p(\Om)}\le C\lk\|v_{0}\|_{L^p(\Om)},
\end{aligned}
$$
where in the last step we used that 
\begin{equation}
\sum_{m=1}^M \frac{k_m}{t_m}\le \int_{k_1}^T\frac{dt}{t}\le C\lk.
\end{equation}
\end{proof}

For sufficiently many time steps, applying Lemma~\ref{lemma: homogeneous smoothing dG_r fully discrete} iteratively, we immediately obtain the following result.
\begin{lmm}\label{lemma: higher smoothing}
 Let $v_k$ and $v_{kh}$ be the semidiscrete and fully discrete solutions  of  \eqref{eq: semidiscrete heat} and \eqref{eq:fully discrete heat}, respectively. For any $m \in \{1,2,\dots M\}$, any $l \le m$, and any  $1\le p\le \infty$ there hold
\[
\sup_{t\in I_m}\|\pa_t(-\Delta)^{l-1} v_{k}(t)\|_{L^p(\Om)}+\sup_{t\in I_m}\|(-\Delta)^l v_{k}(t)\|_{L^p(\Om)}+k_m^{-1}\|[(-\Delta)^{l-1}v_{k}]_{m-1}\|_{L^p(\Om)}\le \frac{C}{t^l_m}\|v_0\|_{L^p(\Om)},
\]
and
\[
\sup_{t\in I_m}\|\pa_t(-\Delta_h)^{l-1} v_{kh}(t)\|_{L^p(\Om)}+\sup_{t\in I_m}\|(-\Delta_h)^l v_{kh}(t)\|_{L^p(\Om)}+k_m^{-1}\|[(-\Delta_h)^{l-1}v_{kh}]_{m-1}\|_{L^p(\Om)}\le \frac{C}{t^l_m}\|v_0\|_{L^p(\Om)},
\]
provided $k \le \frac{t_m}{l+1}$.
\end{lmm}

The next lemma is the semidiscrete analog of Lemma~\ref{lemma:cont_smoothing}.
\begin{lmm}\label{lemma:semidiscrete_smoothing}
Let $v_0 \in L^2(\Omega)$ and $v_k\in \Xk$ be the semidiscrete solution of~\eqref{eq: semidiscrete heat}.  Let $\Omega_0$ be an interior subdomain, i.e. $\bar \Omega_0 \subset \Omega$. Then $v_k(T) \in  W^{1,\infty}(\Omega_0)\cap C^2(\Omega_0)$ and the followings estimates hold
\[
\norm{\na v_k(T)}_{L^{\infty}(\Omega_0)} \le C T^{-\frac{1}{2}-\frac{N}{4}}\ltwonorm{v_0} \quad\text{and}\quad 
\norm{v_k(T)}_{C^2(\Omega_0)} \le C T^{-1-\frac{N}{4}}
\ltwonorm{v_0}.
\]
\end{lmm}
\begin{proof}
The proof is similar to the proof of Lemma~\ref{lemma:cont_smoothing} and uses Lemma~\ref{lemma: higher smoothing}.
\end{proof}

Using the discrete version of the Gagliardo-Nirenberg  inequality 
\begin{equation}\label{eq: discrete Gagliardo-Nirenber}
\|v_h\|_{L^\infty(\Om)}\le C\|\Delta_h v_h\|^{\frac{N}{4}}_{L^2(\Om)}\|v_h\|^{1-\frac{N}{4}}_{L^2(\Om)},\quad \text{for all }\; v_h\in V_h,
\end{equation}
which  for example was established for smooth domains in \cite[Lemma 3.3]{HansboA_2002a}, but the proof is valid for convex domains as well, we immediately obtain the following smoothing result.
\begin{crllr}\label{cor: discrete smoothing in Linfty}
Under the assumption of Lemma \ref{lemma: homogeneous smoothing dG_r fully discrete} for all $m=2,3,\dots, M$, we have
$$
\sup_{t\in I_m}\| v_{kh}(t)\|_{L^\infty(\Om)}\le \frac{C}{t_m^{N/4}}\|v_0\|_{L^2(\Om)}.
$$
\end{crllr}

\subsection{Smoothing pointwise error estimates}

One of the main tools in obtaining error estimates for the optimal control problem under consideration are the pointwise smoothing error estimates that have an independent interest.  The next theorems show that 
for the error at a point $(T, x_0)$ we can obtain  nearly optimal convergence rates in space and superconvergent rates in time. For elliptic problems such  interior pointwise elliptic results are known from~\cite{AHSchatz_LBWahlbin_1977a, AHSchatz_LBWahlbin_1995a}. For homogeneous parabolic problems with smoothing such results are new.  

The first theorem provides an $L^\infty(\Omega)$ error estimate for the semidiscrete error $(v-v_k)(T)$.

\begin{thrm}\label{thm:smoothing_est_time}
Let $v_0\in L^2(\Omega)$, let $v$ and $v_k$ satisfy \eqref{eq:eq_aux} and \eqref{eq: semidiscrete heat}. Then there holds
\[
\linfnorm{(v-v_k)(T)} \le C(T) k^{2r+1}\ltwonorm{v_0},
\]
with $C(T) \sim T^{-(2r+1+\frac{N}{4})}$.
\end{thrm}

Note, that we obtain here a superconvergent estimate of order ${\mathcal O}(k^{2r+1})$ for the discretization with polynomials of order $r$. The proof of this theorem  is given in Section~\ref{sec:proof_thm_smoothing_time}.

\begin{rmrk}\label{remark:application_to_dual}
In the sequel we will apply this and the following theorems for both, a heat equation formulated forward in time~\eqref{eq:eq_aux} and for a heat equation formulated backward in time, i.e.
\[
\begin{aligned}
-\pa_t y-\Delta y &= 0, &&\text{in} \quad (0,T)\times \Om,\;  \\
y                &= 0, &&\text{on}\quad  (0,T)\times\pa\Omega, \\
y(T)             &= y_T, &&\text{in}\quad \Omega
\end{aligned}
\]
for some $y_T \in L^2(\Omega)$. Its semidiscrete approximation $y_k \in \Xk$ solves
\[
B(\varphi_k,y_k) = (y_T,\varphi_k(T)) \quad \text{for all }\; \varphi_k\in \Xk.
\]
For this case the statement of the above theorem reads
\[
\linfnorm{y(0)-y_{k,0}^+} \le C(T) k^{2r+1}\ltwonorm{y_T}.
\]
Correspondingly we will apply also Theorem~\ref{thm:smoothing_est_time_grad} and Theorem~\ref{thm:smoothing_est_space} for this setting.
\end{rmrk}

A corresponding result is true also for the $L^\infty$ norm of the gradient.


\begin{thrm}\label{thm:smoothing_est_time_grad}
	Let $v_0\in L^2(\Omega)$, let $v$ and $v_k$ satisfy \eqref{eq:eq_aux} and \eqref{eq: semidiscrete heat}.  Let moreover $\Omega_0$ with $\bar \Omega_0 \subset \Omega$ be an interior subdomain. Then there holds
	\[
	\norm{\nabla(v-v_k)(T)}_{L^\infty(\Omega_0)} \le C(T) k^{2r+1}\ltwonorm{v_0},
	\]
	with $C(T) \sim T^{-(2r+\frac{3}{2}+\frac{N}{4})}$.
\end{thrm}

The proof of this theorem  is given in Section~\ref{sec:proof_thm_smoothing_time}.
\begin{rmrk}
The result of Theorem~\ref{thm:smoothing_est_time_grad} is valid also on the whole domain $\Omega$ instead of $\Omega_0$ with a slightly different constant $C(T)$.
\end{rmrk}

For the spatial error $(v_k-v_{kh})(T)$ we can not expect an ${\mathcal O}(h^2)$ estimate with respect to the global $L^\infty(\Omega)$ norm. However for a given point $x_0 \in \Omega$ we obtain the following  result.

\begin{thrm}\label{thm:smoothing_est_space}
Let $v_0\in L^2(\Omega)$, let $v_k$ and $v_{kh}$ satisfy \eqref{eq: semidiscrete heat} and \eqref{eq:fully discrete heat}, respectively and let $x_0\in \Om$ such that $dist(x_0,\partial \Om)=d$ with $d>4h$. Then there holds
\[
\abs{(v_k-v_{kh})(T,x_0)} \le C(T,d)\ell_{kh} h^2 \ltwonorm{v_0},
\]
where $\ell_{kh} = \lk + \lh$ and $C(T,d)$ is a constant, which explicit dependence on $T$ and $d$ can be tracked from the proof.
\end{thrm}

The proof of this theorem  is given in Section~\ref{sec:proof_thm_smoothing_space}.
Combining both theorems we immediately obtain an estimate for $(v-v_{kh})(T,x_0)$.
\begin{crllr}
Let $v_0\in L^2(\Omega)$, let $v$ and $v_{kh}$ satisfy \eqref{eq:eq_aux} and \eqref{eq:fully discrete heat}, respectively and let $x_0\in \Om$ such that $dist(x_0,\partial \Om)=d$ with $d>4h$. Then there holds
\[
\abs{(v-v_{kh})(T,x_0)} \le C(T,d)\left(k^{2r+1} + \ell_{kh} h^2\right) \ltwonorm{v_0},
\]
where $\ell_{kh} = \lk + \lh$.
\end{crllr}

\section{Discretization of optimal control problem}\label{sec:disc_optimal_control}
In this section we describe the temporal and spatial discretizations of the optimal control problem~\eqref{eq:problem}.

\subsection{Temporal semidiscretization}
 To introduce the associated semidiscrete state $u_k=u_k(q)$ for a given control $q\in \M(\Omega)$ we consider slightly modified semidiscrete spaces $\hXk \subset \Xk \subset \tXk$ defined by
\[
\tXk=\Set{\varphi_k\in L^2(I \times\Om) | \varphi_k|_{I_1}\in \mathbb{P}_{r}(I_1;W^ {1,s}_0(\Om)),\;\varphi_k|_{I_m}\in \mathbb{P}_{r}(I_m;H^1_0(\Om)) , \ m=2,\dots,M}
\]
and
\[
\hXk=\Set{\varphi_k\in L^2(I \times\Om) | \varphi_k|_{I_1}\in \mathbb{P}_{r}(I_1;W^ {1,s'}_0(\Om)),\;\varphi_k|_{I_m}\in \mathbb{P}_{r}(I_m;H^1_0(\Om)) , \ m=2,\dots,M}
\]
with some $1 <s< \frac{N}{N-1}$ and $s'>N$ with $\frac{1}{s}+\frac{1}{s'} = 1$. For this setting we have $\varphi_{k,0}^+ \in C_0(\Omega)$ for all $\varphi_k\in \hXk$ due to the embedding $W^ {1,s'}_0(\Om) \hookrightarrow C_0(\Omega)$. The bilinear form $B(\cdot,\cdot)$ from \eqref{eq: bilinear form B} can be extended to $\tXk \times \hXk$. This allows us to define the semidiscrete state $u_k(q)\in \tXk$ by
\begin{equation}\label{eq:state_k}
B(u_k,\varphi_k)=\langle q,\varphi_{k,0}^+\rangle \quad \text{for all }\; \varphi_k\in \hXk.
\end{equation}
The corresponding semidiscrete control-to-state mappings $S_k \colon \M(\Omega)\to L^2(\Omega)$ is given by $S_k(q) = u_k(q)(T)$ and the semidiscrete reduced cost functional $j_k \colon \M(\Omega) \to \R$ by
\[
j_k(q) = \half \ltwonorm{S_k(q)-u_d}^2 + \alpha \mnorm{q}.
\]
With this reduced cost functional we formulate the semidiscrete  optimal control problems without discretization of the control space as follows:
\begin{equation}\label{eq:red_problem_k}
\text{Minimize }\; j_k(q), \quad q \in \M(\Omega).
\end{equation}

As on the continuous level we obtain the existence of a solution to~\eqref{eq:red_problem_k}.
\begin{thrm}\label{theorem:existence_k}
The problem~\ref{eq:red_problem_k} possesses at least one solution $\oq_k \in \M(\Omega)$ with corresponding state $\ou_k=u_k(\oq_k)$.  There hold the estimates
\[
\ltwonorm{\ou_k(T)} \le 2 \ltwonorm{u_d} \quad\text{and}\quad \alpha\mnorm{\oq_k} \le \frac{1}{2}\ltwonorm{u_d}^2.
\]
\end{thrm}
\begin{proof}
	The existence and the estimates follow by standard arguments, as on the continuous level. 
\end{proof}

The question of uniqueness of $\oq_k$ is more involved and is discussed after the statement of the optimality system. 

\begin{thrm}\label{theorem:opt_cond_k}
The control $\oq_k \in \M(\Omega)$ is a solution of~\eqref{eq:red_problem_k} if and only if the triple $(\oq_k,\ou_k,\oz_k)$ fulfills the following conditions:
\begin{itemize}
	\item semidiscrete state equation, $\ou_k = u_k(\oq_k)\in \tXk$ in the sense of~\eqref{eq:state_k}.
	\item semidiscrete adjoint equation for $\oz_k \in \hXk$ being the  solution of
	\[
         B(\varphi_k,\oz_k) = (\ou_k(T)-u_d,\varphi_k(T)) \quad \text{for all }\; \varphi_k\in \tXk
	\]
	\item variational inequality
	\[
	-\langle q-\oq_k, \oz_{k,0}^+\rangle \le \alpha \left(\mnorm{q}-\mnorm{\oq_k}\right) \quad \text{for all }\; q \in \M(\Omega).
	\]
\end{itemize}
\end{thrm}
\begin{proof}
The proof is the same as for the continuous problem.
\end{proof}

\begin{crllr}\label{cor:opt_cond_k}
Let $\oq_k \in \M(\Omega)$ be a solution of~\eqref{eq:red_problem_k}, $\ou_k \in \tXk$ be the corresponding state, and $\oz_k \in \hXk$ the corresponding adjoint state. Then there hold
\begin{itemize}
\item[(a)] a bound for the adjoint state $\oz_{k,0}^+$
\[
\abs{\oz_{k,0}^+(x)} \le \alpha \quad \text{for all }\: x \in \bar \Omega,
\]
\item[(b)] a support condition for the positive and the negative parts in the Jordan decomposition of $\oq_k = \oq^+_k-\oq^-_k$
\[
\supp \oq^+_k \subset \Set{x\in \Omega| \oz_{k,0}^+(x)=-\alpha} \:\text{ and }\: \supp \oq^-_k \subset \Set{x\in \Omega| \oz_{k,0}^+(x)=\alpha}.
\]
Moreover there is a subdomain $\Omega_0$ with $\bar \Omega_0 \subset \Omega$ such that
$\supp \oq_k \subset \Omega_0.$
\end{itemize}
\end{crllr}
\begin{proof}
The proof is the same as for the continuous problem.
\end{proof}

The uniqueness of the solution $\oq$ on the continuous level follows (cf.~\cite[Theorem 2.4]{CasasE_VexlerB_ZuazuaE_2015}) by the fact that for the solution of the heat equation~\eqref{eq:eq_aux} we have that $v(T)=0$ implies $v_0=0$. This is also true for the $dG(0)$ discretization but is in general wrong for the dG($r$) semidiscretization with  $r\ge 1$. However, the following technical lemma allows us to prove uniqueness of the semidiscrete control $\oq_k$.

\begin{lmm}\label{lemma:discrete_back_unique}
Let $q \in \M(\Omega)$ and $u_k = u_k(q) \in \tXk$ be the corresponding semidiscrete state defined by~\eqref{eq:state_k}. Let $u_k(T)=0$. Then the holds:
\begin{enumerate}
	\item For $r=0$ we have $q =0$.
	\item Let $r>0$. If there exists an open set $D \subset \Omega$ such that $q|_D = 0$, then $q=0$.
\end{enumerate}
\end{lmm}
\begin{proof}
It is well known, cf., e.g.,~\cite{ErikssonK_JohnsonC_ThomeeV_1985}, that dG($r$) discretization of a homogeneous problem coincides with the corresponding subdiagonal Pad\'e approximation scheme. Therefore, there is a rational function $f_r = a_r/b_r$ with polynomials $a_r \in \mathbb{P}_{r}$, $b_r \in \mathbb{P}_{r+1}$ and $b_r(s) \neq 0$ for $s \in \R_+$, such that
\[
u_{k,1}^- = f_r(-k_1 \Delta)q, \quad \text{and}\quad u_{k,m}^- = f_r(-k_m \Delta)u_{k,m-1}^-, \quad m=2,3,\dots M.
\]
By the assumption of the lemma we have $u_{k,M}^-=u_k(T)=0$.
\begin{enumerate}
\item For $r=0$ we have  $f_0(s) = \frac{1}{1+s}$ and therefore
\[
(\operatorname{Id} - k_M \Delta)u_{k,M}^- = u_{k,M-1}^-,
\]
which implies $u_{k,M-1}^-=0$. Similarly, we obtain $u_{k,m}^-=0$ for all $m=2,3,\dots M$ and consequently $q=0$. 
\item For $r>0$ we argue  differently. We consider the eigenvalues $0 < \lambda_1 \le \lambda_2 \le \lambda_3  \dots$ of $-\Delta$ and the corresponding system of eigenfunctions $w_1,w_2,\dots$ with $(w_i,w_j)= \delta_{ij}$. The initial condition $q\in \M(\Omega)\subset H^{-2}(\Omega)$ can be expanded as
\begin{equation}\label{eq:q_expn}
q = \sum_{n=1}^\infty q_n w_n, \quad \text{with}\quad q_n = \langle q, w_n\rangle
\end{equation}
and the convergence to be understood in $H^{-2}(\Omega)$. We define the polynomials
\[
A_k(s)=\prod_{m=1}^M a_m(k_m s) \quad\text{and}\quad B_k(s)=\prod_{m=1}^M b_m(k_m s). 
\]
With this notation we have
\[
\left(B_k(-\Delta)\right)^{-1} A_k(-\Delta) q = u_{k,M}^- = 0
\]
and consequently $A_k(-\Delta) q = 0$. This results is
\[
\sum_{n=1}^\infty q_n A_k(\lambda_n) w_n = 0
\]
and therefore
\[
q_n A_k(\lambda_n) = 0 \quad \text{for all } n \in {\mathbb N}.
\]
Assume now that $q \neq 0$. Since $A_k \in {\mathbb P}_{rM}$ has no more than $rM$ positive zeros, there are only finitely many $q_n$ with $q_n \neq 0$. For this reason we have that the expansion~\eqref{eq:q_expn} is a finite sum and therefore $q \in H^2(\Omega)\cap H^1_0(\Omega)$ since $w_n \in H^2(\Omega)\cap H^1_0(\Omega)$ for every $n$ by convexity of $\Omega$. We have with some $R \in \mathbb{N}$
\[
q = \sum_{i=1}^R q_{n_i} w_{n_i}, \quad n_1 < n_2 < \dots <n_R, \quad q_{n_i} \neq 0. 
\]
From $q \in H^2(\Omega)$ and $q|_D=0$ we obtain that $(-\Delta)^l q$ also vanishes on $D$ for every $l \in \mathbb{N}$. Therefore, we have
\[
(-\Delta)^l q = \sum_{i=1}^R \lambda_{n_i}^l q_{n_i} w_{n_i} = 0 \quad \text{on }\; D
\]
and dividing by $\lambda_{n_R}^l$ we have
\[
\sum_{i=1}^R \left(\frac{\lambda_{n_i}}{\lambda_{n_R}}\right)^l q_{n_i} w_{n_i} = 0 \quad \text{on }\; D.
\] 
For $l\to \infty$ all summands with $\lambda_{n_i}<\lambda_{n_R}$ converge to zero resulting in
\[
w = \sum_{i :\lambda_{n_i}=\lambda_{n_R} } q_{n_i} w_{n_i} = 0 \quad \text{on }\; D.
\]
This $w \neq 0$ is an eigenfunction of $-\Delta$, which provides a contradiction, since a nontrivial eigenfunction can not vanish on an open set by the unique 
continuation principle, see, e.g.,~\cite[p. 64]{Leis:1997}. This completes the proof.
\end{enumerate}
\end{proof}

\begin{thrm}
The solution $\oq_k \in \M(\Omega)$ of~\eqref{eq:red_problem_k} is unique.
\end{thrm}
\begin{proof}
We first observe the uniqueness of $\ou_k(T)$ by the strict convexity of the tracking term in $j_k(q)$. It remains to show, that this implies the uniqueness of $\oq_k$. Assume there are two optimal controls and consider the difference $q:=\oq_{k,1}-\oq_{k,2} \in \M(\Omega)$. Let $w_k = u_k(\oq_1)-u_k(\oq_2)$, i.e. $w_k = u_k(q)$. Then there holds $u_k(T)=0$. In the case $r=0$ we immediately obtain $q=0$ by the first statement of Lemma~\ref{lemma:discrete_back_unique}. For $r>0$ we obtain from Corollary~\ref{cor:opt_cond_k}, that $\supp \oq_{k,i} \subset \Omega_0$ with $\bar \Omega_0 \subset \Omega$ and therefore $\supp q \subset \Omega_0$. This implies the existence of an open set $D \subset \Omega \setminus \bar \Omega_0$ with $q|_D=0$. Then we obtain $q=0$ from the second statement of Lemma~\ref{lemma:discrete_back_unique}.
\end{proof}

\subsection{Space-time discretization}
 For a given control $q\in \M(\Omega)$ we also introduce the associated fully discrete state $u_{kh}=u_{kh}(q)\in \Xkh$ by
\begin{equation}\label{eq:state_kh}
B(u_{kh},\varphi_{kh})=\langle q,\varphi_{kh}^+\rangle \quad \text{for all }\; \varphi_{kh}\in \Xkh,
\end{equation}
the fully discrete control-to-state mappings $S_{kh} \colon \M(\Omega)\to L^2(\Omega)$ by $S_{kh}(q) = u_{kh}(q)(T)$, and the fully discrete reduced cost functional $j_{kh} \colon \M(\Omega) \to \R$ by
\[
j_{kh}(q) = \half \ltwonorm{S_{kh}(q)-u_d}^2 + \alpha \mnorm{q}.
\]
Based on this definition we formulate the corresponding optimal control problem, where we first look for the control variable in the whole space $\M(\Omega)$. This leads to the following formulation.
\begin{equation}\label{eq:red_problem_kh}
\text{Minimize }\; j_{kh}(q), \quad q \in \M(\Omega).
\end{equation}
 One can not expect, that this problem has a unique solution. For $r=0$ however, where is a unique solution in the properly defined discrete subspace $\M_h$ of $\M(\Omega)$, see the discussion below. To introduce the space $\M_h$, let ${\mathcal N}_h$ be the set of all interior nodes of the mesh $\mathcal{T}$. For $x_i \in {\mathcal N}_h$ let $\delta_{x_i} \in \M(\Omega)$ denote the Dirac measure concentrated in $x_i$ and $\varphi_{h,i} \in V_h$ be the nodal basis function associated to the node $x_i$. Then we define the space $\M_h$ as
\[
\M_h = \operatorname{span}\Set{\delta_{x_i} | x_i \in {\mathcal N}_h} \subset \M(\Omega)
\]
and introduce a projection operator $\Lambda_h \colon \M(\Omega) \to \M_h$ (cf., e.g.,~\cite{CasasClasonKunisch:2012}) by
\[
\Lambda_h q  = \sum_{x_i \in {\mathcal N}_h} \langle q, \varphi_{h,i}\rangle \delta_{x_i}.
\]
The definition implies that
\begin{equation}\label{eq:prop_Lambda_h}
\langle \Lambda_h q, w \rangle = \langle q, i_h w \rangle \quad \text{for all }\; q \in \M(\Omega), \; w \in C_0(\Omega),
\end{equation} 
where $i_h \colon C_0(\Omega) \to V_h$ is the nodal interpolation operator. The following two properties of $\Lambda_h$ can be directly checked.
\begin{lmm}
There holds
\begin{itemize}
\item[(a)] $\mnorm{\Lambda_h q } \le \mnorm{q}$ for all $q \in \M(\Omega)$.
\item[(b)] The fully discrete solutions of the state equation associated with $q$ and with $\Lambda_h q$ are the same, i.e.
\[
u_{kh}(\Lambda q) = u_{kh}(q) \quad \text{for all }\; q \in \M(\Omega).
\]
\end{itemize}
\end{lmm}
\begin{proof}
The proof of (a) follows from~\cite[Thm.~3.1]{CasasE_ClasonC_KunischK_2012} and the proof of (b) uses the definition~\eqref{eq:state_kh} of $u_{kh}$ and~\eqref{eq:prop_Lambda_h}.
\end{proof}

The next theorem provides the existence of a solution to~\eqref{eq:red_problem_kh}.
\begin{thrm}\label{theorem:existence_kh}
There exists a solution of~\eqref{eq:red_problem_kh}.  For each solution $\tilde q_{kh} \in \M(\Omega)$ the projection $\oq_{kh} = \Lambda_h \tilde q_{kh} \in \M_h$ is also a solution of~\eqref{eq:red_problem_kh}. For $r=0$ the solution $\oq_{kh} \in \M_h$ is unique.
 For any solution $\oq_{kh} \in \M_h$ and the corresponding state $\ou_{kh}$ the following estimates hold
\[
\ltwonorm{\ou_{kh}(T)} \le 2 \ltwonorm{u_d} \quad\text{and}\quad \alpha\mnorm{\oq_{kh}} \le \frac{1}{2}\ltwonorm{u_d}^2.
\]
\end{thrm}
\begin{proof}
The existence and the estimates follow as on the continuous level. The fact that $\oq_{kh} = \Lambda_h \tilde q_{kh} \in \M_h$ is also a solution of~\eqref{eq:red_problem_kh} follows directly from Lemma~\ref{eq:prop_Lambda_h}. The uniqueness in the case of $r=0$ follows from the fully discrete analog of the first statement of Lemma~\ref{lemma:discrete_back_unique}, cf. also the proof of~\cite[Theorem 4.8]{CasasE_VexlerB_ZuazuaE_2015}.
\end{proof}

\begin{rmrk}
For $r>0$ it seems that problem~\eqref{eq:red_problem_kh} may in general have multiple solutions in $\M_h$.  The argument we used to prove uniqueness of the semidiscrete solution $\oq_k$ is based on the second statement of Lemma~\ref{lemma:discrete_back_unique}, which does not extend to the fully discrete setting.
\end{rmrk}

In the next theorem we state the optimality system on the fully discrete level. 
\begin{thrm}\label{theorem:opt_cond_kh}
The control $\oq_{kh} \in \M_h$ is a solution of~\eqref{eq:red_problem_kh} in $\M_h$ if and only if the triple $(\oq_{kh},\ou_{kh},\oz_{kh})$ fulfills the following conditions:
\begin{itemize}
	\item fully discrete state equation, $\ou_{kh} = u_{kh}(\oq_{kh})\in \Xkh$ in the sense~\eqref{eq:state_kh}.
	\item fully discrete adjoint equation for $\oz_{kh} \in \Xkh$ being the  solution of
	\[
         B(\varphi_{kh},\oz_{kh}) = (\ou_{kh}(T)-u_d,\varphi_{kh}(T)) \quad \text{for all }\; \varphi_{kh}\in \Xkh.
	\]
	\item variational inequality
	\[
	-\langle q-\oq_{kh}, \oz_{kh,0}^+\rangle \le \alpha \left(\mnorm{q}-\mnorm{\oq_{kh}}\right) \quad \text{for all }\; q \in \M(\Omega).
	\]
\end{itemize}
\end{thrm}
\begin{proof}
The proof is the same as for the continuous problem.
\end{proof}
\begin{rmrk}
Please note, that the variational inequality in the above theorem holds for all variations $q\in \M(\Omega)$ and not only for those from $\M_h$. This is due to the fact that the solution $\oq_{kh} \in \M_h$ solves the problem~\eqref{eq:red_problem_kh}, where the control is not discretized, see Theorem~\ref{theorem:existence_kh}.
\end{rmrk}

\begin{crllr}\label{cor:opt_cond_kh}
Let $\oq_{kh} \in \M_h$ be a solution of ~\eqref{eq:red_problem_kh}, $\ou_{kh} \in \Xkh$ be the corresponding state and $\oz_{kh} \in \Xkh$ the corresponding adjoint state. Then there hold
\begin{itemize}
\item[(a)] a bound for the adjoint state $\oz_{kh,0}^+$
\[
\abs{\oz_{kh,0}^+(x)} \le \alpha \quad \text{for all }\: x \in \bar \Omega,
\]
\item[(b)] a support condition for the positive and the negative parts in the Jordan decomposition of $\oq_{kh} = \oq^+_{kh}-\oq^-_{kh}$
\[
\supp \oq^+_{kh} \subset \Set{x \in {\mathcal N}_h| \oz_{kh,0}^+(x)=-\alpha} \:\text{ and }\: \supp \oq^-_{kh} \subset \Set{x\in {\mathcal N}_h| \oz_{kh,0}^+(x)=\alpha}.
\]
Moreover there is a subdomain $\Omega_0$ independent on $k$ and $h$ with $\bar \Omega_0 \subset \Omega$ such that
$
\supp \oq_{kh} \subset \Omega_0.
$
\end{itemize}
\end{crllr}
\begin{proof}
The proof is the same as for the continuous problem.
\end{proof}
\section{General error estimates for the optimal control problem}\label{sec:general_ee}
In this section we prove an error estimate for the error between the optimal state on the continuous and on the discrete level, which does not require any further assumptions on the structure of the solution.

As the first step we provide an estimate for the error in the state at terminal time for a given control $q\in \M(\Omega)$.

\begin{lmm}\label{lemma:est_state}
Let $q\in \M(\Omega)$ be a given control with $\supp q \subset \Omega_0$ and $\bar \Omega_0 \subset \Omega$. Let $u = u(q)$ be the solution of the state equation~\eqref{eq:state}, $u_k = u_k(q) \in \Xk$ be the semidiscrete approximation~\eqref{eq:state_k} and $u_{kh}=u_{kh}(q)\in \Xkh$ the fully discrete approximation~\eqref{eq:state_kh}. Then there hold
\[
\ltwonorm{(u-u_k)(T)}\le C(T) k^{2r+1} \mnorm{q}
\]
and
\[
\ltwonorm{(u_k-u_{kh})(T)}\le C(T)\ell_{kh} h^2 \mnorm{q},
\]
where $\ell_{kh} = \lk + \lh$.
\end{lmm}
\begin{proof}
To prove the first estimate we consider  the solution $y \in W(0,T)$ of the dual problem
\[
\begin{aligned}
-\pa_t y-\Delta y &= 0, &&\text{in} \quad (0,T)\times \Om,\;  \\
y                &= 0, &&\text{on}\quad  (0,T)\times\pa\Omega, \\
y(T)             &= (u-u_k)(T), &&\text{in}\quad \Omega
\end{aligned}
\]
and its semidiscrete approximation $y_k \in \hXk$ solving
\[
B(\varphi_k,y_k) = ((u-u_k)(T),\varphi_k(T)) \quad \text{for all }\; \varphi_k\in \tXk
\]
There holds
\[
\begin{aligned}
\ltwonorm{(u-u_k)(T)}^2 &= B(u,y) - B(u_k,y_k) = \langle q, y(0)-y_{k,0}^+ \rangle\\
& \le \mnorm{q} \linfnorm{y(0)-y_{k,0}^+} \le C k^{2r+1} \mnorm{q} \ltwonorm{(u-u_k)(T)},
\end{aligned}
\]
where in the last step we used Theorem~\ref{thm:smoothing_est_time} for the error $y(0)-y_{k,0}^+$ in the $L^\infty(\Omega)$ norm, see also Remark~\ref{remark:application_to_dual}.

For the proof of the spatial estimate we consider the dual solution $w_k \in \hXk$ solving
\[
B(\varphi_k,w_k) = ((u_k-u_{kh})(T),\varphi_k(T)) \quad \text{for all }\; \varphi_k\in \tXk.
\]
and $w_{kh}\in \Xkh$ solving
\[
B(\varphi_{kh},w_{kh}) = ((u_k-u_{kh})(T),\varphi_{kh}(T)) \quad \text{for all }\; \varphi_{kh}\in \Xkh.
\]
Then we get
\[
\begin{aligned}
\ltwonorm{(u_k-u_{kh})(T)}^2 &= B(u_k,w_k) - B(u_{kh},w_{kh}) = \langle q, w_{k,0}^+-w_{kh,0}^+ \rangle\\
& \le \mnorm{q} \norm{w_{k,0}^+-w_{kh,0}^+}_{L^\infty(\Omega_0)} \le C \ell_{kh} h^2 \mnorm{q} \ltwonorm{(u_k-u_{kh})(T)},
\end{aligned}
\]
where we used the fact that $\supp q \subset \Omega_0$ and Theorem~\ref{thm:smoothing_est_space} in the last step.

\end{proof}
\begin{rmrk}
Please note that the assumption $\supp q \subset \Omega_0$ with $\bar \Omega_0 \subset \Omega$ in the above theorem is required only for the spatial estimate.
\end{rmrk}

Based on this theorem we can directly obtain estimates for optimal values of the cost functional.

\begin{thrm}
Let $\oq \in \M(\Omega)$ be the optimal solution of~\eqref{eq:red_problem} with the corresponding optimal state $\ou$. Let $\oq_k \in \M(\Omega)$ be the optimal solution of the semidiscrete problem~\eqref{eq:red_problem_k} with the corresponding state $\ou_k \in \Xk$ and let $\oq_{kh} \in \M_h$ be a solution of the fully discrete problem~\eqref{eq:red_problem_kh} with the corresponding state $\ou_{kh} \in \Xkh$. Then there hold:
\[
\abs{j(\oq)-j_k(\oq_k)} \le C k^{2r+1}
\]
and
\[
\abs{j_k(\oq_k)-j_{kh}(\oq_{kh})} \le C \ell_{kh} h^2,
\]
where  $\ell_{kh} = \lk + \lh$ and $C=C(T,u_d)$ depends on $T$ and $\ltwonorm{u_d}$.
\end{thrm}
\begin{proof}
By the optimality of $\oq$ for~\eqref{eq:red_problem} we have
\[
j(\oq)-j_k(\oq_k)\le j(\oq_k) - j_k(\oq_k)
\]
Similarly by the optimality of $\oq_k$ for~\eqref{eq:red_problem} we have
\[
j(\oq)-j_k(\oq_k)\ge j(\oq) - j_k(\oq)
\]
and therefore
\[
\abs{j(\oq)-j_k(\oq_k)} \le \max\left(\abs{j(\oq_k) - j_k(\oq_k)},\abs{j(\oq) - j_k(\oq)}\right).
\]
For both $q=\oq$ and $q=\oq_k$ we estimate
\[
\begin{aligned}
\abs{j(q) - j_k(q)} &= \frac{1}{2}\Abs{\ltwonorm{u(q)-u_d}^2 - \ltwonorm{u_k(q)-u_d}^2}\\
&= \frac{1}{2} \Abs{(u(q) - u_k(q), u(q) + u_k(q) -2 u_d)}\\
&\le\frac{1}{2} \ltwonorm{u(q) - u_k(q)} \ltwonorm{u(q) + u_k(q) -2 u_d}.
\end{aligned}
\]
Then using the first estimate from Lemma~\ref{lemma:est_state}, the estimates 
\[
\ltwonorm{u(q)} \le C \mnorm{q} \quad \text{and} \quad \ltwonorm{u_k(q)} \le C \mnorm{q}
\]
as wells as estimates for $\oq$ and $\oq_k$ from Theorem~\ref{theorem:existence} and Theorem~\ref{theorem:existence_k} we complete the proof for the temporal error. The spatial error is estimated similarly by using the second estimate from Lemma~\ref{lemma:est_state}.
\end{proof}

The next theorem is the main result of this section, which provides an estimate for the error between the optimal states.

\begin{thrm}\label{theorem:general_estimate}
Let $\oq \in \M(\Omega)$ be the optimal solution of~\eqref{eq:red_problem} with the corresponding optimal state $\ou$. Let $\oq_k \in \M(\Omega)$ be the optimal solution of the semidiscrete problem~\eqref{eq:red_problem_k} with the corresponding state $\ou_k \in \Xk$ and let $\oq_{kh} \in \M_h$ be a solution of the fully discrete problem~\eqref{eq:red_problem_kh} with the corresponding state $\ou_{kh} \in \Xkh$. Then there hold:
\[
\ltwonorm{(\ou-\ou_k)(T)}\le C k^{r+\frac{1}{2}}
\]
and
\[
\ltwonorm{(\ou_k-\ou_{kh})(T)}\le C \ell_{kh}^\frac{1}{2} h,
\]
where $\ell_{kh} = \lk + \lh$ and $C=C(T,\alpha,u_d)$ depends on $T$, $\alpha$ and $\ltwonorm{u_d}$.
\end{thrm}
\begin{proof}
To prove the first estimate we use the variational inequality from Theorem~\ref{theorem:opt_cond} with $q = \oq_k$
\[
-\langle \oq_k-\oq, \oz(0)\rangle \le \alpha \left(\mnorm{\oq_k}-\mnorm{\oq}\right)
\]
and the corresponding variational inequality from Theorem~\ref{theorem:opt_cond_k} with $q = \oq$
\[
-\langle \oq-\oq_k, \oz_{k,0}^+\rangle \le \alpha \left(\mnorm{\oq}-\mnorm{\oq_k}\right).
\]
Adding these two inequalities results in
\[
\langle \oq_k-\oq, \oz(0)-\oz_{k,0}^+\rangle \ge 0.
\]
To proceed we introduce $\hat u_k = u_k(\oq)\in \Xk$ as the solution of~\eqref{eq:state_k} for $q = \oq$ and $\hat z_k \in \Xk$ fulfilling
\begin{equation}\label{eq:hat_z_k}
B(\varphi_k,\hat z_k) = (\ou(T)-u_d,\varphi_k(T)) \quad \text{for all }\; \varphi_k\in \Xk.
\end{equation}
Using the semidiscrete state and adjoint equations we obtain
\[
\begin{aligned}
0 &\le \langle \oq_k-\oq, \oz(0)-\oz_{k,0}^+\rangle\\
  & =  \langle \oq_k-\oq, \oz(0)- \hat z_{k,0}^+\rangle + \langle \oq_k-\oq, \hat z_{k,0}^+ -\oz_{k,0}^+\rangle\\
  & =  \langle \oq_k-\oq, \oz(0)- \hat z_{k,0}^+\rangle + B(\ou_k - \hat u_k,\hat z_k - \oz_k )\\
  & =  \langle \oq_k-\oq, \oz(0)- \hat z_{k,0}^+\rangle + ((\ou-\ou_k)(T),(\ou_k-\hat u_k)(T))\\
  & =  \langle \oq_k-\oq, \oz(0)- \hat z_{k,0}^+\rangle - \ltwonorm{(\ou-\ou_k)(T)}^2 + ((\ou-\ou_k)(T),(\ou-\hat u_k)(T)).
\end{aligned}
\]
This results in
\[
\ltwonorm{(\ou-\ou_k)(T)}^2 \le \langle \oq_k-\oq, \oz(0)- \hat z_{k,0}^+\rangle + ((\ou-\ou_k)(T),(\ou-\hat u_k)(T)). 
\]
By the Cauchy-Schwarz inequality in the last term and absorbing $\ltwonorm{(\ou-\ou_k)(T)}$ in the left-hand side we obtain
\begin{equation}\label{eq:starting_estimate_k}
\ltwonorm{(\ou-\ou_k)(T)}^2 \le 2 \langle \oq_k-\oq, \oz(0)- \hat z_{k,0}^+\rangle + \ltwonorm{(\ou-\hat u_k)(T)}^2.
\end{equation}
Using the estimates for $\mnorm{\oq}$ and $\mnorm{\oq_k}$ from Theorem~\ref{theorem:existence} and Theorem~\ref{theorem:existence_k} we get
\[
\ltwonorm{(\ou-\ou_k)(T)}^2 \le C \ltwonorm{u_d}^\frac{1}{2} \linfnorm{\oz(0)- \hat z_{k,0}^+} + \ltwonorm{(\ou-\hat u_k)(T)}^2.
\]
For the term $\linfnorm{\oz(0)- \hat z_{k,0}^+}$ we can directly apply Theorem~\ref{thm:smoothing_est_time} resulting in
\[
\linfnorm{\oz(0)- \hat z_{k,0}^+} \le C(T) k^{2r+1} \ltwonorm{\ou(T)-u_d} \le  C(T) k^{2r+1} \ltwonorm{u_d}.
\]
The term $\ltwonorm{(\ou-\hat u_k)(T)}$ is estimated by the first estimate in Lemma~\ref{lemma:est_state} leading to
\[
\ltwonorm{(\ou-\hat u_k)(T)} \le C(T)k^{2r+1} \mnorm{\oq} \le C(T)k^{2r+1} \ltwonorm{u_d}^\frac{1}{2}.
\]
Putting these estimates together we obtain
\[
\ltwonorm{(\ou-\ou_k)(T)}^2 \le C k^{2r+1},
\]
which is the the first desired estimate. The estimate for $(\ou_k-\ou_{kh})(T)$ is obtained similarly using Theorem~\ref{thm:smoothing_est_space} and the second estimate from Lemma~\ref{lemma:est_state}.
\end{proof}

For the error in the control we can in general only expect a weak star convergence, see the following lemma.
\begin{lmm}\label{lemma:weak_star_conv}
Let $\oq \in \M(\Omega)$ be the optimal solution of~\eqref{eq:red_problem}, $\oq_k \in \M(\Omega)$ be the optimal solution of the semidiscrete problem~\eqref{eq:red_problem_k}, and $\oq_{kh} \in \M_h$ be a solution of the fully discrete problem~\eqref{eq:red_problem_kh}. Then there holds
\[
\oq_k \overset{\ast}{\rightharpoonup} \oq \quad \text{in }\;\M(\Omega) \quad \text{for }\; k \to 0
\]
and for fixed $k>0$
\[
\oq_{kh} \overset{\ast}{\rightharpoonup} \oq_k \quad \text{in }\;\M(\Omega) \quad \text{for }\; h \to 0.
\]
\end{lmm}
\begin{proof}
The proof is similar to the proof of~\cite[Theorem 4.10]{CasasE_VexlerB_ZuazuaE_2015}.
\end{proof}

Under an additional assumption stronger results are discussed in the next section.
\section{Improved error estimates for the optimal state and control}\label{sec:improved_ee}
In the previous section we provided error estimates for the error in the cost functional and for the optimal states at the terminal time. In general we can not expect an error estimate for the control, $\oq-\oq_{kh}$, with respect to the norm in $\M(\Omega)$, since only weak star convergence of the controls can be expected, cf. the corresponding discussion in~\cite{CasasE_VexlerB_ZuazuaE_2015}. However, if the optimal control consists of finitely many Diracs, error rates for the positions and the coefficients of these Diracs can be shown. To prove such error estimates and to improve the estimate for the state from Theorem~\ref{theorem:general_estimate} we make the following assumption.
\begin{assumption}\label{ass:finite}
Let $\oq$ be the solution of the problem of~\eqref{eq:red_problem} with the corresponding optimal state $\ou$ and adjoint state $\oz$. 
We assume that
\[
\supp \oq = \Set{x \in \Omega | \abs{\oz(0,x)}=\alpha} = \Set{\bar x_1,\bar x_2,\dots,\bar x_K}
\]
with $K \in \mathbb{N}$ and $\oxi \in \Omega$ for $i=1,2,\dots, K$ are pairwise disjoint points. Moreover, there holds
\[
\nabla^2 \oz(0,\oxi) \;\text{is positive definite for } \oxi \text{ with } \oz(0,\oxi) = -\alpha 
\]
and
\[
\nabla^2 \oz(0,\oxi) \;\text{is negative definite for } \oxi \text{ with } \oz(0,\oxi) = \alpha,
\]
where $\nabla^2 \oz(0,\oxi)$ denotes the Hessian matrix of $\oz$ with respect to the spatial variable.
\end{assumption}
\begin{rmrk}
\begin{itemize}
\item From Corollary~\ref{cor:structure_opt_control} (b) we have that
\[
\supp \oq \subset \Set{x \in \Omega | \abs{\oz(0,x)}=\alpha}.
\]
Here, we assume equality of these two sets and the finite cardinality of them.
\item Due to the fact $\abs{\oz(0,x)} \le \alpha$ by  Corollary~\ref{cor:structure_opt_control} (a), the points $\oxi$ with $\oz(0,\oxi) = -\alpha$ are the minimizers and the points $\oxi$ with $\oz(0,\oxi) = \alpha$ are the maximizers of $\oz(0)$. Therefore, we have $\nabla \oz(0,\oxi) = 0$ and the corresponding Hessian matrices are positive semidefinite in the former and negative semidefinite in the later case. In addition we assume positive and negative definiteness respectively. This assumption corresponds to sufficient second order optimality conditions for minimizers and maximizers of $\oz(0)$.
\item Similar assumptions can be found in the literature, see~\cite{MerinoNeitzelTroeltzsch:2011, shapiroextremalvalue} in the context of semi-infinite programming and the notion of non-degeneracy in super-resolution~\cite{Duval2015}.
\end{itemize}
\end{rmrk}

Under the above assumption the optimal control $\oq$ consists of finitely many Diracs and has the form
\begin{equation}\label{eq:str_oq_ass}
\oq = \sum_{i=1}^K \ob_i \delta_{\oxi}
\end{equation}
with $\ob = \{\ob_i\}\in \R^K$, where $\ob_i >0$ for $\oxi$ with $\oz(0,\oxi) = -\alpha$ and $\ob_i <0$ for $\oxi$ with $\oz(0,\oxi) = \alpha$.

\subsection{Error estimates for the temporal error}
We first prove that under Assumption~\ref{ass:finite} the structure of the semidiscrete control $\oq_k$ is similar to that of $\oq$~\eqref{eq:str_oq_ass}. To this end we first show that Hessian matrix of the discrete adjoint state $\oz_k$ has the same definiteness properties as of the continuous adjoint state $\oz$ in the neighborhoods of the points $\oxi$.
\begin{lmm}\label{lemma:definitness_hessian_z_k}
Let Assumption~\ref{ass:finite} be fulfilled.  Let $\oq_k \in \M(\Omega)$ be the optimal solution of the semidiscrete problem~\eqref{eq:red_problem_k} with the corresponding state $\ou_k \in \Xk$ and the adjoint state $\oz_k \in \Xk$. Then there exist $\varepsilon >0$, $k_0>0$, and $\gamma>0$ such that
\[
\lambda_{\min}(\nabla^2 \oz_{k,0}^+)(x) \ge \gamma 
\]
for all $x \in B_\varepsilon(\oxi)$ and all $k \le k_0$ for $\oxi$ with $\oz(0,\oxi) = -\alpha$, where $\lambda_{\min}(\cdot)$ denotes the smallest eigenvalue of the corresponding matrix. Similarly,
\[
\lambda_{\min}(-\nabla^2 \oz_{k,0}^+)(x) \ge \gamma 
\]
for all $x \in B_\varepsilon(\oxi)$ and all $k \le k_0$ for $\oxi$ with $\oz(0,\oxi) = \alpha$.
\end{lmm}
\begin{proof}
We consider $\oxi$ with $\oz(0,\oxi) = -\alpha$. The Hessian matrix $\nabla^2 \oz(0,\oxi)$ is positive definite by Assumption~\ref{ass:finite}. Moreover $\oz(0) \in C^2(\Omega_0)$ by Lemma~\ref{lemma:higher_reg_z}. Therefore, there exists a neighborhood $B_\varepsilon(\oxi)$ such that $\nabla^2 \oz(0,x)$ is uniformly positive definite for $x \in B_\varepsilon(\oxi)$. It remains to prove that
\[
\norm{\oz(0)-\oz_{k,0}^+}_{C^2(\Omega_0)} \to 0 \quad \text{as }\; k \to 0.
\]
There holds
\[
\norm{\oz(0)-\oz_{k,0}^+}_{C^2(\Omega_0)} \le c \norm{\oz(0)-\oz_{k,0}^+}_{H^4(\Omega_0)} \le c \norm{\Delta(\oz(0)-\oz_{k,0}^+)}_{H^2(\Omega)} \le c
\ltwonorm{\Delta^2(\oz(0)-\oz_{k,0}^+)}
\]
by the embedding $H^4(\Omega_0) \hookrightarrow C^2(\Omega_0)$, the interior regularity result~\cite[Chapter 6.3,Theorem 2]{LCEvans_2010} and convexity of $\Omega$. To proceed we insert $\hat z_k \in X_k$ defined by~\eqref{eq:hat_z_k} leading to
\[
\norm{\oz(0)-\oz_{k,0}^+}_{C^2(\Omega_0)} \le c \ltwonorm{\Delta^2(\oz(0)-\hat z_{k,0}^+)} + c \ltwonorm{\Delta^2(\hat z_{k,0}^+-\oz_{k,0}^+)}.
\]
The first term is directly estimated by Lemma~\ref{lemma:higher_estimates_k} (below) with $j=2$ resulting in
\[
\ltwonorm{\Delta^2(\oz(0)-\hat z_{k,0}^+)} \le C k^{2r+1}
\]
and for the second term we have by the smoothing estimate  Lemma~\ref{lemma: higher smoothing}
\[
\ltwonorm{\Delta^2(\hat z_{k,0}^+-\oz_{k,0}^+)} \le C \ltwonorm{(\ou-\ou_k)(T)} \le C k^{r+\frac{1}{2}},
\]
where in the last step we used the first estimate from Theorem~\ref{theorem:general_estimate}. This completes the proof.
\end{proof}

\begin{lmm}\label{lemma:oq_finite_structure_k}
Let Assumption~\ref{ass:finite} be fulfilled.  Let $\oq_k \in \M(\Omega)$ be the optimal solution of the semidiscrete problem~\eqref{eq:red_problem_k} with the corresponding state $\ou_k \in \Xk$ and the adjoint state $\oz_k \in \Xk$. Then there is an $\varepsilon >0$ and $k_0>0$ such that the neighborhoods $B_\varepsilon(\oxi)$ are pairwise disjoint and for each $i$ and $k \le k_0$ there is a unique $\oxki \in B_\varepsilon(\oxi)$ such that
\[
\oz_{k,0}^+(\oxki) = \alpha \;\text{ if }\; \oz(0,\oxi) = \alpha
\]
and
\[
\oz_{k,0}^+(\oxki) = -\alpha \;\text{ if }\; \oz(0,\oxi) = -\alpha.
\]
Moreover there are no further points $x \in \Omega$ with $\oz_{k,0}^+(x) = \pm \alpha$ and the semidiscrete control has the structure
\[
\oq_k = \sum_{i=1}^K \ob_{k,i} \delta_{\oxki}
\]
with $\ob_k = \{\ob_{k,i}\} \in \R^K$, where $\ob_{k,i} > 0$ for $\oxki$ with $\oz_{k,0}^+(\oxki) = -\alpha$ and $\ob_{k,i} < 0$ for $\oxki$ with $\oz_{k,0}^+(\oxki) = \alpha$.
\end{lmm}
\begin{proof}
First we choose $\varepsilon >0$ such that the statement of Lemma~\ref{lemma:definitness_hessian_z_k} are fulfilled for all $i$  and the balls $ \bar B_\varepsilon(\oxi)$ are pairwise disjoint. Let $i$ be fixed with $\oz_{k,0}^+(\oxki) = -\alpha$. The case of  $\oz_{k,0}^+(\oxki) = \alpha$ is discussed in the same fashion. From Lemma~\ref{lemma:weak_star_conv} we have $\oq_k \overset{\ast}{\rightharpoonup} \oq$ in $\M(\Omega)$. We choose a smooth cut-off function $\omega$ with $\omega(x)=1$ for all $x \in B_{\varepsilon/2}(\oxi)$ and with $\supp \omega \subset B_\varepsilon(\oxi)$. From the weak star convergence we obtain
\[
\langle \oq_k, \omega \rangle \to \langle \oq,\omega \rangle = \ob_i \omega(\oxi) = \ob_i > 0.
\]
Therefore, there exists $k_0>0$ such that $\langle \oq_k, \omega \rangle>0$ for all $k \le k_0$, which proves that $\supp \oq_k \cap B_\varepsilon(\oxi)$ is not empty. The support condition for $\oq_k$ from Corollary~\ref{cor:opt_cond_k} implies the existence of at least one $\oxki \in B_\varepsilon(\oxi)$ with $\oz_{k,0}^+(\oxki) = -\alpha$.
By Lemma~\ref{lemma:definitness_hessian_z_k} $\oz_{k,0}^+$ is strictly convex on $B_\varepsilon(\oxi)$. This implies the uniqueness of the minimizer $\oxki$ in $B_\varepsilon(\oxi)$. In order to show that there are no further points $x$ with $\oz_{k,0}^+(x) = -\alpha$ in the complement of the union of all $B_\varepsilon(\oxi)$, it is sufficient to show that $\linfnorm{\oz(0)-\oz_{k,0}^+} \to 0$ for $k \to 0$. We have as in the proof of the previous lemma
\[
\linfnorm{\oz(0)-\oz_{k,0}^+} \le \linfnorm{\oz(0)-\hat z_{k,0}^+}+\linfnorm{\hat z_{k,0}^+-\oz_{k,0}^+},
\]
where $\hat z_k \in X_k$ is defined by~\eqref{eq:hat_z_k}. For the first term we obtain by Theorem~\ref{thm:smoothing_est_time}
\[
\linfnorm{\oz(0)-\hat z_{k,0}^+} \le C(T) k^{r+1} \ltwonorm{\bar u (T) -u_d}
\]
and for the second one
\[
\linfnorm{\hat z_{k,0}^+-\oz_{k,0}^+} \le \ltwonorm{(\ou-\ou_k)(T)} \le C k^{r+\frac{1}{2}},
\]
where in the last step we used the first estimate from Theorem~\ref{theorem:general_estimate}. This completes the proof.
\end{proof}

The main result of this section provides optimal order estimates for the error in the position $\oxi -\oxki$, the coefficients $\ob_i-\ob_{k,i}$ and improves the first estimate from Theorem~\ref{theorem:general_estimate} for  the state error $\ltwonorm{(\ou-\ou_k)(T)}$.

\begin{thrm}\label{th:improved_k_error}
Let $\oq$ be the solution of the problem of~\eqref{eq:red_problem} with the corresponding optimal state $\ou$ and the adjoint state $\oz$ and let Assumption~\ref{ass:finite} be fulfilled.  Let moreover $\oq_k \in \M(\Omega)$ be the optimal solution of the semidiscrete problem~\eqref{eq:red_problem_k} with the corresponding state $\ou_k \in \Xk$ and the adjoint state $\oz_k \in \Xk$. Then there exists $k_0>0$ such that for $k \le k_0$ there hold
\begin{itemize}
\item[(a)]
\[
\ltwonorm{(\ou-\ou_k)(T)} \le  C k^{2r+1},
\]
\item[(b)]
\[
\max_{1\le i \le K} \abs{\oxi-\oxki} \le C k^{2r+1}
\]
\item[(c)]
\[
\abs{\ob-\ob_k} \le C k^{2r+1}
\]
\item[(d)]
\[
\norm{\oq-\oq_k}_{\left(W^{1,\infty}(\Omega)\right)^*}\le C k^{2r+1}
\]
\end{itemize}
where $C=C(T,\alpha,u_d)$ depends on $T$, $\alpha$ and $\ltwonorm{u_d}$.
\end{thrm}
\begin{rmrk} \label{rem:tempestforKR}
From the equivalence relation in~\eqref{eqofwinfkr} we directly infer 
\begin{align*}
\norm{\oq-\oq_k}_{\text{KR}}  \le C k^{2r+1}.
\end{align*}
from statement (d) above.
\end{rmrk}
To prepare the proof of Theorem~\ref{th:improved_k_error} we first  estimate the error in the position $\oxi -\oxki$ and in the coefficients $\ob_i-\ob_{k,i}$ in terms of the state error $\ltwonorm{(\ou-\ou_k)(T)}$.

\begin{lmm}\label{lemma:xi-xki}
Let Assumption~\ref{ass:finite} be fulfilled. Then there exists $k_0>0$ such that for $k \le k_0$ there holds
\[
\abs{\oxi-\oxki} \le C k^{2r+1} + C \ltwonorm{(\ou-\ou_k)(T)}
\]
for all $i=1,2,\dots,K$, where $C = C(T,u_d)$ depends on $T$ and $\ltwonorm{u_d}$.
\end{lmm}
\begin{proof}
For a fixed $i \in \{1,2,\dots, K\}$ we assume without restriction that $\oz(\oxi)=-\alpha$ (the case $\oz(\oxi) = \alpha$ can be treated similarly). Then we have that $\oz_{k,0}^+(\oxki)=-\alpha$ by Lemma~\ref{lemma:oq_finite_structure_k}. The point $\oxi$ is a minimizer of $\oz(0)$ and the point $\oxki$ is a minimizer of $\oz_{k,0}^+$. Therefore, there holds
\[
\nabla \oz(0,\oxi) = 0 \quad \text{and} \quad \nabla \oz_{k,0}^+(\oxki) = 0.
\]
Due to the fact that $\na^2\oz(0)$ is uniformly positive definite on $B_\varepsilon(\oxi)$ and $\oxki \in B_\varepsilon(\oxi)$ we have
\[
\abs{\oxi-\oxki} \le  C \abs{\nabla \oz(0,\oxi)-\nabla \oz(0,\oxki)} = C\abs{\nabla \oz_{k,0}^+(\oxki)-\nabla \oz(0,\oxki)}
\]
and therefore
\[
\abs{\oxi-\oxki} \le C \linfnormn{\nabla(\oz(0)-\oz_{k,0}^+)} \le  C \|\nabla(\oz(0)-\hat z_{k,0}^+)\|_{L^\infty{(\Om_0)}} +  C \|\nabla(\hat z_{k,0}^+-\oz_{k,0}^+)\|_{L^\infty{(\Om_0)}},
\]
where $\hat z_{k,0}^+ \in \Xk$ is the solution of the intermediate discrete adjoint equation~\eqref{eq:hat_z_k} and $\Omega_0$ is an interior subdomain with
\[
\bigcup\limits_{1\le i \le K} B_\varepsilon(\oxi) \subset \Omega_0, \quad \bar \Omega_0 \subset \Omega.
\]
By Theorem~\ref{thm:smoothing_est_time_grad} we have
\[
\|\nabla(\oz(0)-\hat z_{k,0}^+)\|_{L^\infty{(\Om_0)}} \le C k^{2r+1} \ltwonorm{\ou(T)-u_d} \le  C k^{2r+1} \ltwonorm{u_d}
\]
and by the smoothing property from Lemma~\ref{lemma:semidiscrete_smoothing}
\[
\|\nabla(\hat z_{k,0}^+-\oz_{k,0}^+)\|_{L^\infty{(\Om_0)}} \le C \ltwonorm{(\ou-\ou_k)(T)}.
\]
This completes the proof.
\end{proof}

To proceed we introduce the operators $\widetilde G, \widetilde G^k, \widetilde G_k \colon \R^K \to L^2(\Omega)$ by
\[
\widetilde G(\beta) = S \left(\sum_{i=1}^K \beta_i \delta_{\oxi}\right), \; \widetilde G^k(\beta) = S \left(\sum_{i=1}^K \beta_i \delta_{\oxki}\right), \;  \widetilde G_k(\beta) = S_k \left(\sum_{i=1}^K \beta_i \delta_{\oxki}\right),
\]
where $S$ and $S_k$ are the continuous and the semidiscrete solution operators defined above. Moreover we restrict the codomains of these operators to the corresponding image sets and call the resulting operator $G$, $G^k$ and $G_k$ with
\[
G \colon \R^K \to \operatorname{Im}(\widetilde G)\subset L^2(\Omega), \quad G^k \colon \R^K \to \operatorname{Im}(\widetilde G^k)\subset L^2(\Omega),\quad  G_k \colon \R^K \to \operatorname{Im}(\widetilde G_k)\subset L^2(\Omega)
\]
and
\begin{equation}\label{eq:GandG_k}
G(\beta) = \widetilde G(\beta), \; G^k(\beta)= \widetilde G^k(\beta), \; G_k(\beta) = \widetilde G_k(\beta), \quad \text{for all }\; \beta \in \R^K.
\end{equation}
In the next lemma we estimate the errors between these operators.
\begin{lmm}\label{lemma:G-G_k}
Let Assumption~\ref{ass:finite} be fulfilled. There hold
\[
\ltwonorm{G(\beta)-G^k(\beta)} \le C \abs{\beta} \max_{1 \le i \le K} \abs{\oxi-\oxki}
\]
and
\[
\ltwonorm{G^k(\beta)-G_k(\beta)} \le C \abs{\beta} k^{2r+1}.
\]
\end{lmm}
\begin{proof}
For given $\beta \in \R^K$ we consider $q,q^k\in \M(\Omega)$ defined by
\[
q = \sum_{i=1}^K \beta_i \delta_{\oxi} \quad \text{and} \quad q^k = \sum_{i=1}^K \beta_i \delta_{\oxki}
\] 
as well as the corresponding states $u=u(q)$, $u^k = u(q^k)$ in the sense of Proposition~\ref{prop:state} and $u_k = u_k(q^k) \in \tXk$ in the sense of~\eqref{eq:state_k}. The second statement is then directly given by Lemma~\ref{lemma:est_state}, since
\[
\ltwonorm{G^k(\beta)-G_k(\beta)} = \ltwonorm{(u^k-u_k)(T)} \le C k^{2r+1} \mnorm{q^k} \le C \abs{\beta} k^{2r+1}.
\]
To prove the first statement we consider a dual problem for $y \in W(0,T)$ solving
\[
\begin{aligned}
-\pa_t y-\Delta y &= 0, &&\text{in} \quad (0,T)\times \Om,\;  \\
y                &= 0, &&\text{on}\quad  (0,T)\times\pa\Omega, \\
y(T)             &= (u-u^k)(T), &&\text{in}\quad \Omega
\end{aligned}
\]
and obtain
\[
\begin{aligned}
\ltwonorm{G(\beta)-G^k(\beta)}^2 &= \ltwonorm{(u-u^k)(T)}^2 = \langle q -q^k,y(0) \rangle\\
& = \sum_{i=1}^K \beta_i \left(y(0,\oxi) - y(0,\oxki)\right)\\
& \le C \abs{\beta} \|\nabla y(0)\|_{L^\infty{(\Om_0)}}  \max_{1 \le i \le K} \abs{\oxi-\oxki}\\
& \le C \abs{\beta} \ltwonorm{(u-u^k)(T)}  \max_{1 \le i \le K} \abs{\oxi-\oxki},
\end{aligned}
\]
where in the last step we used smoothing estimate from Lemma~\ref{lemma:cont_smoothing} for $y$. This completes the proof.
\end{proof}

\begin{lmm}\label{lemma:G_G_k_inverse}
Let Assumption~\ref{ass:finite} be fulfilled. The operators $G$, $G^k$ are bijective and there is a constant $c>0$ such that
\[
\abs{\beta} \le c \ltwonorm{G(\beta)} \quad \text{and} \quad \abs{\beta} \le c \ltwonorm{G^k(\beta)}
\]
hold for all $\beta \in \R^K$. Moreover, there is  $k_0>0$ such that $G_k$ is bijective and the estimate
\[
\abs{\beta} \le c \ltwonorm{G_k(\beta)}
\]
holds for all $\beta \in \R^K$ and all $k \le k_0$.
\end{lmm}
\begin{proof}
All three operators $G$, $G^k$ and $G_k$  are surjective by definition. We first argue the injectivity of $G$. Let $G(\beta)=0$ for some $\beta \in \R^K$. This means that for the solution $v$ of the heat equation with the initial condition given as the corresponding linear combination of Diracs, i.e. $v = S\left(\sum_{i=1}^K \beta_i \delta_{\oxi}\right)$, we have $v(T)=0$. By a similar argument as in proof of uniqueness of the optimal control $\oq$, cf.~\cite[Theorem 2.4]{CasasE_VexlerB_ZuazuaE_2015}, we obtain
\[
\sum_{i=1}^K \beta_i \delta_{\oxi} = 0.
\]
This results in $\beta = 0$ by the fact that the points $\oxi$ are pairwise disjoint. This provides the existence of an inverse mapping $G^{-1} \colon \operatorname{Im}(\widetilde G)\subset L^2(\Omega) \to \R^K$ and the estimate
\begin{equation}\label{eq:est_G_-1}
\abs{\beta} \le  \norm{G^{-1}}_{L^2(\Omega)\to \R^K} \ltwonorm{G(\beta)}
\end{equation}
holds. For the operator $G^k$ we can argue similarly. It remains to show that $G_k$ is bijective and $\norm{G_k^{-1}}_{L^2(\Omega)\to \R^K}$ is bounded independently of $k$. Let $\beta \in \R^K$ be arbitrary. There holds by~\eqref{eq:est_G_-1} and Lemma~\ref{lemma:G-G_k}
\[
\begin{aligned}
\abs{\beta} &\le c \ltwonorm{G(\beta)} \le c \ltwonorm{G_k(\beta)} + c \ltwonorm{G(\beta)-G^k(\beta)} + c \ltwonorm{G^k(\beta)-G_k(\beta)}\\
&\le c \ltwonorm{G_k(\beta)} + c \abs{\beta} \max_{1 \le i \le K} \abs{\oxi-\oxki} + c \abs{\beta} k^{2r+1}\\
& \le c\ltwonorm{G_k(\beta)} + c \abs{\beta} k^{r+\frac{1}{2}} + c \abs{\beta} k^{2r+1},
\end{aligned}
\]
where in the last step we used Lemma~\ref{lemma:xi-xki} and Theorem~\ref{theorem:general_estimate}. Choosing $k_0$ small enough we obtain
\[
\abs{\beta} \le c \ltwonorm{G_k(\beta)} + \frac{1}{2} \abs{\beta},
\]
which completes the proof.
\end{proof}

\begin{lmm}\label{lemma:ob-obk}
Let Assumption~\ref{ass:finite} be fulfilled. Then there exists $k_0>0$ such that for $k \le k_0$ there holds
\[
\abs{\ob-\ob_k} \le C k^{2r+1} + C \ltwonorm{(\ou-\ou_k)(T)},
\]
where $C = C(T,u_d)$ depends on $T$ and $\ltwonorm{u_d}$.
\end{lmm}
\begin{proof}
There holds by Lemma~\ref{lemma:G_G_k_inverse}
\[
\begin{aligned}
\abs{\ob-\ob_k} &\le c \ltwonorm{G(\ob)-G(\ob_k)}\\
& \le c \ltwonorm{G(\ob)-G_k(\ob_k)} + c\ltwonorm{G_k(\ob_k)-G^k(\ob_k)} + c\ltwonorm{G^k(\ob_k)-G(\ob_k)}.
\end{aligned}
\]
By the definition we have $G(\ob) = \ou(T)$ and $G_k(\ob_k) = \ou_k(T)$. Using Lemma~\ref{lemma:G-G_k} and Lemma~\ref{lemma:xi-xki} we obtain
\[
\begin{aligned}
\abs{\ob-\ob_k} &\le c \ltwonorm{(\ou-\ou_k)(T)} + c \abs{\ob_k} k^{2r+1} + c \abs{\beta} \max_{1 \le i \le K} \abs{\oxi-\oxki}\\
& \le c \ltwonorm{(\ou-\ou_k)(T)} + c \abs{\ob_k} k^{2r+1}.
\end{aligned}
\]
The fact that
\[
\abs{\ob_k} \le c \mnorm{\oq_k} \le c \ltwonorm{u_d}
\]
completes the proof.
\end{proof}

Previous lemmas allow us to obtain the corresponding estimate for a negative norm of $\oq-\oq_k$ in terms of  $\ltwonorm{(\ou-\ou_k)(T)}$.


\begin{lmm}\label{lemma:q-qk_dual}
Let Assumption~\ref{ass:finite} be fulfilled.  Let $\oq\in \M(\Omega)$ be the solution of~\eqref{eq:red_problem} and $\oq_k \in \M(\Omega)$ be the solution of the semidiscrete problem~\eqref{eq:red_problem_k}. Then there holds
\[
\norm{\oq-\oq_k}_{\left(W^{1,\infty}(\Omega)\right)^*} \le C k^{2r+1} + C \ltwonorm{(\ou-\ou_k)(T)}.
\]
\end{lmm}
\begin{proof}
Let $\varphi \in W^{1,\infty}(\Omega)$ with $\norm{\varphi}_{W^{1,\infty}(\Omega)} \le 1$. We have to estimate
\[
\begin{aligned}
\langle \oq-\oq_k,\varphi\rangle &= \sum_{i=1}^K \left(\ob_i \varphi(\oxi)-\ob_{k,i} \varphi(\oxki)\right)\\
&=\sum_{i=1}^K (\ob_i-\ob_{k,i}) \varphi(\oxi) + \sum_{i=1}^K \ob_{k,i} (\varphi(\oxi)-\varphi(\oxki))\\
&\le C \abs{\ob-\ob_k} \linfnorm{\varphi} + \abs{\ob_k} \max_{1 \le i \le K} \abs{\oxi-\oxki} \linfnorm{\nabla \varphi}.
\end{aligned}
\]
Using Lemma~\ref{lemma:xi-xki} and Lemma~\ref{lemma:ob-obk} we obtain
\[
\langle \oq-\oq_k,\varphi\rangle \le C k^{2r+1} + C \ltwonorm{(\ou-\ou_k)(T)},
\]
which completes the proof.
\end{proof}

We proceed with the proof of Theorem~\ref{th:improved_k_error}.
\begin{proof}
We start with the estimate~\eqref{eq:starting_estimate_k} from the proof of Theorem~\ref{theorem:general_estimate}, i.e.
\begin{equation}\label{eq:interm_est_improved_k}
\ltwonorm{(\ou-\ou_k)(T)}^2 \le 2 \langle \oq_k-\oq, \oz(0)- \hat z_{k,0}^+\rangle + \ltwonorm{(\ou-\hat u_k)(T)}^2,
\end{equation}
where $\hat u_k = u_k(\oq) \in \Xk$ and $\hat z_k \in \Xk$ is the solution of~\eqref{eq:hat_z_k}. The second term can be estimated as in the proof of Theorem~\ref{theorem:general_estimate} by the first estimate in Lemma~\ref{lemma:est_state} leading to
\[
\ltwonorm{(\ou-\hat u_k)(T)} \le Ck^{2r+1} \mnorm{\oq} \le Ck^{2r+1} \ltwonorm{u_d}^2.
\]
It remains to estimate the duality product from~\eqref{eq:interm_est_improved_k}. We have by Lemma~\ref{lemma:q-qk_dual}
\[
\begin{aligned}
\langle \oq_k-\oq, \oz(0)- \hat z_{k,0}^+\rangle &\le \norm{\oq-\oq_k}_{\left(W^{1,\infty}(\Omega)\right)^*} \norm{\oz(0)- \hat z_{k,0}^+}_{W^{1,\infty}(\Omega_0)}\\
\le & \left(C k^{2r+1} + C \ltwonorm{(\ou-\ou_k)(T)}\right) \norm{\oz(0)- \hat z_{k,0}^+}_{W^{1,\infty}(\Omega_0)},
\end{aligned}
\]
where we have used the fact that $\supp \oq, \supp \oq_k \subset \Omega_0$.
Using Theorem~\ref{thm:smoothing_est_time} and Theorem~\ref{thm:smoothing_est_time_grad} we have
\[
\norm{\oz(0)- \hat z_{k,0}^+}_{W^{1,\infty}(\Omega)}  \le C k^{2r+1} \ltwonorm{\ou(T)-u_d} \le C k^{2r+1} \ltwonorm{u_d}.
\]
Putting all terms in~\eqref{eq:interm_est_improved_k} together we get
\[
\ltwonorm{(\ou-\ou_k)(T)}^2 \le C \left(C k^{2r+1} + C \ltwonorm{(\ou-\ou_k)(T)}\right) k^{2r+1} + C k^{4r+2}.
\]
Absorbing $\ltwonorm{(\ou-\ou_k)(T)}$ in the left-hand side we obtain the estimate (a) in Theorem~\ref{th:improved_k_error}. The estimates (b), (c), and (d) are obtained from (a) using Lemma~\ref{lemma:xi-xki}, Lemma~\ref{lemma:ob-obk} as well as Lemma~\ref{lemma:q-qk_dual}.
\end{proof}
\subsection{Error estimates for the spatial error}

For the semidiscretization we have shown (see Lemma~\ref{lemma:oq_finite_structure_k}) that the number of support points of the semidiscrete control $\oq_k \in \M(\Omega)$ is the same as on the continuous level in Assumption~\ref{ass:finite}. For the fully discrete control $\oq_{kh} \in \M_h$ the situation is different. We will show (see next lemma) that in the neighborhood of each support point $\oxki$ of $\oq_k$ there is at least one support point of $\oq_{kh}$, but there could be more than one such point. This phenomena is also observed in our numerical experiments.

\begin{lmm}\label{lemma:oq_finite_structure_h}
Let Assumption~\ref{ass:finite} be fulfilled.  Let $\oq_k \in \M(\Omega)$ be the optimal solution of the semidiscrete problem~\eqref{eq:red_problem_k} with the corresponding state $\ou_k \in \Xk$ and the adjoint state $\oz_k \in \Xk$. Let $\oq_{kh} \in \M_h$ be an optimal solution of the fully discrete problem~\eqref{eq:red_problem_kh} with the corresponding state $\ou_{kh} \in \Xkh$ and the adjoint state $\oz_{kh} \in \Xkh$. Then there is $k_0>0$ such that for any fixed $k<k_0$ the following holds. There is an $\varepsilon >0$ and  $h_0>0$ such that the neighborhoods $B_\varepsilon(\oxki)$ are pairwise disjoint and for each $i$ and $h \le h_0$ there is at least one $\oxkhij \in B_\varepsilon(\oxki)\cap {\mathcal N}_h$ such that
\[
\oz_{kh,0}^+(\oxkhij) = \alpha \;\text{ if }\;  \oz_{k,0}^+(\oxki) = \alpha
\]
and
\[
\oz_{kh,0}^+(\oxkhij) = -\alpha \;\text{ if }\; \oz_{k,0}^+(\oxki) = -\alpha.
\]
Moreover there are no points $x \in \Omega \setminus \bigcup_i B_\varepsilon(\oxki)$ with $\oz_{kh,0}^+(x) = \pm \alpha$.
\end{lmm}
\begin{proof}
The proof is similar to the proof of Lemma~\ref{lemma:oq_finite_structure_k}.
\end{proof}

Under the conditions of Lemma~\ref{lemma:oq_finite_structure_h} a fully discrete control $\oq_{kh}$ consists of groups of Dirac functionals for each single Dirac $\delta_{\oxki}$ on the semidiscrete level. This means, that $\oq_{kh}$ is given as
\[
\oq_{kh} = \sum_{i=1}^K \sum_{j=1}^{n_i} \ob_{kh,ij} \delta_{\oxkhij},
\]
where $n_i \in {\mathbb N}$ describes the cardinality of $\supp \oq_{kh}|_{B_\varepsilon(\oxki)}$. The cardinality of $\supp \oq_{kh}$ is then $K_h = \sum_{i=1}^K n_i \ge K$. In order to compare the vector of coefficients $\ob_{kh} =\{\ob_{kh,ij}\} \in \R^{K_h}$ with the vector $\ob_k \in \R^K$ on the semidiscrete level, we define $\hat \beta_{kh} \in \R^K$ by
\[
\hat \beta_{kh,i} = \sum_{j=1}^{n_i} \ob_{kh,ij}.
\]

The next theorem is the main result of this section.
\begin{thrm}\label{th:improved_h_error}
Under the conditions of Lemma~\ref{lemma:oq_finite_structure_h} there holds
\begin{itemize}
\item[(a)]
\[
\ltwonorm{{\ou_k-\ou_{kh}}} \le C\ell_{kh}^\frac{1}{2} h,
\]
\item[(b)] 
\[
\abs{\oxki-\oxkhij} \le C\ell_{kh}^\frac{1}{2} h
\]
for all $1 \le i \le K$ and $1 \le j \le n_i$.
\item[(c)]
\[
\abs{\ob_k-\hat \beta_{kh}} \le C\ell_{kh}^\frac{1}{2} h,
\]
\item[(d)]
\[
\norm{\oq_k-\oq_{kh}}_{\left(W^{1,\infty}(\Omega)\right)^*}\le C\ell_{kh}^{\frac{1}{2}}h.
\]
\end{itemize}
where $\ell_{kh} = \lk + \lh$.
\end{thrm}
\begin{rmrk}
As in Remark~\ref{rem:tempestforKR}, we directly conclude the a priori estimate
\begin{align*}
\norm{\oq_k-\oq_{kh}}_{\text{KR}} \leq C \ell^{\frac{1}{2}}_{kh} h
\end{align*}
from statement (d) in Theorem~\ref{th:improved_h_error}.
\end{rmrk}

To prove Theorem~\ref{th:improved_h_error} we start with the lemma providing a sub-optimal estimate for the distance between the support points of $\oq_k$ and $\oq_{kh}$.

\begin{lmm}\label{lemma:error_x_h_suboptimal}
Under the conditions of Lemma~\ref{lemma:oq_finite_structure_h} there is a constant $C>0$ such that for each $\oxkhi \in B_\varepsilon(\oxki)$ with  $\oz_{kh,0}^+(\oxkhi) = \pm \alpha$ there holds
\[
\abs{\oxki-\oxkhi} \le C  \ell_{kh}^\frac{1}{4} h^\frac{1}{2},
\]
where $\ell_{kh} = \lk + \lh$.
\end{lmm}
\begin{proof}
We consider the Taylor expansion with an appropriate $\xi \in (\oxki,\oxkhi)$
\[
\begin{aligned}
\oz_{k,0}^+(\oxkhi) &=  \oz_{k,0}^+(\oxki) + \nabla  \oz_{k,0}^+(\oxki)^T (\oxki-\oxkhi) + \frac{1}{2} (\oxki-\oxkhi)^T \nabla^2 \oz_{k,0}^+(\xi) (\oxki-\oxkhi) \\
&= \oz_{kh,0}^+(\oxkhi)  + \frac{1}{2} (\oxki-\oxkhi)^T \nabla^2 \oz_{k,0}^+(\xi) (\oxki-\oxkhi),
\end{aligned}
\]
where we used that $\oz_{k,0}^+(\oxki) = \oz_{kh,0}^+(\oxkhi)$ by the previous lemma and $\nabla  \oz_{k,0}^+(\oxki) = 0$ by the optimality of $\oxki$ for $\oz_{k,0}^+$, see Corollary~\ref{cor:opt_cond_k} and Lemma~\ref{lemma:oq_finite_structure_k}. Using uniform definiteness of $ \nabla^2 \oz_{k,0}^+$, see Lemma~\ref{lemma:definitness_hessian_z_k}, we obtain
\[
\gamma \abs{\oxki-\oxkhi}^2 \le \linfnormn{\oz_{k,0}^+-\oz_{kh,0}^+},
\]
where $\Omega_0$ is an interior subdomain with
\[
\supp \oq_k \cup \supp \oq_{kh} \subset \Omega_0, \qquad \bar \Omega_0 \subset \Omega,
\]
see Corollary~\ref{cor:opt_cond_k} and  Corollary~\ref{cor:opt_cond_kh}.
To proceed we introduce an intermediate discrete adjoint state $\hat z_{kh} \in \Xkh$ fulfilling
\begin{equation}\label{eq:hat_z_kh}
B(\varphi_{kh},\hat z_{kh}) = (\ou_k(T)-u_d,\varphi_{kh}(T)) \quad \text{for all }\; \varphi_{kh}\in \Xkh.
\end{equation}
We obtain
\[
\gamma \abs{\oxki-\oxkhi}^2 \le \linfnormn{\oz_{k,0}^+-\hat z_{kh,0}^+} +  \linfnormn{\hat z_{kh,0}^+-\oz_{kh,0}^+}.
\]
The first term is estimated by Theorem~\ref{thm:smoothing_est_space} leading to
\[
\linfnormn{\oz_{k,0}^+-\hat z_{kh,0}^+} \le C \ell_{kh} h^2
\]
and the second term by the smoothing property from Corollary \ref{cor: discrete smoothing in Linfty} and by the second estimate from Theorem~\ref{theorem:general_estimate}
\[
\linfnormn{\hat z_{kh,0}^+-\oz_{kh,0}^+} \le C \ltwonorm{(\ou_k-\ou_{kh})(T)} \le C \ell_{kh}^\frac{1}{2} h.
\]
This completes the proof.
\end{proof}

For the proof of Theorem~\ref{th:improved_h_error} we introduce a further intermediate adjoint state $\tilde z_k \in \Xk$ defined by
\begin{equation}\label{eq:tilde_z_k}
B(\varphi_k,\tilde z_k) = (\ou_{kh}(T)-u_d,\varphi_k(T)) \quad \text{for all }\; \varphi_k\in \Xk.
\end{equation}

\begin{lmm}\label{lemma:tilde_xki}
Let the conditions of Lemma~\ref{lemma:oq_finite_structure_h} be fulfilled and let $\tilde z_k \in \Xk$ be defined by~\eqref{eq:tilde_z_k}. Then for each $\oxki$ with $\oz_{k,0}^+(\oxki)=-\alpha$ there is a  minimizer $\tilde x_{k,i} \in B_\varepsilon(\oxki)$ of $\tilde z_{k,0}^+$ and for each $\oxki$ with $\oz_{k,0}^+(\oxki)=\alpha$ there is a  maximizer $\tilde x_{k,i} \in B_\varepsilon(\oxki)$ of $\tilde z_{k,0}^+$.  Moreover there holds
\[
\abs{\oxki - \tilde x_{k,i}} \le C\ell_{kh}^\frac{1}{2}h.
\]
\end{lmm}
\begin{proof}
Without restriction we assume that $\varepsilon>0$ and $k_0>0$ from Lemma~\ref{lemma:oq_finite_structure_h} are chosen such that the statement of Lemma~\ref{lemma:definitness_hessian_z_k} holds. We fix $i$ with  $\oz_{k,0}^+(\oxki)=-\alpha$ and introduce two functions $F,F_h \colon B_\varepsilon(\oxki) \to \R^N$ by
\[
F(x) = \nabla \oz_{k,0}^+(x) \quad \text{and} \quad F_h(x) = \nabla \tilde z_{k,0}^+(x).
\]
There holds by the optimality of $\oxki$ for $\oz_{k,0}^+(x)$ that $F(\oxki) = 0$ and by Lemma~\ref{lemma:definitness_hessian_z_k} that $F'(\oxki) = \nabla^2\oz_{k,0}^+(\oxki) $ is positive definite. Moreover we have
\[
\abs{F(\oxki)-F_h(\oxki)} \le \|\nabla(\oz_{k,0}^+ - \tilde z_{k,0}^+)\|_{L^\infty(\Om_0)} \le C \ltwonorm{(\ou_k-\ou_{kh})(T)} \le C\ell_{kh}^\frac{1}{2} h
\]
and
\begin{equation}\label{eq:Fp-Fph}
\abs{F'(\oxki)-F'_h(\oxki)} \le \norm{\oz_{k,0}^+ - \tilde z_{k,0}^+}_{C^2(\Omega_0)} \le C \ltwonorm{(\ou_k-\ou_{kh})(T)} \le C\ell_{kh}^\frac{1}{2} h
\end{equation}
by the smoothing property from Lemma~\ref{lemma:semidiscrete_smoothing} and the estimate from Theorem~\ref{theorem:general_estimate}. In addition $F'_h$ is Lipschitz continuous on $B_\varepsilon(\oxki)$ with the Lipschitz constant
\[
L = \norm{\tilde z_{k,0}^+}_{C^3(\Omega_0)} \le C \ltwonorm{\ou_{kh}(T)-u_d} \le C \ltwonorm{u_d},
\]
where we have used interior estimate as in the proof of Lemma~\ref{lemma:cont_smoothing}.
In this setting we can apply~\cite[Theorem 3.1]{RannacherVexler:2005} to get the existence of $\tilde x_{k,i} \in B_\varepsilon(\oxki)$ (for $h < h_0$) with $F_h(\tilde x_{k,i})=0$ and a positive definite $F'_h(\tilde x_{k,i})$, such that
\[
\abs{\oxki - \tilde x_{k,i}} \le C \abs{F(\oxki)-F_h(\oxki)} \le C\ell_{kh}^\frac{1}{2}h.
\]
This completes the proof.
\end{proof}

In the next lemma we improve the estimate from Lemma~\ref{lemma:error_x_h_suboptimal}.
\begin{lmm}\label{lemma:tilde_xki_step2}
Let the conditions of Lemma~\ref{lemma:tilde_xki} be fulfilled. Then there holds for all $1 \le i \le K$ and all $1 \le j \le n_i$
\[
\abs{\tilde x_{k,i} - \oxkhij} \le  C\ell_{kh}^\frac{1}{2}h.
\]
\end{lmm}
\begin{proof}
We fix an $i$ with  $\oz_{k,0}^+(\oxki)=-\alpha$ and an $1\le j\le n_i$. For $\tilde z_k \in \Xk$ defined in~\eqref{eq:tilde_z_k}  we observe
\[
\linfnorm{\tilde z_{k,0}^+-\oz_{kh,0}^+} \le C\ell_{kh} h^2 \ltwonorm{\ou_{kh}(T)-u_d} \le C\ell_{kh} h^2
\]
by Theorem~\ref{thm:smoothing_est_space}. Due to Corollary~\ref{cor:opt_cond_k} we have $\oz_{k,0}^+(x) \ge -\alpha$ for all $x \in \bar \Omega$ and therefore
\[
\tilde z_{k,0}^+(\tilde x_{k,i}) \ge -\alpha- C\ell_{kh}h^2 = \oz_{kh,0}^+(\oxkhij) - C\ell_{kh}h^2.
\]
Using the Taylor expansion and the fact that $\nabla \tilde z_{k,0}^+(\tilde x_{k,i}) = 0$ we get with some $\xi \in B_\varepsilon(\oxki)$
\[
\begin{aligned}
\tilde z_{k,0}^+(\oxkhij) &= \tilde z_{k,0}^+(\tilde x_{k,i}) + \frac{1}{2} (\tilde x_{k,i} - \oxkhij)^T \nabla^2 \tilde z_{k,0}(\xi) (\tilde x_{k,i} - \oxkhij)\\
&\ge \oz_{kh,0}^+(\oxkhij) - C\ell_{kh}h^2 + \frac{\gamma}{2} \abs{\tilde x_{k,i} - \oxkhij}^2,
\end{aligned}
\]
where we have used that $\nabla^2 \tilde z_{k,0}$ is uniformly positive definite on $B_\varepsilon(\oxki)$ by the positive definiteness of $\nabla^2 \oz_{k,0}$ and~\eqref{eq:Fp-Fph}. This results in
\[
\begin{aligned}
\abs{\tilde x_{k,i} - \oxkhij}^2 &\le \abs{\tilde z_{k,0}^+(\oxkhij)-\oz_{kh,0}^+(\oxkhij)} + C\ell_{kh}h^2 \\
&\le \linfnormn{\tilde z_{k,0}^+-\oz_{kh,0}^+}+ C\ell_{kh}h^2\le C\ell_{kh}h^2. 
\end{aligned}
\]
This completes the proof.
\end{proof}

We proceed with the proof of Theorem~\ref{th:improved_h_error}.
\begin{proof}
The first statement is already shown in Theorem~\ref{theorem:general_estimate}. The second statement follows directly from Lemma~\ref{lemma:tilde_xki} and Lemma~\ref{lemma:tilde_xki_step2} by the triangle inequality. It remains to prove statements (c) and (d). To this end we will use the operator $G_k$ introduced in~\eqref{eq:GandG_k}. Similarly we introduce the operator $G_k^h \colon \R^K \to L^2(\Omega)$ by
\[
G_k^h(\beta) = S_{kh} \left(\sum_{i=1}^K \beta_i \delta_{\oxki} \right).
\]
Without restriction we assume that $k_0>0$ from Lemma~\ref{lemma:oq_finite_structure_h} is chosen such that the statement of Lemma~\ref{lemma:G_G_k_inverse} holds. Then we obtain similarly to the proof of Lemma~\ref{lemma:ob-obk} using $G_k(\ob_k)=\ou_k(T)$
\begin{equation}\label{eq:mid_ed_proof_th_impr_h}
\begin{aligned}
\abs{\ob_k-\hat \beta_{kh}} &\le C \ltwonorm{G_k(\ob_k)-G_k(\hat \beta_{kh})}\\
& \le C\ltwonorm{(\ou_k-\ou_{kh})(T)} + C\ltwonorm{\ou_{kh}(T)-G_k^h(\hat \beta_{kh})} + C\ltwonorm{G_k^h(\hat \beta_{kh})-G_k(\hat \beta_{kh})}.
\end{aligned}
\end{equation}
The first term is estimated by Theorem~\ref{theorem:general_estimate}
\[
\ltwonorm{(\ou_k-\ou_{kh})(T)} \le C \ell_{kh}^{\frac{1}{2}} h,
\]
the last term is estimated by Lemma~\ref{lemma:est_state} leading to
\[
\ltwonorm{G_k^h(\hat \beta_{kh})-G_k(\hat \beta_{kh})} \le C \ell_{kh} h^2 \abs{\hat \beta_{kh}} \le  C \ell_{kh} h^2.
\]
To estimate the second term in~\eqref{eq:mid_ed_proof_th_impr_h} we observe that
\[
\ou_{kh}(T)-G_k^h(\hat \beta_{kh}) = S_{kh} \left(\sum_{i=1}^K \Bigl(\sum_{j=1}^{n_i} \ob_{kh,ij} \delta_{\oxkhij} - \hat \beta_{kh,i} \delta_{\oxki} \Bigr) \right).
\]
For a given $\psi \in W^{1,\infty}(\Omega)$ we get for the inner difference using $\hat \beta_{kh,i} =\sum_{j=1}^{n_i} \ob_{kh,ij}$
\[
\left\langle\sum_{j=1}^{n_i} \ob_{kh,ij} \delta_{\oxkhij} - \hat \beta_{kh,i} \delta_{\oxki},\psi \right\rangle = \sum_{j=1}^{n_i}  \ob_{kh,ij} (\psi(\oxkhij)-\psi(\oxki))
\le C \linfnorm{\nabla \psi} \max_{1\le j\le n_i} \abs{\oxkhij-\oxki}.
\]
Then by a duality argument as in the proof of Lemma~\ref{lemma:G-G_k} we obtain
\[
\ltwonorm{\ou_{kh}(T)-G_k^h(\hat \beta_{kh})} \le C \max_{1 \le i \le K} \max_{1\le j\le n_i} \abs{\oxkhij-\oxki}
\]
resulting in
\[
\ltwonorm{\ou_{kh}(T)-G_k^h(\hat \beta_{kh})} \le C \ell_{kh}^\frac{1}{2}h
\]
by the statement (b). Inserting this into~\eqref{eq:mid_ed_proof_th_impr_h} completes the proof of the statement (c). To prove statement (d) let $\varphi \in W^{1,\infty}(\Omega)$ with $\norm{\varphi}_{W^{1,\infty}(\Omega)} \le 1$ be given. We estimate
\[
\begin{aligned}
\langle \oq_k-\oq_{kh},\varphi\rangle &= \sum_{i=1}^K \left(\ob_{k,i} \varphi(\oxki)- \sum_{j=1}^{n_i} \ob_{kh,ij}  \varphi(\oxkhij)\right)\\
&=\sum_{i=1}^K \left( (\ob_{k,i}-\hat \beta_{kh,i}) \varphi(\oxki) + \sum_{j=1}^{n_i}  \ob_{kh,ij} (\varphi(\oxki)-\varphi(\oxkhij))\right )\\
&\le C \left(\abs{\ob_k-\hat \beta_{kh}} \linfnorm{\varphi} + \abs{\hat \beta_{kh}} \max_{\substack{1 \le i \le K\\1 \leq j \leq n_i}} \abs{\oxki-\oxkhi} \linfnorm{\nabla \varphi}\right).
\end{aligned}
\]
Using statements~$(b)$ and~$(c)$ from Theorem~\ref{th:improved_h_error} as well as the boundedness of~$\abs{\hat \beta_{kh}}$ we obtain
\[
\langle \oq-\oq_k,\varphi\rangle \le C \ell^{\frac{1}{2}}_{kh} h,
\]
which completes the proof.
\end{proof}

\section{Proof of smoothing error estimates in time}\label{sec:proof_thm_smoothing_time}
In this section we prove Theorem~\ref{thm:smoothing_est_time} and Theorem~\ref{thm:smoothing_est_time_grad}. First we establish the following result.

\begin{lmm} \label{lemma: negative smooth superconvergence}
Let $v_0\in L^2(\Omega)$, let $v$ and $v_k$ satisfy \eqref{eq:eq_aux} and \eqref{eq: semidiscrete heat}. Then for $ l=0,1,\dots, r$, there exists a constant $C$ independent of $k$ and $T$ such that
$$
\|(-\Delta)^{-2l-1}(v-v_k)(T)\|_{L^2(\Om)}\le C k^{2l+1}\|v_0\|_{L^2(\Om)}.
$$
\end{lmm}
\begin{proof}
For any  $l\in \{0,1,\dots, r\}$, let $y$ be the solution to the following backward problem
$$
\begin{aligned}
-\pa_t y-\Delta y &= 0, &&\text{in} \quad (0,T)\times \Om,\;  \\
y                &= 0, &&\text{on}\quad  (0,T)\times\pa\Omega, \\
y(T)             &= (-\Delta)^{-2l-1}(v-v_k)(T,x), &&\text{in}\quad \Omega,
\end{aligned}
$$
and $y_k$ be its dG($r$) approximation, i.e.
$$
B(\psi_k,y_k) = ((-\Delta)^{-2l-1}(v-v_k)(T),\psi_k(T)), \quad \text{for all }\; \psi\in \Xk.
$$
Using the orthogonality conditions \eqref{eq: orthogonality semidiscrete} and the dual representation of the bilinear form \eqref{eq:B_dual}
$$
\begin{aligned}
\|(-\Delta)^{-2l-1}&(v-v_k)(T)\|^2_{L^2(\Om)} = B((-\Delta)^{-2l-1}(v-v_k),y)\\
& = B(v-v_k,(-\Delta)^{-2l-1}y)\\
& = B(v-v_k,(-\Delta)^{-2l-1}(y-y_k))\\
& = B(v-\pi_kv,(-\Delta)^{-2l-1}(y-y_k))\\
& = -\sum_{m=1}^M(v-\pi_kv,\pa_t(-\Delta)^{-2l-1}(y-y_k))_{I_m\times\Om}+(\na (v-\pi_kv),\na(-\Delta)^{-2l-1}(y-y_k))_{I\times\Om}\\
&  = -\sum_{m=1}^M(\na(v-\pi_kv),\na(-\Delta)^{-2l-2}\pa_t(y-y_k))_{I_m\times\Om}+(\na (v-\pi_kv),\na(-\Delta)^{-2l-1}(y-y_k))_{I\times\Om}\\
&:=J_1+J_2,
\end{aligned}
$$
where $\pi_k$ is the projection defined in \eqref{eq: projection pi_k}. Note, the jump terms vanish by the definition of $\pi_k$.
We set $\eta := y -\pi_k y$ and $\xi_k=\pi_ky-y_k$. 
Using the approximation and the standard energy estimate we have
\begin{equation}\label{eq: estimate for eta}
\begin{aligned}
\|\na(-\Delta)^{-l-1}\eta\|_{L^2(I\times\Om)}\le Ck^{l+1}\|\na(-\Delta)^{-l-1}\pa_t^{l+1} y\|_{L^2(I\times\Om)}=Ck^{l+1}\|\na y\|_{L^2(I\times\Om)}\le C k^{l+1}\|y(T)\|_{L^2(\Om)}.
\end{aligned}
\end{equation}
Using the properties of the bilinear form \eqref{eq: bilinear form B}, we have
$$
\begin{aligned}
\|\na(-\Delta)^{-l-1}\xi_k\|_{L^2(I\times\Om)}^2 &\le B((-\Delta)^{-l-1}\xi_k,(-\Delta)^{-l-1}\xi_k) = B((-\Delta)^{-2l-2}\xi_k,\xi_k)\\
&=-B((-\Delta)^{-2l-2}\xi_k,\eta)=-B((-\Delta)^{-l-1}\xi_k,(-\Delta)^{-l-1}\eta)\\
& = -(\na(-\Delta)^{-l-1}\xi_k,\na(-\Delta)^{-l-1}\eta)_{I\times\Om}\\
& \le\|\na(-\Delta)^{-l-1}\xi_k\|_{L^2(I\times\Om)} \|\na(-\Delta)^{-l-1}\eta\|_{L^2(I\times\Om)}. 
\end{aligned}
$$
Canceling and using \eqref{eq: estimate for eta} we obtain
\begin{equation}\label{eq: estimate for xi_k}
\|\na(-\Delta)^{-l-1}\xi_k\|_{L^2(I\times\Om)}\le C k^{l+1}\|y(T)\|_{L^2(\Om)}.
\end{equation}
Combining \eqref{eq: estimate for eta} and \eqref{eq: estimate for xi_k} we also have
\begin{equation}\label{eq: estimate for y-y_k}
\|\na(-\Delta)^{-l-1}(y-y_k)\|_{L^2(I\times\Om)}\le C k^{l+1}\|y(T)\|_{L^2(\Om)}.
\end{equation}
Next we estimate $\|\na (-\Delta)^{-l-1}\pa_t( y-y_k)\|_{L^2(I_m\times\Om)}$.  By the triangle inequality, inverse inequality and \eqref{eq: estimate for eta}, \eqref{eq: estimate for xi_k}, and \eqref{eq: estimate for y-y_k} we obtain
\begin{equation}
\begin{aligned}
\|\na &(-\Delta)^{-l-1}\pa_t( y-y_k)\|_{L^2(I_m\times\Om)}\\
&\le \|\na (-\Delta)^{-l-1}\pa_t(\pi_k y-y_k)\|_{L^2(I_m\times\Om)}+\|\na (-\Delta)^{-l-1}\pa_t( y-\pi_k  y)\|_{L^2(I_m\times\Om)}\\
& \le Ck^{-1}\|\na (-\Delta)^{-l-1}(\pi_k y-y_k)\|_{L^2(I_m\times\Om)}+Ck^{l}\|\na (-\Delta)^{-l-1}\pa^{l+1}_t y\|_{L^2(I_m\times\Om)}\\
& \le Ck^{-1}\left(\|\na (-\Delta)^{-l-1}(\pi_k y-y)\|_{L^2(I_m\times\Om)}+\|\na (-\Delta)^{-l-1}(y_k -y)\|_{L^2(I_m\times\Om)}\right)+Ck^{l}\|\na y\|_{L^2(I_m\times\Om)}\\
& \le Ck^{l}\|\na (-\Delta)^{-l-1}\pa_t^{l+1} y\|_{L^2(I_m\times\Om)}+Ck^{l}\|\na y\|_{L^2(I_m\times\Om)}\\
&\le Ck^{l}\|\na y\|_{L^2(I_m\times\Om)}.
\end{aligned}
\end{equation}
This allows to estimate $J_1$ as follows
$$
\begin{aligned}
J_1 &\le \sum_{m=1}^M\|\na (-\Delta)^{-l-1}( v-\pi_kv)\|_{L^2(I_m\times\Om)}\|\na (-\Delta)^{-l-1}\pa_t( y-y_k)\|_{L^2(I_m\times\Om)}\\
& \le  Ck^{2l+1} \sum_{m=1}^M\|\na (-\Delta)^{-l-1}\pa^{l+1}_t v\|_{L^2(I_m\times\Om)}\|\na (-\Delta)^{-l-1}\pa_t^{l+1} y\|_{L^2(I_m\times\Om)}\\
& =  Ck^{2l+1} \sum_{m=1}^M\|\na  v\|_{L^2(I_m\times\Om)}\|\na  y\|_{L^2(I_m\times\Om)}\\
& =  Ck^{2l+1} \|  v_0\|_{L^2(\Om)}\|(-\Delta)^{-2l-1}(v-v_k)(T)\|_{L^2(\Om)}.
\end{aligned}
$$
Similarly, using approximation and \eqref{eq: estimate for y-y_k} we obtain for $J_2$
$$
\begin{aligned}
J_2 &\le \|\na (-\Delta)^{-l-1}( v-\pi_kv)\|_{L^2(I\times\Om)}\|\na (-\Delta)^{-l}( y-y_k)\|_{L^2(I\times\Om)}\\
& \le  Ck^{2l+1} \|\na (-\Delta)^{-l-1}\pa_t^{l+1} v\|_{L^2(I\times\Om)}\| y(T)\|_{L^2(\Om)}\\
& \le   Ck^{2l+1} \|  v_0\|_{L^2(\Om)}\|(-\Delta)^{-2l-1}(v-v_k)(T)\|_{L^2(\Om)}.
\end{aligned}
$$
Combining the estimates for $J_1$ and $J_2$ and canceling $\|(-\Delta)^{-2l-1}(v-v_k)(T)\|_{L^2(\Om)}$, we obtain the lemma. 
\end{proof}
Now we show the next result.

\begin{lmm}\label{lemma:higher_estimates_k}
Let $v_0\in L^2(\Omega)$, let $v$ and $v_k$ satisfy \eqref{eq:eq_aux} and \eqref{eq: semidiscrete heat}.  Then for  $j \in \N_0$ provided  $k \le \frac{T}{2r+j+2}$ and $M>2r+j+2$, there exists a constant $C(T)$ independent of $k$  such that
\[
\ltwonorm{(-\Delta)^j(v-v_k)(T)} \le C(T) k^{2r+1}\ltwonorm{v_0},
\]
where $C(T)\sim T^{-2r-j-1}$.
\end{lmm}
\begin{proof}
For any  $j\in \N_0$, let $y$ and $y_k$ be the solutions to the continuous and to the semidiscrete dual problems  with $y_k(T) = y(T)=(-\Delta)^j (v-v_{k})(T)$, i.e. $y \in H^1(I; L^2(\Om))\cap L^2(I; H^1_0(\Om))$ solving
\begin{equation}\label{eq: dual continuous with yT=Delta v-v_k}
\begin{aligned}
-\pa_t y-\Delta y &= 0, &&\text{in} \quad (0,T)\times \Om,\;  \\
y                &= 0, &&\text{on}\quad  (0,T)\times\pa\Omega, \\
y(T)             &=(-\Delta)^j (v-v_{k})(T) , &&\text{in}\quad \Omega
\end{aligned}
\end{equation}
and $y_k \in \Xk$ satisfying 
\begin{equation}\label{eq: dual semidiscrete with yT=Delta v-v_k}
B(\varphi_{k},y_{k}) = (\varphi_{k}(T),(-\Delta)^j(v-v_{k})(T)), \quad \text{for all }\; \varphi_{k}\in \Xk.
\end{equation}
We choose  $\tilde{m}$ such that $\frac{T}{2}\in I_{\tilde{m}}$ and define $\tilde{v}:=\chi_{(t_{\tilde{m}},T]}v$ as well as $\tilde{v}_{k}=\chi_{(t_{\tilde{m}},T]}v_{k}$, i.e.  $\tilde{v}$ and $\tilde{v}_{k}$ are zero on $I_1\cup \cdots \cup I_{\tilde{m}}$ and $\tilde{v}=v$ and $\tilde{v}_{k}=v_{k}$ on the remaining time intervals.
Then we test~\eqref{eq: dual continuous with yT=Delta v-v_k} with $\varphi=(-\Delta)^j\tilde{v}$  and  
choose  $\varphi_k=(-\Delta)^j \tilde{v}_k$ in \eqref{eq: dual semidiscrete with yT=Delta v-v_k}. Using~\eqref{eq: Bilinear tilde}, we have
 $$
  \begin{aligned}
 \|(-\Delta)^j (v-v_{k})(T)&\|^2_{L^2(\Omega)}=B((-\Delta)^j\tilde{v},y)-B((-\Delta)^j\tilde{v}_k,y_k)\\
 =&B(\tilde{v},(-\Delta)^j y)-B(\tilde{v}_k,(-\Delta)^j y_k)\\
 =& B({v},(-\Delta)^j \tilde y)+(v(t_{\tilde{m}}),(-\Delta)^j y(t_{\tilde{m}}))-B({v}_k,(-\Delta)^j \tilde y_k)-(v_{k,\tilde{m}}^-,(-\Delta)^j y_{k,\tilde{m}}^+)\\ =&(v(t_{\tilde{m}}),(-\Delta)^j y(t_{\tilde{m}}))-(v_{k,\tilde{m}}^-,(-\Delta)^j y_{k,\tilde{m}}^+)\\
 =&(v(t_{\tilde{m}})-v_{k,\tilde{m}}^-,(-\Delta)^j y(t_{\tilde{m}}))+(v_{k,\tilde{m}}^-,(-\Delta)^j (y(t_{\tilde{m}})-y_{k,\tilde{m}}^+))\\
  =& I_1 +I_2.
 \end{aligned}
  $$
Note, that $B({v},(-\Delta)^j \tilde y)=0$ and $B({v}_k,(-\Delta)^j \tilde y_k)=0$ by construction.
  Using the Cauchy-Schwarz inequality, Lemma \ref{lemma: negative smooth superconvergence} with $l=r$ and $T=t_{\tilde{m}}$ and the smoothing estimate \eqref{eq: continuous smoothing} 
 $$
 \begin{aligned}
 I_1&\le \|(-\Delta)^{-2r-1} (v(t_{\tilde{m}})-v_{k,\tilde{m}}^-)\|_{L^2(\Om)}\|(-\Delta)^{2r+1+j}y(t_{\tilde{m}})\|_{L^2(\Om)}\\
 &\le C\frac{k^{2r+1}}{T^{2r+j+1}}\|v_0\|_{L^2(\Om)}\|(-\Delta)^j(v-v_{k})(T)\|_{L^2(\Om)}.
\end{aligned}
 $$  
Similarly,  using the Cauchy-Schwarz inequality, Lemma \ref{lemma: negative smooth superconvergence} with $l=r$ and $T=t_{\tilde{m}}$ and the semidiscrete smoothing estimate in Lemma~\ref{lemma: higher smoothing}
 $$
 \begin{aligned}
 I_2&\le \|(-\Delta)^{2r+j+1} v_{k,\tilde{m}}^-\|_{L^2(\Om)}\|(-\Delta)^{-2r-1}(y(t_{\tilde{m}})-y_{k,\tilde{m}}^+)\|_{L^2(\Om)}\\
 &\le C\frac{k^{2r+1}}{T^{2r+j+1}}\|v_0\|_{L^2(\Om)}\|(-\Delta)^j(v-v_{k})(T)\|_{L^2(\Om)}.
\end{aligned}
 $$
Combining the estimates for $I_1$ and $I_2$ and canceling $\|(-\Delta)^j(v-v_{k})(T)\|_{L^2(\Om)}$ on both sides,  we obtain the lemma.
\end{proof}
\subsection{Proof of Theorem \ref{thm:smoothing_est_time} }
We use the Gagliardo-Nirenberg inequality~\eqref{eq:Gagliardo-Nirenberg_Omega} and obtain
\[
\norm{v-v_k}_{L^\infty(\Omega)}\le C \norm{\Delta(v-v_k)}^{\frac{N}{4}}_{L^2(\Omega)} \norm{(v-v_k)}^{1-\frac{N}{4}}_{L^2(\Omega)}.
\]
Application of Lemma~\ref{lemma:higher_estimates_k} with $j=0$ and $j=1$ yields the result.

\subsection{Proof of Theorem \ref{thm:smoothing_est_time_grad}}
We  use  Gagliardo-Nirenberg inequality~\eqref{eq:general_Gagliardo-Nirenberg} with $B = \Omega_0$ as in the proof of Lemma \ref{lemma:cont_smoothing} and obtain
\[
\begin{aligned}
\linfnormn{\nabla (v-v_k)(T)} &\le  C \ltwonorm{\nabla \Delta (v-v_k)(T)}^{\frac{N}{4}} \ltwonorm{\nabla (v-v_k)(T)}^{1-\frac{N}{4}}\\
& \le C \ltwonorm{(v-v_k)(T)}^{\frac{1}{2}-\frac{N}{8}} \ltwonorm{\Delta (v-v_k)(T)}^{\frac{1}{2}} \ltwonorm{\Delta^2(v-v_k)(T)}^{\frac{N}{8}}. 
\end{aligned}
\]
Application of Lemma~\ref{lemma:higher_estimates_k} with $j=0$, $j=1$, and $j=2$ yields the result.

\section{Proof of smoothing error estimates in space}\label{sec:proof_thm_smoothing_space}
In this section we prove Theorem~\ref{thm:smoothing_est_space}. Before we provide the proof  we show the following results.
\begin{lmm}\label{prop: Phuk-ukh}
Let $v_k \in \Xk$ and $v_{kh}\in \Xkh$ be the semidiscrete and fully discrete solutions  of~\eqref{eq: semidiscrete heat} and~\eqref{eq:fully discrete heat}, respectively.
Then there exists a constant $C$ independent of $h$, $k$, and $T$ such that
$$
\|\Delta_h^{-1}(P_hv_{k}-v_{kh})(T)\|_{L^2(\Om)} \le Ch^2\lk\| v_0\|_{L^2(\Om)}.
$$
\end{lmm}
\begin{proof}
 Let $z_{kh} \in \Xkh$ be  the solution to a dual problem  with $z_{kh}(T)=\Delta_h^{-1}(P_h v_{k}-v_{kh})(T)$, i.e.
$$
B(\chi_{kh},z_{kh})=(\chi_{kh}(T),\Delta_h^{-1}(P_hv_{k}-v_{kh})(T)) \quad \text{for all } \; \chi_{kh} \in \Xkh.
$$
Then taking $\chi_{kh}=\Delta_h^{-1}(P_hv_{k}-v_{kh})$ by the Galerkin orthogonality, the stability of the $L^2$ projection, the standard elliptic error estimates, Lemma~\ref{lemma: stability dG_r fully discrete} and Corollary \ref{cor: maximal parabolic initial in L1}, we obtain
$$
\begin{aligned}
\|\Delta_h^{-1}(P_hv_{k}-v_{kh})(T)\|^2_{L^2(\Om)}&=B(\Delta_h^{-1}(P_hv_{k}-v_{kh}),z_{kh})\\
&=B(P_hv_{k}-v_{kh},\Delta_h^{-1}z_{kh})\\
&=B(P_hv_{k}-v_{k},\Delta_h^{-1}z_{kh})\\
&=(\na(P_hv_{k}-R_hv_{k}), \na\Delta_h^{-1}z_{kh})_{\IOm}\\
&=-(P_hv_{k}-R_hv_{k}, z_{kh})_{\IOm}\\
&\le \|P_hv_{k}-R_hv_{k}\|_{L^1(I;L^2(\Om))}\|z_{kh}\|_{L^\infty(I;L^2(\Om))}\\
&\le C\|v_{k}-R_hv_{k}\|_{L^1(I;L^2(\Om))}\|z_{kh}\|_{L^\infty(I;L^2(\Om))}\\
&\le Ch^2\|\Delta v_{k}\|_{L^1(I;L^2(\Om))}\|\Delta_h^{-1}(P_hv_{k}-v_{kh})(T)\|_{L^2(\Om)}\\
&\le Ch^2\lk\| v_0\|_{L^2(\Om)}\|\Delta_h^{-1}(P_hv_{k}-v_{kh})(T)\|_{L^2(\Om)}.
\end{aligned}
$$
Canceling, we obtain the result. 
\end{proof}

In order to establish optimal pointwise error estimates for $R_hv_k-v_{kh}$, we first show the corresponding estimate with respect to the $L^2(\Omega)$ norm and then for $\Delta_h(R_hv_k-v_{kh})$  in the $L^2(\Omega)$ norm likewise.
\begin{lmm}\label{lemma: u_k-u_kh}
Let $v_k \in \Xk$ and $v_{kh} \in \Xkh$ be the semidiscrete and fully discrete solutions  of~\eqref{eq: semidiscrete heat} and~\eqref{eq:fully discrete heat}, respectively.
There exists a constant $C$ independent of $k$, $h$, and $T$ such that
$$
\|(R_hv_k-v_{kh})(T)\|_{L^2(\Om)}\le \frac{Ch^2}{T}\lk\|v_0\|_{L^2(\Om)}.
$$
\end{lmm}
\begin{proof}
Let $y_{kh}\in \Xkh$ be the solution to a dual problem  with $y_{kh}(T)=(R_hv_k-v_{kh})(T)$, i.e. $y_{kh} \in \Xkh$ satisfies 
\begin{equation}\label{eq: dual with yT=R_hv_k-v_kh}
B(\varphi_{kh},y_{kh}) = (\varphi_{kh}(T),(R_hv_k-v_{kh})(T)), \quad \text{for all }\; \varphi_{kh}\in \Xkh.
\end{equation}
We abbriviate $\psi_{kh} = R_hv_k-v_{kh} \in \Xkh$ and set $\tilde \psi_{kh}$ to be zero on $I_1\cup \cdots \cup I_{\tilde{m}}$ for $\tilde{m}$ such that $\frac{T}{2}\in I_{\tilde{m}}$ and $\tilde \psi_{kh}= \psi_{kh}$ on the remaining time intervals. Similarly we define $\tilde y_{kh}$.  
Then by \eqref{eq: Bilinear tilde} and using the Galerkin orthogonality, we have 
$$
\begin{aligned}
\| (R_hv_k-v_{kh})(T)\|^2_{L^2(\Omega)}=B(\tilde \psi_{kh},y_{kh}) &=B(\psi_{kh},\tilde y_{kh})+(\psi_{kh,\tilde m}^-,y_{kh,\tilde{m}}^+)\\
&=B(R_hv_k-v_{kh},\tilde y_{kh})+((R_hv_{k,\tilde{m}}-v_{kh,\tilde{m}})^-,y_{kh,\tilde{m}}^+)\\
&=B(R_hv_k-v_{k},\tilde y_{kh})+((R_hv_{k,\tilde{m}}-v_{kh,\tilde{m}})^-,y_{kh,\tilde{m}}^+)\\
&=J_1+J_2.
\end{aligned}
$$
Using~\eqref{eq:B_dual} and the property of the Ritz projection, we have
$$
\begin{aligned}
J_1 &= - \sum_{m=\tilde{m}+1}^M ( R_hv_k-v_{k},\partial_t y_{kh})_{I_m \times \Omega}-\sum_{m=\tilde{m}+1}^{M} (R_hv_{k,m}^--v_{k,m}^-,[y_{kh}]_m)
-(R_hv_{k,\tilde{m}}^--v_{k,\tilde{m}}^-,y_{kh,\tilde{m}}^+) \\
&\le \|R_hv_k-v_{k}\|_{L^\infty((t_{\tilde{m}-1},T);L^2(\Om))}\left(\|\partial_t y_{kh}\|_{L^1(I;L^2(\Om))}+\sum_{m=1}^M\|[y_{kh}]_m\|_{L^2(\Om)}+\|y_{kh,\tilde{m}}^+\|_{L^2(\Om)}\right).
\end{aligned}
$$
By the Lemma \ref{lemma: homogeneous smoothing dG_r fully discrete} and Lemma~\ref{lemma: stability dG_r fully discrete}
we obtain
$$
\|\partial_t y_{kh}\|_{L^1(I;L^2(\Om))}+\sum_{m=1}^M\|[y_{kh}]_m\|_{L^2(\Om)}+\|y_{kh,\tilde{m}}^+\|_{L^2(\Om)}\le C\lk\| (R_hv_k-v_{kh})(T)\|_{L^2(\Omega)}.
$$
By the approximation properties of the Ritz projection, $H^2$ regularity, and using the fact that $\frac{T}{2}\in I_{\tilde{m}}$ we have
$$
\|R_hv_k-v_{k}\|_{L^\infty((t_{\tilde{m}-1},T);L^2(\Om))}\le Ch^2\|\Delta u_k\|_{L^\infty((t_{\tilde{m}-1},T);L^2(\Om))}\le C\frac{h^2}{T}\| v_0\|_{L^2(\Om)}.
$$
Canceling, be obtain the result for $J_1$. 

To estimate $J_2$ we add and subtract $v_{k,m}$. Thus we obtain
$$
\begin{aligned}
J_2 &= ((R_hv_{k,\tilde{m}}-v_{kh,\tilde{m}})^-,y_{kh,\tilde{m}}^+)\\
 &=((R_hv_{k,\tilde{m}}-v_{k,\tilde{m}})^-,y_{kh,\tilde{m}}^+)+((v_{k,\tilde{m}}-v_{kh,\tilde{m}})^-,y_{kh,\tilde{m}}^+) :=J_{21}+J_{22}.
\end{aligned}
$$
Similarly to the above, using Lemma \ref{lemma: homogeneous smoothing dG_r fully discrete} and Lemma~\ref{lemma: stability dG_r fully discrete} we obtain
$$
\begin{aligned}
J_{21}&\le \|R_hv_k-v_{k}\|_{L^\infty((t_{\tilde{m}-1},T);L^2(\Om))}\|y_{kh,\tilde{m}}^+\|_{L^2(\Om)}\\
&\le Ch^2\|\Delta v_{k}\|_{L^\infty((t_{\tilde{m}-1},T);L^2(\Om))}\| (R_hv_k-v_{kh})(T)\|_{L^2(\Omega)}\\
&\le  \frac{Ch^2}{T}\lk\| v_0\|_{L^2(\Om)}\| (R_hv_k-v_{kh})(T)\|_{L^2(\Omega)}.
\end{aligned}
$$
To estimate $J_{22}$ we use Lemma~\ref{prop: Phuk-ukh} with $T = t_{\tilde m}$ and the fact the constant there does not depend on $T$ together with Lemma~\ref{lemma: stability dG_r fully discrete}. Hence, 
$$
\begin{aligned}
J_{22} &= ((P_h v_{k,\tilde{m}}-v_{kh,\tilde{m}})^-,y_{kh,\tilde{m}}^+)\\
&\le \|\Delta_h^{-1}(P_hv_{k,\tilde{m}}-v_{kh,\tilde{m}})^-\|_{L^2(\Om)}\|\Delta_h y_{kh,\tilde{m}}^+\|_{L^2(\Om)}\\
& \le  \frac{Ch^2}{T}\lk\| v_0\|_{L^2(\Om)}\| (R_hv_k-v_{kh})(T)\|_{L^2(\Omega)}.
\end{aligned}
$$
Canceling, we obtain the lemma.
\end{proof}

Next we establish the following smoothing result in $L^2$ norm with discrete Laplacian.
\begin{lmm}\label{lemma: Delta_h u_k-u_kh}
Let $v_k \in \Xk$ and $v_{kh} \in \Xkh$ be the semidiscrete and fully discrete solutions  of~\eqref{eq: semidiscrete heat} and~\eqref{eq:fully discrete heat}, respectively.
There exists a constant $C$ independent of $k$, $h$, and $T$ such that
$$
\|\Delta_h(R_hv_k-v_{kh})(T)\|_{L^2(\Om)}\le \frac{Ch^2}{T^{2}}\lk\|v_0\|_{L^2(\Om)}.
$$
\end{lmm}
\begin{proof}
Let $y_{kh} \in \Xkh$ be the solution to a dual problem  with $y_{kh}(T)=\Delta_h(R_hv_k-v_{kh})(T)$, i.e. $y_{kh}$ satisfies 
\begin{equation}\label{eq: dual with yT=Deltah_R_hu_k-u_kh}
B(\varphi_{kh},y_{kh}) = (\varphi_{kh}(T),\Delta_h(R_hv_k-v_{kh})(T)), \quad \text{for all } \; \varphi_{kh}\in \Xkh.
\end{equation}
As in the proof of the previous lemma we abbriviate $\psi_{kh} = R_hv_k-v_{kh}$ and set $\tilde \psi_{kh}$ to be zero on $I_1\cup \cdots \cup I_{\tilde{m}}$ for some $\tilde{m}$ to be specified later and $\tilde \psi_{kh} = \psi_{kh}$ on the remaining time intervals.  Similarly we define $\tilde y_{kh}$.
Then setting $\varphi_{kh} = \Delta_h \tilde \psi_{kh}$ and using the Galerkin orthogonality and~\eqref{eq: Bilinear tilde}, we have 
$$
\begin{aligned}
\| \Delta_h(R_hv_k-v_{kh})(T)\|^2_{L^2(\Omega)}&=B(\Delta_h \tilde \psi_{kh},y_{kh})\\
&=B(\tilde \psi_{kh},\Delta_h y_{kh})\\
&=B(\psi_{kh},\Delta_h \tilde y_{kh}) + (\psi_{kh,\tilde m}^-,y_{kh,\tilde m}^+)\\
&=B(R_hv_k-v_{kh},\Delta_h\tilde y_{kh})+((R_hv_{k,\tilde{m}}-v_{kh,\tilde{m}})^-,\Delta_h y_{kh,\tilde{m}}^+)\\
&=B(R_hv_k-v_{k},\Delta_h\tilde y_{kh})+((R_hv_{k,\tilde{m}}-v_{kh,\tilde{m}})^-,\Delta_h y_{kh,\tilde{m}}^+)\\
&=J_1+J_2.
\end{aligned}
$$
Choosing $\tilde{m}$ such that $\frac{T}{4}\in I_{\tilde{m}}$, using the definition of the bilinear form  $B(\cdot,\cdot)$, the discrete maximal parabolic regularity from Corollary~\ref{cor: maximal parabolic initial in L1} and  Lemma~\ref{lemma: higher smoothing} we obtain
$$
\begin{aligned}
J_1&=  \sum_{m=\tilde{m}}^M (\partial_t (R_hv_k-v_{k}), \Delta_hy_{kh})_{I_m \times \Omega}+\sum_{m=\tilde{m}}^{M} ([R_hv_{k}-v_{k}]_m,\Delta_h y_{kh,m}^+) \\
&\le \|\pa_t(R_hv_k-v_{k})\|_{L^\infty((t_{\tilde{m}},T);L^2(\Om))}\|\Delta_h y_{kh}\|_{L^1((t_{\tilde{m}},T);L^2(\Om))}\\
&+\max_{\tilde{m}\le m\le M}\left\{k_m^{-1}\|[R_hv_{k}-v_{k}]_m\|_{L^2(\Om)}\right\}\sum_{m=\tilde{m}}^Mk_m\|\Delta_h y_{kh,m}^+\|_{L^2(\Om)}\\
&\le C \lk h^2\left(\|\pa_t\Delta v_{k}\|_{L^\infty((t_{\tilde{m}},T);L^2(\Om))}
+\max_{\tilde{m}\le m\le M}\left\{k_m^{-1}\|[\Delta v_{k}]_m\|_{L^2(\Om)}\right\}\right) \|\Delta_h(R_hv_k-v_{kh})(T)\|_{L^2(\Om)}\\
&\le \frac{Ch^2}{T^2}\lk\|v_0\|_{L^2(\Omega)}\|\Delta_h(R_hu_k-u_{kh})(T)\|_{L^2(\Om)}.
\end{aligned}
$$
Canceling we obtain the desired estimate for $J_1$.

To estimate $J_2$ we proceed as in the proof of the previous lemma,
\[
J_2=((R_hv_{k,\tilde{m}}-v_{k,\tilde{m}})^-,\Delta_h y_{kh,\tilde{m}}^+)+((v_{k,\tilde{m}}-v_{kh,\tilde{m}})^-,\Delta_h y_{kh,\tilde{m}}^+) :=J_{21}+J_{22}.
\]
Similarly to the above,
$$
\begin{aligned}
J_{21}&\le \|R_hv_k-v_{k}\|_{L^\infty((t_{\tilde{m}-1},T);L^2(\Om))}\|\Delta_h y_{kh,\tilde{m}}^+\|_{L^2(\Om)}\\
&\le \frac{Ch^2}{T}\|\Delta v_{k}\|_{L^\infty((t_{\tilde{m}-1},T);L^2(\Om))}\| \Delta_h(R_hv_k-v_{kh})(T)\|_{L^2(\Omega)}\\
&\le  \frac{Ch^2}{T^2} \| v_0\|_{L^2(\Om)}\| \Delta_h(R_hv_k-v_{kh})(T)\|_{L^2(\Omega)}.
\end{aligned}
$$
To estimate $J_{22}$, we proceed as in the proof of the  previous lemma, and using Lemma~\ref{prop: Phuk-ukh}, we have
$$
\begin{aligned}
J_{22} &= ((P_h v_{k,\tilde{m}}-v_{kh,\tilde{m}})^-,\Delta_h y_{kh,\tilde{m}}^+)\\
&\le \|\Delta_h^{-1}(P_hv_{k,\tilde{m}}-v_{kh,\tilde{m}})^-\|_{L^2(\Om)}\|\Delta^2_h y_{kh,\tilde{m}}^+\|_{L^2(\Om)}\\
& \le \frac{Ch^2}{T^2}\lk\| v_0\|_{L^2(\Om)}\|\Delta_h(R_hv_k-v_{kh})(T)\|_{L^2(\Omega)}.
\end{aligned}
$$
Canceling we obtain the lemma.
\end{proof}

As a consequence of the two lemmas above and the discrete Gagliardo-Nirenberg  inequality \eqref{eq: discrete Gagliardo-Nirenber}
we immediately obtain the following result.
\begin{lmm}\label{lemma: R_hu_k-u_kh in Linfty}
Let $v_k$ and $v_{kh}$ be the semidiscrete and fully discrete solutions  of   \eqref{eq: semidiscrete heat} and \eqref{eq:fully discrete heat}, respectively. Then, 
there exists a constant $C$ independent of $k$ and $h$ such that
$$
\|(R_hv_k-v_{kh})(T)\|_{L^\infty(\Om)}\le \frac{Ch^2}{T^{1+\frac{N}{4}}}\lk\|v_0\|_{L^2(\Om)}.
$$
\end{lmm}
\begin{rmrk}
The result in Lemma \ref{lemma: R_hu_k-u_kh in Linfty} is rather interesting and of independent interest. It shows that the $L^\infty(\Omega)$ error between the semidiscrete solution and its Ritz projection for piecewise linear elements is of optimal second order even if the exact solution at a final time $T$ is not in $W^{2,\infty}(\Om)$  as for example in the case if the domain $\Om$ has strong corner singularities. This in particular shows a well-known fact that the presence of corner singularities is essentially an elliptic problem.   

Similar to the results in \cite{ThomeeV_XuJ_NaiY_1989}, it also can be used to  show a superconvergent result for the gradient.  Thus for $N=2$ using the discrete  Sobolev inequality 
\begin{equation}\label{eq: grad discrete Sobolev}
    \|\nabla \chi\|_{L^\infty(\Om)}\le C\lh^{1/2}\|\Delta_h \chi\|_{L^2(\Om)},\quad \text{for all }\; \chi\in V_h,
\end{equation}
we also have the following superconvergent estimate
$$
\|\na(R_hv_k-v_{kh})(T)\|_{L^\infty(\Om)}\le \frac{Ch^2}{T^2}\lh^{1/2}\lk\|v_0\|_{L^2(\Om)}.
$$
\end{rmrk}

\subsection{Proof of Theorem~\ref{thm:smoothing_est_space} }
Adding and subtracting $R_hv_k$, we have
\[
\abs{(v_k-v_{kh})(T,x_0)} \le \abs{(R_hv_k-v_{kh})(T,x_0)}+\abs{(v_k-R_hv_k)(T,x_0)}.
\]
From Lemma~\ref{lemma: R_hu_k-u_kh in Linfty} we have
$$
\abs{(R_hv_k-v_{kh})(T,x_0)}\le C(T)h^2\lk\|v_0\|_{L^2(\Om)}.
$$
Using the pointwise interior elliptic results from \cite{AHSchatz_LBWahlbin_1977a}
we have 
$$
\begin{aligned}
\abs{(v_k-R_hv_k)(T,x_0)}&\le C\lh\min_{\chi\in V_h}\|v_k(T)-\chi\|_{L^\infty(B_d(x_0))}+Cd^{-\frac{N}{2}}\|(v_k-R_hv_k)(T)\|_{L^2(\Omega)}\\
&\le Ch^{2}\lh\|v_k(T)\|_{W^{2,\infty}(B_d(x_0))}+Ch^2d^{-\frac{N}{2}}\|v_k(T)\|_{H^2(\Omega)}\\
&\le Ch^{2}\lh\| v_k(T)\|_{H^{4}(B_d(x_0))}+Ch^2d^{-\frac{N}{2}}\|\Delta v_k(T)\|_{L^2(\Omega)}\\
&\le Ch^{2}\lh\|\Delta v_k(T)\|_{H^{2}(\Omega)}+Ch^2d^{-\frac{N}{2}}\|\Delta v_k(T)\|_{L^2(\Omega)}\\
&\le Ch^{2}\lh\|\Delta^2 v_k(T)\|_{L^2(\Omega)}+Ch^2d^{-\frac{N}{2}}\|\Delta v_k(T)\|_{L^2(\Omega)}\\
&\le Ch^{2}\left(T^{-2}\lh+d^{-\frac{N}{2}}T^{-1}\right)\|v_0\|_{L^2(\Om)}
,
\end{aligned}
$$
where we used  the embedding $H^4(B_d(x_0)) \hookrightarrow C^2(B_d(x_0))$, the interior regularity result~\cite[Chapter 6.3,Theorem 2]{LCEvans_2010} and convexity of $\Omega$.

\section{Algorithmic treatment}\label{sec:PDAP}
This section is devoted to the algorithmic solution of the sparse initial data identification problem under consideration.  Let us first note that by Theorem~\ref{theorem:existence_kh} we can look for a minimizer $\oq_{kh}$ in the space $\M_h$ consisting of linear combinations of Diracs concentrated in the interior nodes $\mathcal{N}_h$ of the mesh, i.e.
\[
\oq_{kh} = \sum_{x_i \in \mathcal{N}_h}\bar{\gamma}_{kh,i} \delta_{x_i}, 
\]
where $\bar \gamma_{kh} \in \R^{\# \mathcal{N}_h}$  is a vector of optimal coefficients. Thus, the fully discrete problem~\eqref{eq:red_problem_kh} can be equvalently reformulated as a finite dimensional problem (of dimension $\# \mathcal{N}_h$) in the coefficients $\gamma_{kh}$ with an $l_1$ regularization term leading to
\begin{align*}
\operatorname{Minimize }\; \frac{1}{2} \left\|\sum_{x_i \in \mathcal{N}_h}\gamma_{kh,i} S_{kh}(\delta_{x_i})-u_d\right\|^2_{L^2(\Omega)}+ \alpha \sum^{\# \mathcal{N}_h}_{i=1} |\gamma_{kh,i}| ,\quad \gamma_{kh} \in \R^{\#\mathcal{N}_h}.
\end{align*}
This problem can be solved  by a variety of efficient solution algorithms, e.g., semi-smooth Newton methods,~\cite{milzarekfilter}, or FISTA,~\cite{fista}. 
However, a direct application of finite dimensional optimization algorithms to this problem may lead to mesh-dependent methods, whose convergence behavior critically depends on the fineness of the discretization. In contrast, we employ an optimization algorithm, which can be described on the continuous level, as a solution algorithm for the problem~\eqref{eq:red_problem}. This algorithm can be then directly adapted to the discretized problem~\eqref{eq:red_problem_kh}. Since the convergence properties of the presented algorithm can be analyzed on the continuous level, see~\cite{Walter:2019}, one expects mesh independent behavior for the discretized problem, which is also confirmed by our numerical results.
 
We propose a version of the Primal-Dual-Active-Point (PDAP) method from~\cite{Walter:2019} which iteratively generates a sequence of finite linear combinations of Dirac delta functions. The algorithm on the continuous level is briefly described and its convergence properties are summarized below.
Given an ordered set of finitely many points~$\mathcal{A}=\{x_i\}^K_{i=1}$ define the parametrization
\begin{align*}
Q_{\mathcal{A}} \colon \R^{\# \mathcal{A}} \to \M(\Omega), \quad \beta \mapsto \sum_{x_i \in \mathcal{A}} \beta_i \delta_{x_i} 
\end{align*}
as well as the finite dimensional subproblem
\begin{align} \label{def:finitesubprob}
\text{Minimize}~j \left( Q_{\mathcal{A}}(\beta) \right),~\beta \in \mathbb{R}^{\#\mathcal{A}}.
\end{align}

We initialize the proposed algorithm with a sparse initial iterate~$q_0 \in \M(\Omega)$,~$\# \supp q_0 < \infty$. In the~$n$-th iteration, a new support point~$\hat{x}^n \in \Omega$ is determined based on the violation of the condition (a) in Corollary~\ref{cor:structure_opt_control} by the current adjoint state~$z^n(0)= S^*(Sq_n-u_d)$. Subsequently, the new iterate is found as~$q_{n+1}=Q_{\mathcal{A}_n}(\bar{\beta}^{n+1})$ where~$\bar{\beta}^{n+1} \in \R^{\# \mathcal{A}_n}$ is a solution to~\eqref{def:finitesubprob} for~$\mathcal{A}=\mathcal{A}_n$. Thus, the method alternates between updating the active set~$\mathcal{A}_n$ by adding~$\hat{x}^n$ to the support of the current iterate~$q_n$ and computing a minimizer of~$j$ over~$\M(\mathcal{A}_n)$. 
The procedure is summarized in Algorithm~\ref{alg:PDAP}.
\begin{algorithm}
\begin{algorithmic}
\STATE 1. Choose $q_0 \in \M(\Omega)$, $\# \supp q_0 < \infty$. Set $M_0= j(q_0)/ \beta$.
 \WHILE {$\Phi(q_n)\geq \mathrm{TOL}$}
 \STATE 2.
 Compute $z_n(0)= S^*(Sq_n-u_d) \in C_0 (\Omega)$. Determine $\hat{x}^n \in \argmax_{x \in \Omega} |z_n(0,x)|$.
 \STATE 3.
 Set $\mathcal{A}_n= \supp q_n \cup \{\hat{x}^n\}$ and compute a solution~$\bar{\beta}^{n+1}$ to~\eqref{def:finitesubprob} with $\mathcal{A}=\mathcal{A}_n$.
 \STATE 4.
 Set~$q_{n+1}= Q_{\mathcal{A}_n}(\bar{\beta}^{n+1})$.
\ENDWHILE
\end{algorithmic}
\caption{Primal-Dual-Active-Point method}\label{alg:PDAP}
\end{algorithm}
Note that the support of~$q_n$ is pruned after each iteration i.e. Dirac delta functions with zero coefficients are removed from the iterate. Additionally, we observe that Algorithm~\ref{alg:PDAP} is monotonous, i.e.~$j(q_{n+1})\leq j(q_n)$, and thus also~$\mnorm{q_n}\leq M_0$ for all~$n \in \N$. To monitor the convergence of the algorithm we consider the primal-dual-gap functional $\Phi \colon \M(\Omega) \to \R_+$ which is defined as
\begin{align*}
\Phi(q)= \max_{ \mnorm{\delta q} \leq M_0} \lbrack \langle z[q](0), q-\delta q \rangle +\beta \mnorm{q}-\beta \mnorm{\delta q} \rbrack    \quad \text{where} \quad  z[q](0)= S^*(Sq-u_d)
\end{align*}
and~$M_0= j(q_0)/\alpha$.
This is justified by the following lemma, see~\cite[Lemma 6.12 and Lemma 6.41]{Walter:2019}.
\begin{lmm} \label{prop:primaldualgap}
There holds~$\Phi(q)\geq 0$ for all~$q\in \M(\Omega)$ with equality if and only if~$q=\bar{q}$ is the optimal solution of~\eqref{eq:red_problem}. Furthermore we have
\begin{align*}
j(q)-j(\bar{q}) \leq \Phi(q) \quad \forall q \in \M(\Omega).  
\end{align*}
Let~$\{q_n\}$ denote the sequence generated by Algorithm~\ref{alg:PDAP}. Then we obtain
\begin{align} \label{eq:primaldualforit}
\Phi(q_n)=M_0(\norm{z_n(0)}_{C(\Omega)}-\beta) \quad \text{where} \quad z_n(0)=S^*(Sq_n-u_d),
\end{align}
for~$n \geq 1$.
\end{lmm}

We point out that, due to~\eqref{eq:primaldualforit}, the primal-dual-gap~$\Phi(q_n)$ can be cheaply computed as a byproduct of step $2$. 

The following theorem, see~\cite[Theorem 6.43]{Walter:2019}, provides two convergence results. For the general case we obtain sub-linear convergence of the the cost functional. Under Assumption~\ref{ass:finite}, we obtain linear convergence for the functional, positions of the Diracs and for the corresponding coefficients.
\begin{thrm} \label{thm:convofpdap}
Let the sequence~$\{q_n\} \subset \M(\Omega)$ be generated by Algorithm~\ref{alg:PDAP} starting from~$q_0 \in \M(\Omega)$. Then we have
\begin{align*}
q_n \overset{\ast}{\rightharpoonup} \bar{q}, \quad \Phi(q_n) \rightarrow 0, \quad j(q_n)-j(\bar{q})\leq \frac{j(q_0)-j(\bar{q})}{1+c_1 n}
\end{align*}
for all~$n \in \mathbb{N}$ and a constant~$c_1>0$. If Assumption~\ref{ass:finite} holds, then~$\bar{q}= \sum^K_{i=1} \bar{\beta}_i \delta_{\bar{x}_i}$ and there exist $R,~c_2 >0$, $\zeta \in (0,1)$ with
\begin{align*}
\supp q_n \subset \bigcup^K_{i=1} \bar{B}_{R}(\bar{x}_i), \quad \bar{B}_{R}(\bar{x}_i) \cap \bar{B}_{R}(\bar{x}_j)=\emptyset, \quad \supp q_n\cap \bar{B}_R(\bar{x}_i) \neq \emptyset, \quad i,j= 1,2, \dots,K
\end{align*}
as well as
\begin{align*}
j(q_n)-j(\bar{q})+ \max_{i=1, \dots, K} \left \{ |q_n\left(\bar{B}_{R}(\bar{x}_i)\right)-\bar{\beta}_i|  +  \max_{x \in \supp q_n \cap \bar{B}_{R}(\bar{x}_i) } |x- \bar{x}_i|_{\R^N}\right\}
 \leq c_2 \zeta^n 
\end{align*}
for all~$n \in \mathbb{N}$ large enough.
\end{thrm}

We emphasize that the adaption of Algorithm~\ref{alg:PDAP} to the discrete problem~\eqref{eq:red_problem_kh} is straightforward. In detail, we replace the control-to-state operator~$S$ by its fully discrete counterpart~$S_{kh}$ and compute $z^{n,+}_{kh,0}=S_{kh}^*(S_{kh}q_n-u_d)$.
Moreover, in view of Theorem~\ref{theorem:existence_kh}, the search for the maximizer~$\hat{x}^n$ in step~$2$ can be restricted to the set of interior nodes~$\mathcal{N}_h$. The new coefficient vector~$\bar{\beta}^{n+1}$ is then found as solution to the finite-dimensional subproblem
\begin{align} \label{def:fulldiscretesub}
\text{Minimize}~j_{kh} \left( Q_{\mathcal{A}_n}(\beta) \right),~\beta \in \mathbb{R}^{\#\mathcal{A}_n}.
\end{align}
Note, that the support of $q_n$ usually consists of only few points, i.e. the dimension of the subproblem~\eqref{def:fulldiscretesub} is small the this subproblem can be solved efficiently by existing finite dimensional algorithms. In our numerical realization we use the semi-smooth Newton method for the solution of~\eqref{def:fulldiscretesub}.
\section{Numerical examples}
The final section is devoted to the presentation of numerical experiments which serve to underline the practical applicability of the proposed sparse control approach as well as to verify the derived theoretical results.
Throughout the section, the spatial domain is given by the unit square~$\Omega=(0,1)\times (0,1)$ and the final time is set to~$T=0.1$. All arising discrete optimal control problems are solved by an adaptation of the PDAP method, Algorithm~\ref{alg:PDAP}, as described at the end of the previous section.
\subsection{Identification of point sources} \label{subsec:inverseprob}
\begin{figure}[htb] 
\centering
\begin{subfigure}[t]{.45\linewidth}
\centering
\scalebox{.538}
{
%
%
\definecolor{mycolor1}{rgb}{0.00000,0.44700,0.74100}%
\begin{tikzpicture}

\begin{axis}[%
width=4.521in,
height=3.527in,
at={(0.758in,0.519in)},
scale only axis,
xmin=0,
xmax=1,
xlabel style={font=\color{white!15!black}},
xlabel={$\text{x}_\text{1}$},
ymin=0,
ymax=1,
ylabel style={font=\color{white!15!black}},
ylabel={$\text{x}_\text{2}$},
axis background/.style={fill=white}
]
\addplot [color=mycolor1, forget plot]
  table[row sep=crcr]{%
0	0\\
};
\addplot [color=blue, draw=none, mark=x, mark options={solid, blue}, forget plot]
  table[row sep=crcr]{%
0.263091083266217	0.258378565204941\\
};
\node[above, align=center]
at (axis cs:0.263,0.258) {-10.00};
\addplot [color=blue, draw=none, mark=x, mark options={solid, blue}, forget plot]
  table[row sep=crcr]{%
0.76061544960808	0.734190309666141\\
};
\node[above, align=center]
at (axis cs:0.761,0.734) {25.00};
\end{axis}
\end{tikzpicture}%
}\\
\caption{Reference source~$q^\dagger$}
\end{subfigure}
\begin{subfigure}[t]{.45\linewidth}
\centering
\scalebox{.47}
{\includegraphics{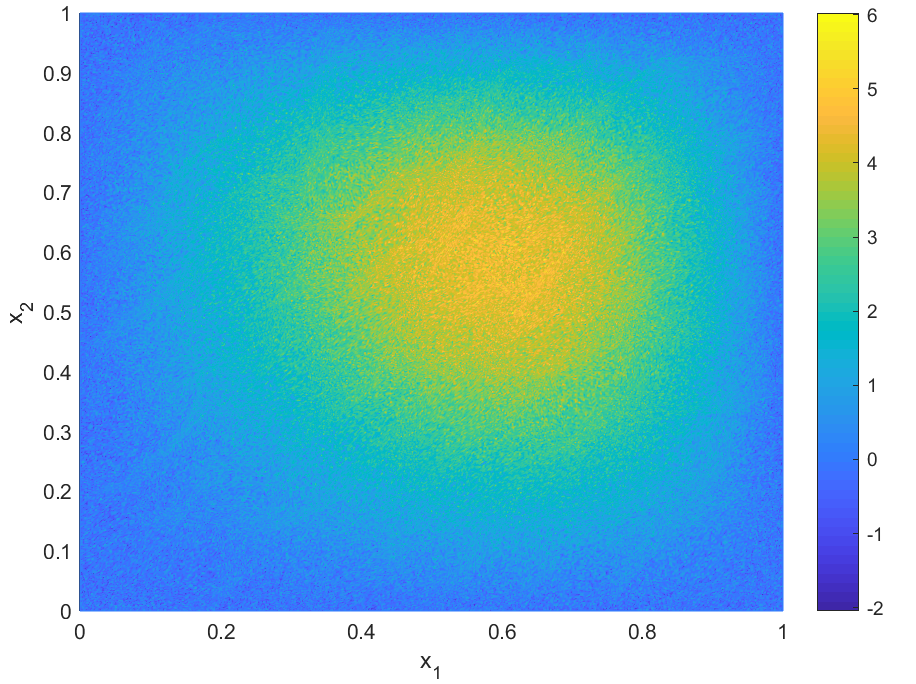}}\\
\caption{Noisy observation~$u_{\text{obs}}$}
\end{subfigure}
\caption{Inverse problem setup}
\label{fig:setup}
\end{figure}
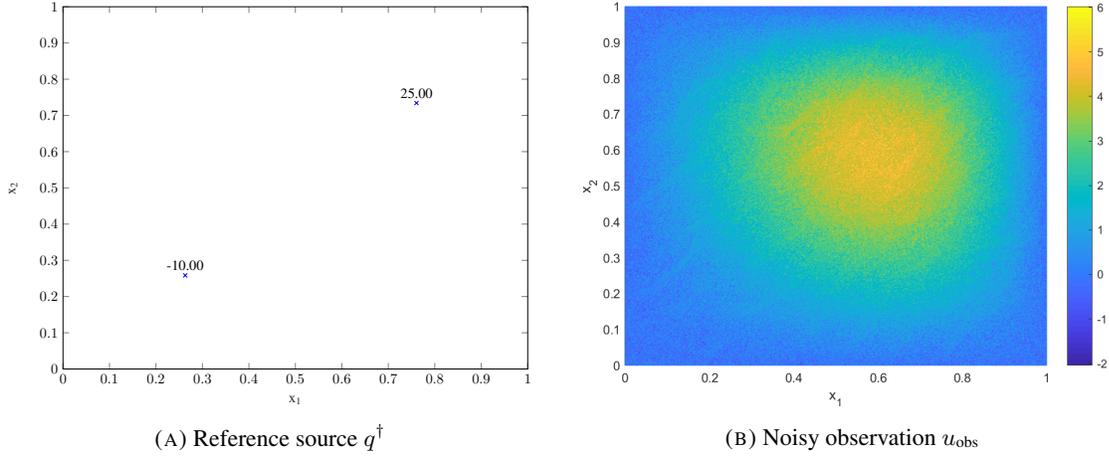

\begin{figure}[htb] 
\centering
\begin{subfigure}[t]{.45\linewidth}
\centering
\scalebox{.5}
{
%
%
\definecolor{mycolor1}{rgb}{0.00000,0.44700,0.74100}%
\begin{tikzpicture}

\begin{axis}[%
width=4.521in,
height=3.527in,
at={(0.758in,0.519in)},
scale only axis,
xmin=0,
xmax=1,
xlabel style={font=\color{white!15!black}},
xlabel={$\text{x}_\text{1}$},
ymin=0,
ymax=1,
ylabel style={font=\color{white!15!black}},
ylabel={$\text{x}_\text{2}$},
axis background/.style={fill=white}
]
\addplot [color=mycolor1, forget plot]
  table[row sep=crcr]{%
0	0\\
};
\addplot [color=blue, draw=none, mark=x, mark options={solid, blue}, forget plot]
  table[row sep=crcr]{%
0.699360444627145	0.689105463441521\\
};
\node[above, align=center]
at (axis cs:0.699,0.689) {19.31};
\addplot [color=blue, draw=none, mark=x, mark options={solid, blue}, forget plot]
  table[row sep=crcr]{%
0.295575037477538	0.308762864982223\\
};
\node[above, align=center]
at (axis cs:0.296,0.309) {-9.54};
\end{axis}
\end{tikzpicture}%
}\\
\caption{Recovered measure~$\oqkh$}
\end{subfigure}
\begin{subfigure}[t]{.45\linewidth}
\centering
\scalebox{.47}
{\includegraphics{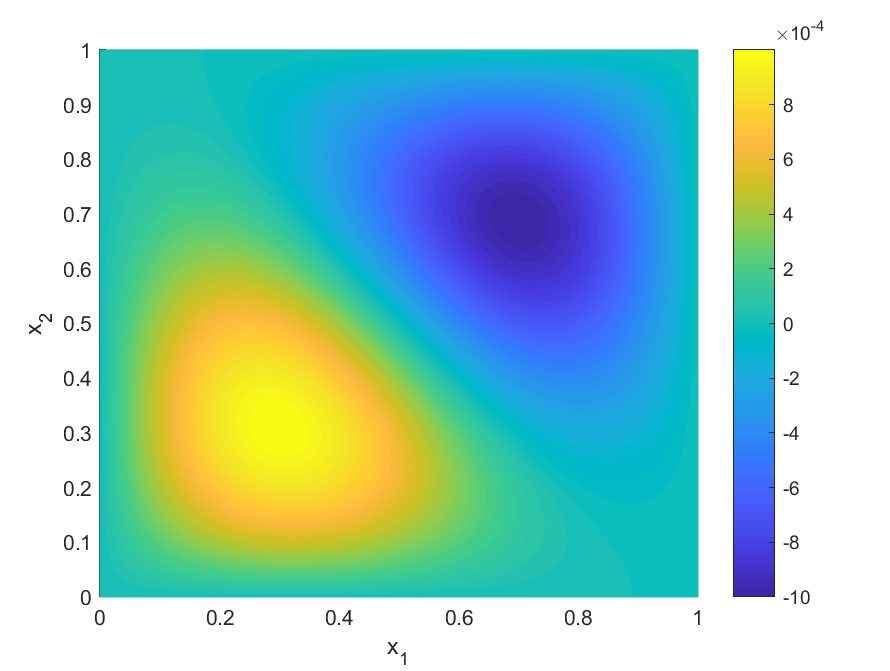}}\\
\caption{Initial adjoint state~$\oz_{kh,0}^+$}
\end{subfigure}
\begin{subfigure}[t]{.45\linewidth}
\centering
\scalebox{.5}
{\includegraphics{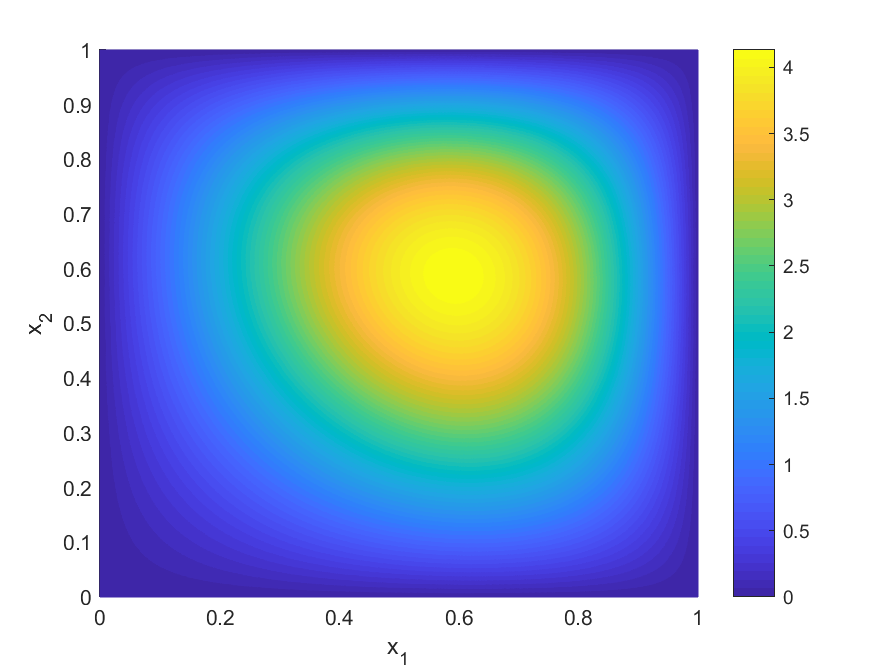}}\\
\caption{Recovered state~$S_{kh}(\oukh)$}
\end{subfigure}
\begin{subfigure}[t]{.45\linewidth}
\centering
\scalebox{.5}
{\includegraphics{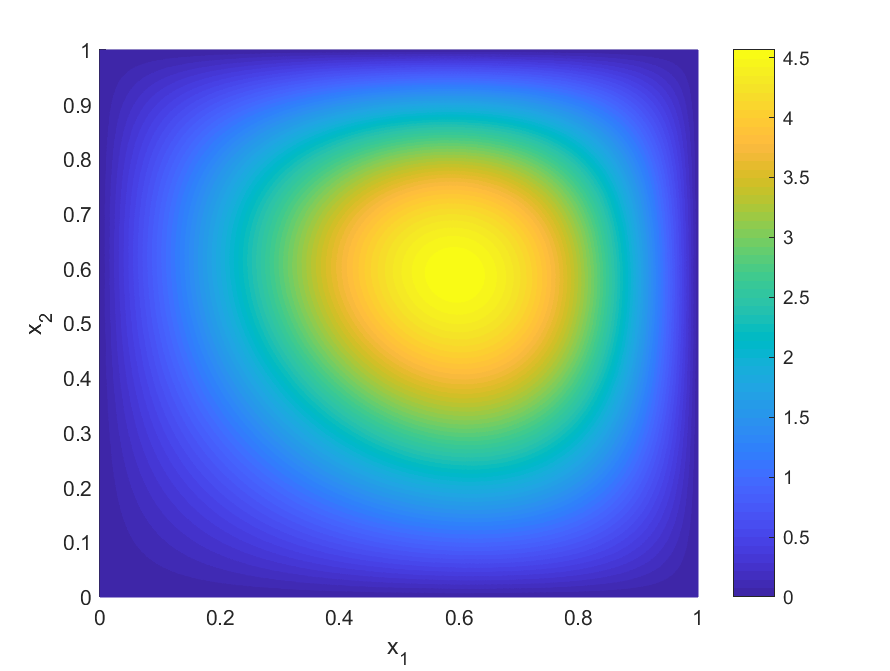}}\\
\caption{Reference state~$S_{kh}(q^\dagger)$}
\end{subfigure}
\caption{Reconstruction results}
\label{fig:reconstructionresults}
\end{figure}
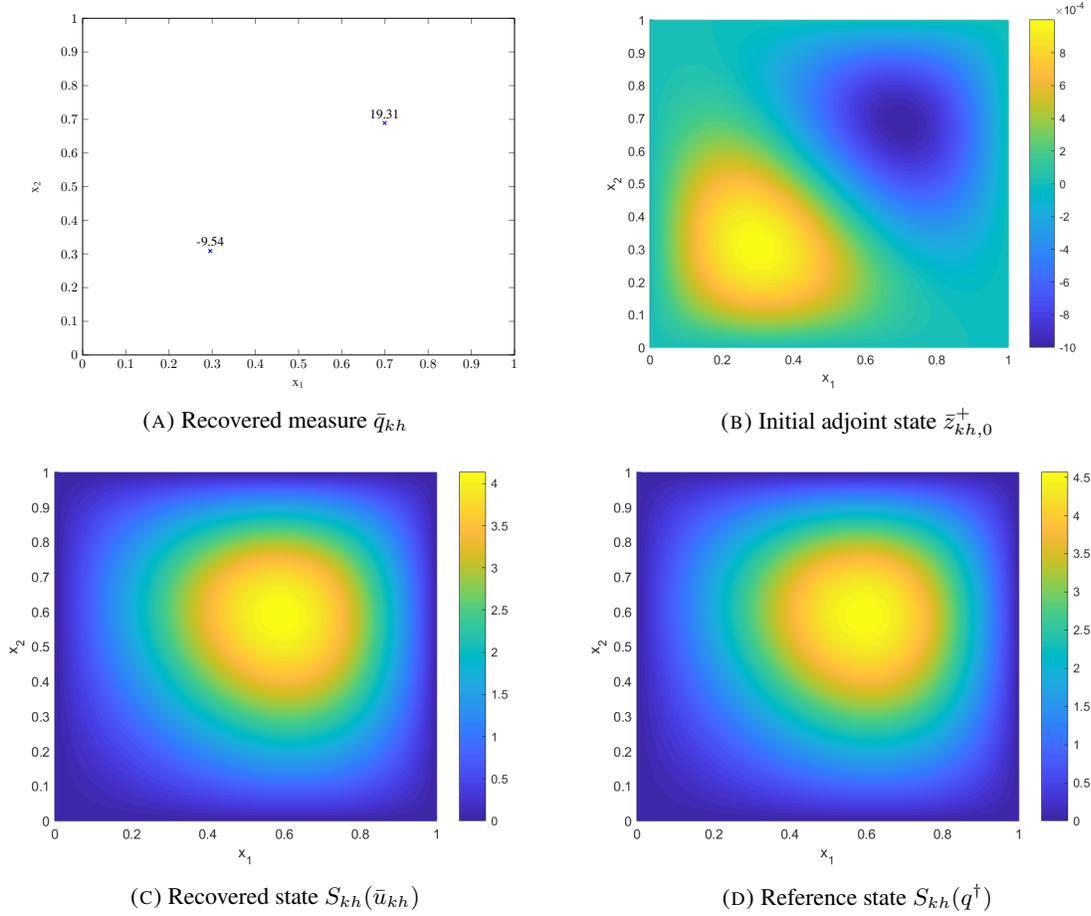

First, we aim to identify a sparse source term~$q^\dagger=-10 \delta_{x_1}+25 \delta_{x_2}$ from noisy observations of~$u(T)=S(q^\dagger)$. The time interval~$(0,T]$ is uniformly partitioned into~$M=256$ subintervals, the spatial domain~$\Omega$ is divided into triangles, see the description in Section~\ref{sec:discretization}. We emphasize that the support points~$x_1$ and~$x_2$, respectively, correspond to nodes of the triangulation. For the discretization of the state equation a~cG($1$)dG($0$) (i.e $r=0$) approximation is considered. The observations are given by~$u_{obs}= S_{kh}(q^\dagger)+\delta$ where~$\delta \in L^2(\Omega)$ is a given noise term. We plot~$u_{\text{obs}}$ alongside~$q^\dagger$ in Figure~\ref{fig:setup}.
To reconstruct~$q^\dagger$ from the given final time observation we propose to solve~$\eqref{eq:red_problem_kh}$ with~$u_d=u_{\text{obs}}$. For the described example we empirically determine~$\alpha=0.001$ as a suitable regularization parameter. Applying the Primal-Dual-Active-Point method to~$\eqref{eq:red_problem_kh}$ yields a reconstruction~$\oqkh \in ~\mathcal{M}_h $ with~$\#\supp \oqkh=3$. By a closer inspection, two of its support points are located in adjacent nodes of the triangulation. A possible explanation for this clustering of support points is provided by Theorem~\ref{theorem:existence_kh}. More in detail, a spike appearing in a discrete optimal solution~$\bar{q}\not \in \mathcal{M}_h$ at an off-grid location will appear as several nodal Dirac delta functions in the projected solution~$\Lambda_h \bar{q}$. For a better visualization of the results we replace the Dirac delta functions associated to the clustering support points by a single one of the same combined mass located at the center of gravity of the original positions. The post-processed measure is depicted in Figure~\ref{fig:reconstructionresults} together with~$\oz_{kh,0}^+=S^*_{kh}(S_{kh}\oqkh-u_d)$.

As predicted by Corollary~\ref{cor:opt_cond_kh} we have~$|\oz_{kh,0}^+(x)|\leq \alpha$ for all~$x\in \Omega$ and equality holds at the support points of~$\oqkh$. Moreover, we also plot the final state~$S_{kh}(q^\dagger)$ corresponding to the initial source~$q^\dagger$ as well as the reconstructed final state~$S_{kh}(\oqkh)$. We see that the proposed sparse control approach together with the lumping of clustering support points recovers the main structural features of the source~$q^\dagger$. In particular, we point out to the correct number of points sources as well as quantitatively reasonable estimates of their locations and coefficients. Note that we cannot expect the exact recovery of~$q^\dagger$ due to the appearance of the noise term~$\delta$ as well as the nonzero regularization parameter~$\alpha$. We specifically stress that~$\supp \oqkh \cap \supp q^\dagger =\emptyset $.

\subsection{Space refinement}
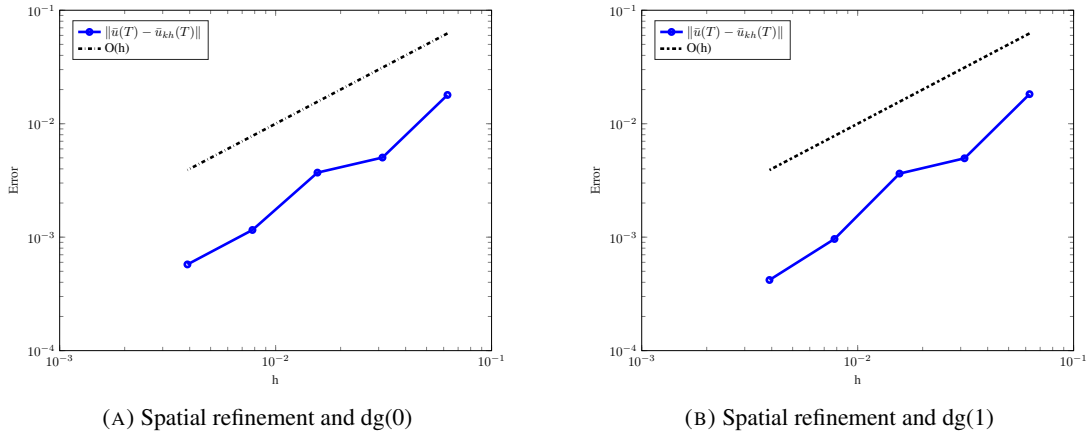
\begin{figure}[htb] 
\centering
\begin{subfigure}[t]{.45\linewidth}
\centering
\scalebox{.5}
{
%
%
\begin{tikzpicture}

\begin{axis}[%
width=4.521in,
height=3.566in,
at={(0.758in,0.481in)},
scale only axis,
xmode=log,
xmin=0.001,
xmax=0.1,
xminorticks=true,
xlabel style={font=\color{white!15!black}},
xlabel={h},
ymode=log,
ymin=0.0001,
ymax=0.1,
yminorticks=true,
ylabel style={font=\color{white!15!black}},
ylabel={Error},
axis background/.style={fill=white},
legend style={at={(0.03,0.97)}, anchor=north west, legend cell align=left, align=left, draw=white!15!black}
]
\addplot [color=blue, mark=o, mark options={solid, blue},line width=2pt]
  table[row sep=crcr]{%
0.0625	0.0179248745804004\\
0.03125	0.00503225031064678\\
0.015625	0.00370734466383733\\
0.0078125	0.00115706927193495\\
0.00390625	0.000575085125397159\\
};
\addlegendentry{$\|\bar{u}(T)-\bar{u}_{kh}(T)\|$}

\addplot [color=black, dashdotted, line width=2pt]
  table[row sep=crcr]{%
0.0625	0.0625\\
0.03125	0.03125\\
0.015625	0.015625\\
0.0078125	0.0078125\\
0.00390625	0.00390625\\
};
\addlegendentry{O(h)}

\end{axis}
\end{tikzpicture}%
}\\
\caption{Spatial refinement and dg(0)}
\end{subfigure}
\begin{subfigure}[t]{.45\linewidth}
\centering
\scalebox{.5}
{
%
%
\begin{tikzpicture}

\begin{axis}[%
width=4.521in,
height=3.566in,
at={(0.758in,0.481in)},
scale only axis,
xmode=log,
xmin=0.001,
xmax=0.1,
xminorticks=true,
xlabel style={font=\color{white!15!black}},
xlabel={h},
ymode=log,
ymin=0.0001,
ymax=0.1,
yminorticks=true,
ylabel style={font=\color{white!15!black}},
ylabel={Error},
axis background/.style={fill=white},
legend style={at={(0.03,0.97)}, anchor=north west, legend cell align=left, align=left, draw=white!15!black}
]
\addplot [color=blue, mark=o,line width=2pt, mark options={solid, blue}]
  table[row sep=crcr]{%
0.0625	0.0182221567674502\\
0.03125	0.00495579396121221\\
0.015625	0.00363396151224937\\
0.0078125	0.000964053706530775\\
0.00390625	0.000419441712729643\\
};
\addlegendentry{$\|\bar{u}(T)-\bar{u}_{kh}(T)\|$}

\addplot [color=black,line width=2pt, dashdotted]
  table[row sep=crcr]{%
0.0625	0.0625\\
0.03125	0.03125\\
0.015625	0.015625\\
0.0078125	0.0078125\\
0.00390625	0.00390625\\
};
\addlegendentry{O(h)}

\end{axis}
\end{tikzpicture}%
}\\
\caption{Spatial refinement and dg(1)}
\end{subfigure}
\caption{Convergence results for spatial error}
\label{fig:spacerefinement}
\end{figure}

Next we practically verify the derived a priori error estimates for the optimal states. Let us first discuss spatial refinement. To this end we consider cG($1$)dG($r$) approximations for both $r=0$ and $r=1$ of the state equation on an equidistant grid in time with~$M=256$ steps and a sequence~$\{\mathcal{T}_i\}^6_{i=1}$ of spatial triangulations. Here,~$\mathcal{T}_{i+1}$ is obtained by one global uniform refinement of~$\mathcal{T}_i$,~$1 \leq i \leq 5$. The desired state~$u_d$ and the regularization parameter~$\alpha$ are chosen as in Section~\ref{subsec:inverseprob}. Since no analytic solution for this problem is known we take the optimal state on the finest spatial grid as a reference~$\ou$. On each refinement level, the optimal state~$\oukh$ is computed using the PDAP algorithm. The convergence plots are given in Figure~\ref{fig:spacerefinement}. For visual comparison we also plot the corresponding rate of convergence as given in Theorem~\ref{th:improved_h_error} without the logarithmic factor. We clearly see that the computed rates for the optimal states match the predicted order of~$\mathcal{O}(h)$ for both temporal approximation schemes. 
\subsection{Time refinement}
\begin{figure}[htb] 
\centering
\begin{subfigure}[t]{.45\linewidth}
\centering
\scalebox{.5}
{
%
%
\begin{tikzpicture}

\begin{axis}[%
width=4.521in,
height=3.566in,
at={(0.758in,0.481in)},
scale only axis,
xmode=log,
xmin=0.0001,
xmax=0.01,
xminorticks=true,
xlabel style={font=\color{white!15!black}},
xlabel={k},
ymode=log,
ymin=0.0001,
ymax=0.1,
yminorticks=true,
ylabel style={font=\color{white!15!black}},
ylabel={Error},
axis background/.style={fill=white},
legend style={at={(0.03,0.97)}, anchor=north west, legend cell align=left, align=left, draw=white!15!black}
]
\addplot [color=blue, mark=o, mark options={solid, blue}, line width=2pt]
  table[row sep=crcr]{%
0.00625	0.0322689764389044\\
0.003125	0.0176537112989011\\
0.0015625	0.00823675863102631\\
0.00078125	0.00289224928514931\\
};
\addlegendentry{$\|\bar{u}(T)-\bar{u}_{kh}(T)\|$}

\addplot [color=black, dashdotted, , line width=2pt]
  table[row sep=crcr]{%
0.00625	0.00625\\
0.003125	0.003125\\
0.0015625	0.0015625\\
0.00078125	0.00078125\\
};
\addlegendentry{$O(k)$}

\end{axis}
\end{tikzpicture}%
}\\
\caption{Temporal refinement for dg(0)}
\end{subfigure}
\begin{subfigure}[t]{.45\linewidth}
\centering
\scalebox{.5}
{
%
%
\begin{tikzpicture}

\begin{axis}[%
width=4.521in,
height=3.566in,
at={(0.758in,0.481in)},
scale only axis,
xmode=log,
xmin=0.0001,
xmax=0.01,
xminorticks=true,
xlabel style={font=\color{white!15!black}},
xlabel={k},
ymode=log,
ymin=1e-10,
ymax=0.001,
yminorticks=true,
ylabel style={font=\color{white!15!black}},
ylabel={Error},
axis background/.style={fill=white},
legend style={at={(0.03,0.97)}, anchor=north west, legend cell align=left, align=left, draw=white!15!black}
]
\addplot [color=blue, mark=o, mark options={solid, blue}, , line width=2pt]
  table[row sep=crcr]{%
0.00625	0.000147842508654195\\
0.003125	1.92434726969871e-05\\
0.0015625	2.42731068578816e-06\\
0.00078125	2.72808039410665e-07\\
};
\addlegendentry{$\|\bar{u}(T)-\bar{u}_{kh}(T)\|$}

\addplot [color=black, dashdotted, , line width=2pt]
  table[row sep=crcr]{%
0.00625	2.44140625e-07\\
0.003125	3.0517578125e-08\\
0.0015625	3.814697265625e-09\\
0.00078125	4.76837158203125e-10\\
};
\addlegendentry{$O(k^3)$}

\end{axis}
\end{tikzpicture}%
}\\
\caption{Temporal refinement for dg(1)}
\end{subfigure}
\caption{Convergence results for temporal error}
\label{fig:timerefinement}
\end{figure}
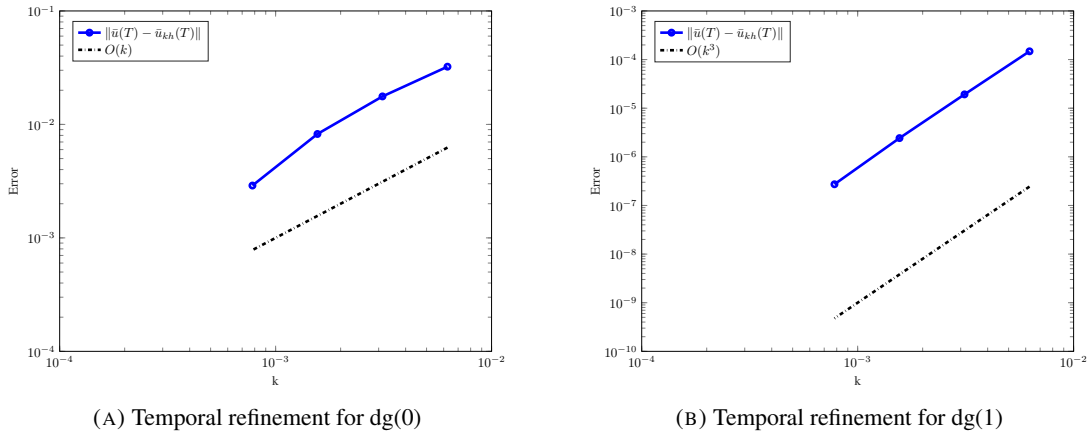

In order to verify the temporal error estimate we discretize the state equation again by the cG($1$)dG($r$) scheme for both $r=0$ and $r=1$, on equidistant time grids with~$2^i$ steps,~$i=4, \dots,8$, and a fixed triangulation of the spatial domain. The desired state~$u_d$ is chosen as the discrete final state corresponding to the measure~$q^\dagger$ on the finest discretization. The regularization parameter is set to~$\alpha=0.001$. Again, the optimal state on the finest discretization is considered as reference~$\ou$. The computed convergence results for the optimal states are plotted in Figure~\ref{fig:timerefinement} alongside the rates of convergence derived in Theorem~\ref{th:improved_k_error}. As predicted by the theory, we observe a linear~$\mathcal{O}(k)$ rate for~dG($0$) and a cubic~$\mathcal{O}(k^3)$ rate of convergence for~dG($1$).   

\bibliographystyle{siam}
\bibliography{smoothing}

\begin{thebibliography}{10}

\bibitem{AdamsFournier:1977}
{\sc R.~A. Adams and J.~Fournier}, {\em Cone conditions and properties of
  {S}obolev spaces}, J. Math. Anal. Appl., 61 (1977), pp.~713--734.

\bibitem{fista}
{\sc A.~Beck and M.~Teboulle}, {\em A fast iterative shrinkage-thresholding
  algorithm for linear inverse problems}, SIAM J. Imaging Sci., 2 (2009),
  pp.~183--202.

\bibitem{bogachev}
{\sc V.~I. Bogachev}, {\em Measure theory. {V}ol. {I}, {II}}, Springer-Verlag,
  Berlin, 2007.

\bibitem{Casas:1997}
{\sc E.~Casas}, {\em Pontryagin's principle for state-constrained boundary
  control problems of semilinear parabolic equations}, SIAM J. Control Optim.,
  35 (1997), pp.~1297--1327.

\bibitem{CasasClasonKunisch:2012}
{\sc E.~Casas, C.~Clason, and K.~Kunisch}, {\em Approximation of elliptic
  control problems in measure spaces with sparse solutions}, SIAM J. Control
  Optim., 50 (2012), pp.~1735--1752.

\bibitem{CasasE_ClasonC_KunischK_2012}
\leavevmode\vrule height 2pt depth -1.6pt width 23pt, {\em Approximation of
  elliptic control problems in measure spaces with sparse solutions}, SIAM J.
  Control Optim., 50 (2012), pp.~1735--1752.

\bibitem{CasasE_VexlerB_ZuazuaE_2015}
{\sc E.~Casas, B.~Vexler, and E.~Zuazua}, {\em Sparse initial data
  identification for parabolic {PDE} and its finite element approximations},
  Math. Control Relat. Fields, 5 (2015), pp.~377--399.

\bibitem{CasasZuazua:2013}
{\sc E.~Casas and E.~Zuazua}, {\em Spike controls for elliptic and parabolic
  {PDE}s}, Systems Control Lett., 62 (2013), pp.~311--318.

\bibitem{ClasonKunisch:2011}
{\sc C.~Clason and K.~Kunisch}, {\em A duality-based approach to elliptic
  control problems in non-reflexive {B}anach spaces}, ESAIM Control Optim.
  Calc. Var., 17 (2011), pp.~243--266.

\bibitem{Duval2015}
{\sc V.~Duval and G.~Peyr\'e}, {\em Exact support recovery for sparse spikes
  deconvolution}, Found. Comput. Math., 15 (2015), pp.~1315--1355.

\bibitem{ErikssonK_JohnsonC_LarssonS_1998a}
{\sc K.~Eriksson, C.~Johnson, and S.~Larsson}, {\em Adaptive finite element
  methods for parabolic problems. {VI}. {A}nalytic semigroups}, SIAM J. Numer.
  Anal., 35 (1998), pp.~1315--1325 (electronic).

\bibitem{ErikssonK_JohnsonC_ThomeeV_1985}
{\sc K.~Eriksson, C.~Johnson, and V.~Thom\'ee}, {\em Time discretization of
  parabolic problems by the discontinuous {G}alerkin method}, RAIRO Mod\'el.
  Math. Anal. Num\'er., 19 (1985), pp.~611--643.

\bibitem{LCEvans_2010}
{\sc L.~C. Evans}, {\em Partial differential equations}, vol.~19 of Graduate
  Studies in Mathematics, American Mathematical Society, Providence, RI,
  second~ed., 2010.

\bibitem{FabrePUelZuazua:1995}
{\sc C.~Fabre, J.-P. Puel, and E.~Zuazua}, {\em On the density of the range of
  the semigroup for semilinear heat equations}, in Control and optimal design
  of distributed parameter systems ({M}inneapolis, {MN}, 1992), vol.~70 of IMA
  Vol. Math. Appl., Springer, New York, 1995, pp.~73--91.

\bibitem{Griepentrog:2007}
{\sc J.~A. Griepentrog}, {\em Maximal regularity for nonsmooth parabolic
  problems in {S}obolev-{M}orrey spaces}, Adv. Differential Equations, 12
  (2007), pp.~1031--1078.

\bibitem{HansboA_2002a}
{\sc A.~Hansbo}, {\em Strong stability and non-smooth data error estimates for
  discretizations of linear parabolic problems}, BIT, 42 (2002), pp.~351--379.

\bibitem{IsakovV_2017}
{\sc V.~Isakov}, {\em Inverse problems for partial differential equations},
  vol.~127 of Applied Mathematical Sciences, Springer, Cham, third~ed., 2017.

\bibitem{kantorovich}
{\sc L.~V. Kantorovic and G.~v.~S. Rubinstein}, {\em On a space of completely
  additive functions}, Vestnik Leningrad. Univ., 13 (1958), pp.~52--59.

\bibitem{Kovats:2002}
{\sc J.~Kovats}, {\em Real analytic solutions of parabolic equations with
  time-measurable coefficients}, Proc. Amer. Math. Soc., 130 (2002),
  pp.~1055--1064.

\bibitem{KunischPieperVexler:2014}
{\sc K.~Kunisch, K.~Pieper, and B.~Vexler}, {\em Measure valued directional
  sparsity for parabolic optimal control problems}, SIAM J. Control Optim., 52
  (2014), pp.~3078--3108.

\bibitem{Leis:1997}
{\sc R.~Leis}, {\em Initial-boundary value problems in mathematical physics},
  in Modern mathematical methods in diffraction theory and its applications in
  engineering ({F}reudenstadt, 1996), vol.~42 of Methoden Verfahren Math.
  Phys., Peter Lang, Frankfurt am Main, 1997, pp.~125--144.

\bibitem{LeykekhmanD_VexlerB_2016a}
{\sc D.~Leykekhman and B.~Vexler}, {\em Pointwise best approximation results
  for {G}alerkin finite element solutions of parabolic problems}, SIAM J.
  Numer. Anal., 54 (2016), pp.~1365--1384.

\bibitem{LeykekhmanD_VexlerB_2017a}
\leavevmode\vrule height 2pt depth -1.6pt width 23pt, {\em Discrete maximal
  parabolic regularity for {G}alerkin finite element methods}, Numer. Math.,
  135 (2017), pp.~923--952.

\bibitem{DMeidner_RRannacher_BVexler_2011a}
{\sc D.~Meidner, R.~Rannacher, and B.~Vexler}, {\em A priori error estimates
  for finite element discretizations of parabolic optimization problems with
  pointwise state constraints in time}, SIAM J. Control Optim., 49 (2011),
  pp.~1961--1997.

\bibitem{MerinoNeitzelTroeltzsch:2011}
{\sc P.~Merino, I.~Neitzel, and F.~Tr\"oltzsch}, {\em On linear-quadratic
  elliptic control problems of semi-infinite type}, Appl. Anal., 90 (2011),
  pp.~1047--1074.

\bibitem{milzarekfilter}
{\sc A.~Milzarek and M.~Ulbrich}, {\em A semismooth {N}ewton method with
  multidimensional filter globalization for {$l_1$}-optimization}, SIAM J.
  Optim., 24 (2014), pp.~298--333.

\bibitem{Palencia:1994}
{\sc C.~Palencia}, {\em On the stability of variable stepsize rational
  approximations of holomorphic semigroups}, Math. Comp., 62 (1994),
  pp.~93--103.

\bibitem{PieperVexler:2013}
{\sc K.~Pieper and B.~Vexler}, {\em A priori error analysis for discretization
  of sparse elliptic optimal control problems in measure space}, SIAM J.
  Control Optim., 51 (2013), pp.~2788--2808.

\bibitem{RannacherVexler:2005}
{\sc R.~Rannacher and B.~Vexler}, {\em A priori error estimates for the finite
  element discretization of elliptic parameter identification problems with
  pointwise measurements}, SIAM J. Control Optim., 44 (2005), pp.~1844--1863.

\bibitem{AHSchatz_LBWahlbin_1977a}
{\sc A.~H. Schatz and L.~B. Wahlbin}, {\em Interior maximum norm estimates for
  finite element methods}, Math. Comp., 31 (1977), pp.~414--442.

\bibitem{AHSchatz_LBWahlbin_1995a}
\leavevmode\vrule height 2pt depth -1.6pt width 23pt, {\em Interior
  maximum-norm estimates for finite element methods. {II}}, Math. Comp., 64
  (1995), pp.~907--928.

\bibitem{shapiroextremalvalue}
{\sc A.~Shapiro}, {\em Second-order derivatives of extremal-value functions and
  optimality conditions for semi-infinite programs}, Math. Oper. Res., 10
  (1985), pp.~207--219.

\bibitem{ThomeeV_XuJ_NaiY_1989}
{\sc V.~Thom\'ee, J.~Xu, and N.~Y. Zhang}, {\em Superconvergence of the
  gradient in piecewise linear finite-element approximation to a parabolic
  problem}, SIAM J. Numer. Anal., 26 (1989), pp.~553--573.

\bibitem{Walter:2019}
{\sc D.~Walter}, {\em On sparse sensor placement for parameter identification
  problems with partial dierential equations}, PhD thesis, Technische
  Universit{\"a}t M{\"u}nchen, 2019.

\end{thebibliography}
\end{document}